\documentclass[11pt]{article}%
\usepackage{graphicx,subfigure}
\usepackage{amsfonts}
\usepackage{mathrsfs}
\usepackage{amssymb}
\usepackage{chngpage}
\usepackage{amsmath}
\usepackage{graphicx}
\usepackage{amsmath,amssymb}
\usepackage[mathscr]{eucal}
\usepackage{xcolor}%
\providecommand{\U}[1]{\protect\rule{.1in}{.1in}}
\newtheorem{theorem}{Theorem}

\newtheorem{corollary}[theorem]{Corollary}
\newtheorem{lemma}[theorem]{Lemma}
\newtheorem{remark}[theorem]{Remark}

\numberwithin{equation}{section}
\numberwithin{theorem}{section}

\def\({\left(}
\def\){\right)}

\begin{document}
\centerline{\Large\bf Dynamics and Wong-Zakai approximations of stochastic  }
\vspace{0.3cm} \centerline {\Large\bf  nonlocal PDEs with long time memory }
%==========================================

\vskip 16pt\centerline{\bf Jiaohui Xu$^1$, {\bf Tom$\acute{\bf a}$s   Caraballo$^{2,3,}$\footnote{Corresponding author}\footnote{E-mail addresses: jxu@us.es (J. Xu),
caraball@us.es (T. Caraballo), jvalero@umh.es (J. Valero).}, \bf Jos\'e Valero$^4$}}

\centerline {\footnotesize \it $^1$Center for Nonlinear Studies, School of
Mathematics,} \centerline{\footnotesize\it Northwest University, Xi'an
710127, P. R. China}

%%\centerline{{\bf Tom$\acute{\bf a}$s   Caraballo$^{2,3}$}\footnote{Corresponding author}
%%\footnote{ E-mail addresses: jiaxu1@alum.us.es (J. Xu),
%%caraball@us.es (T. Caraballo), jvalero@umh.es (J. Valero).}}
\centerline{\footnotesize\it $^2$Dpto. Ecuaciones Diferenciales y
An$\acute{a}$lisis Num$\acute{e}$rico, } \centerline{\footnotesize\it
Facultad de Matem\'aticas, Universidad de Sevilla,  c/ Tarfia s/n,
41012-Sevilla, Spain} \centerline{\footnotesize\it $^3$Department of
Mathematics, Wenzhou University, Wenzhou, Zhejiang Province, 325035, P. R.
China}

%%\medskip
%%
%%\centerline{\bf Jos\'e Valero$^4$}
%%
\centerline{\footnotesize\it  $^4$Centro de Investigaci\'on Operativa,
Universidad Miguel Hern\'andez de Elche,} \centerline{\footnotesize\it
Avenida de la Universidad s/n, 03202-Elche, Spain}

\medskip

\bigbreak\noindent\textbf{Abstract} In this paper, a combination of Galerkin's
method and Dafermos' transformation is first used to prove the existence and
uniqueness of solutions for a class of stochastic nonlocal PDEs with long time
memory driven by additive noise. Next, the existence of tempered random
attractors for such equations is established in an appropriate space for the
analysis of problems with delay and memory. Eventually, the convergence of
solutions of Wong-Zakai approximations and upper semicontinuity of random
attractors of the approximate random system, as the step sizes of
approximations approach zero, are analyzed in a detailed way. \vskip5pt

\medskip

\noindent\textit{Keywords:} {Long time memory, Wong-Zakai approximation,
Dafermos transformation, Random attractors, Upper semicontinuity.}

%\noindent{\it AMS subject classifications}. 35R11, 35Q30, 65F08, 60H15, 65F10

\section{Introduction}

Motivated by some interesting physical problems related to thermal memory or
materials with memory, several papers have been published (see \cite{C1, G3,
G1, G2, X0, X2} and the references therein) concerning a semilinear partial
differential equation to model the heat flow in a rigid, isotropic,
homogeneous heat conductor with linear memory. The equation is the following,
\begin{equation}
\hskip-1cm%
\begin{cases}
\displaystyle c_{0}\partial_{t}u-k_{0}\Delta
u-{\textcolor{red}\displaystyle\int_{-\infty}^{t}k(t-s)\Delta u(s)ds}%
+f(u)=h,\\
u(x,t)=0,\\[0.8ex]%
{u(x,\tau+t)=u_{0}(x,t)},\\
\end{cases}
\begin{aligned} &\hskip-.3 cm\mbox{in}~~\mathcal{O}\times (\tau,+\infty),\\ &\hskip-.3 cm\mbox{on}~\partial\mathcal{O}\times(\tau,+\infty),\\[0.8ex] &\hskip-.3 cm\mbox{in}~~\mathcal{O}\times(-\infty,0],\\ \end{aligned} \label{Eq1-1}%
\end{equation}
where $\tau\in\mathbb{R}$, $\mathcal{O}\subset\mathbb{R}^{N}$ is a bounded
domain with regular boundary, {{$u:\mathcal{O}\times\mathbb{R}\rightarrow
\mathbb{R}$ is the temperature field}}, {{$k:\mathbb{R}^{+}\rightarrow
\mathbb{R}$ is the heat flux memory kernel}}, $\mathbb{R}^{+}$ denotes the
interval $(0,+\infty)$, $c_{0}$ and $k_{0}$ denote the specific heat and the
instantaneous conductivity, respectively.

To solve \eqref{Eq1-1} successfully, one can make the past history of $u$,
from $-\infty$ to $0$, be part of the forcing term given by the causal
function $g$ as follows,
\[
g(x,t)={h(x,t)}+{\int_{-\infty}^{\tau}k(t-s)\Delta u(x,s)ds},\qquad
x\in\mathcal{O},~~t\geq\tau.
\]
In this way, \eqref{Eq1-1} becomes an initial value problem without delay or
memory,
\begin{equation}
\hskip-1.5cm%
\begin{cases}
c_{0}\partial_{t}u-k_{0}\Delta u-\displaystyle\int_{\tau}^{t}k(t-s)\Delta
u(s)ds+f(u)=g,\\
u(x,t)=0,\\
{u(x,\tau)=u_{0}(x)},\\
\end{cases}
\begin{aligned} &\mbox{in}~\mathcal{O}\times (\tau,+\infty),\\ &\mbox{on}~\partial\mathcal{O}\times (\tau,+\infty),\\ &\mbox{in}~\mathcal{O}.\\ \end{aligned} \label{eq1-1bisd}%
\end{equation}
However, proceeding in this way, we cannot construct a dynamical system
generated by the solutions of the original problem \eqref{Eq1-1} in a correct
way, since the history part of the function $u$ is necessary to solve problem \eqref{eq1-1bisd}.

Therefore, two alternatives have been carried out to handle the problem in a
correct mathematical way.

$\bullet$ {Alternative 1:} Based on {Dafermos'} idea for linear
viscoelasticity problems (see, e.g., \cite{C1}), one can define the new
variables,
\[
u^{t}(x,s)=u(x,t-s),\qquad s\geq0,
\]%
\begin{equation}
\eta^{t}(x,s)=\int_{0}^{s}u^{t}(x,r)dr=\int_{t-s}^{t}u(x,r)dr,\qquad s\geq0.
\label{eq3.2a}%
\end{equation}
Assuming $k(\infty)=0$, thanks to a change of variable and a formal
integration by parts, we obtain
\[
\int_{-\infty}^{t}k(t-s)\Delta{u(s)}ds=-\int_{0}^{\infty}k^{\prime}%
(s)\Delta{\eta^{t}(s)}ds.
\]
Here and in the sequel, the \textit{prime} denotes derivation with respect to
the variable $s$. Setting
\[
\mu(s)=-k^{\prime}(s),
\]
the original problem \eqref{eq1-1bisd} becomes an autonomous one without
delay,
\begin{equation}
\hskip-1.5cm%
\begin{cases}
c_{0}\partial_{t}u-k_{0}\Delta u-\displaystyle\int_{0}^{\infty}\mu
(s)\Delta\eta^{t}(s)ds+f(u)=h,\\
{\eta_{t}^{t}}(s)=-{\eta_{s}^{t}}(s)+u(t),\\
u(x,t)=\eta^{t}(s,x)=0,\\
u(x,\tau)={u_{0}(x)},\\
\eta^{\tau}(x,s)={\eta_{0}(s)},
\end{cases}
\begin{aligned} &\mbox{in}~~\mathcal{O}\times (\tau,\infty),\\ &\mbox{in}~~\mathcal{O}\times (\tau,\infty), s>0,\\ &\mbox{on}~~\partial \mathcal{O}\times(\tau,\infty), s>0,\\ &\mbox{in}~~ \mathcal{O},\\ &\mbox{in}~~\mathcal{O}\times\mathbb{R}^+. \end{aligned} \label{eq3.3a}%
\end{equation}
From the definition of $\eta^{t}(x,s)$ (see (\ref{eq3.2a})), we see that
\begin{equation}
\eta_{0}(s)=\int_{\tau-s}^{\tau}u(r)dr{=\int_{\tau-s}^{\tau}u_{0}%
(r-\tau)dr=\int_{-s}^{0}u_{0}(r)dr}, \label{eq2.5a}%
\end{equation}
which is the initial integrated past history of $u$ with vanishing boundary.
Consequently, any solution to \eqref{eq1-1bisd} is a solution to
\eqref{eq3.3a} for the corresponding initial values $(u_{0},\eta_{0})$ given
by \eqref{eq2.5a}.

Observe that problem \eqref{eq3.3a} can be solved for arbitrary initial values
$(u_{0},\eta_{0})$ in a proper phase space $L^{2}(\mathcal{O})\times L_{\mu
}^{2}(\mathbb{R}^{+};H_{0}^{1}(\mathcal{O}))$, i.e., the second component
$\eta_{0}$ does not necessarily depend on $u_{0}(\cdot)$, where {$L_{\mu}%
^{2}(\mathbb{R}^{+};H_{0}^{1}(\mathcal{O}))$} is {a Hilbert space specified
later.} \medskip

Let $\mu$ satisfy the hypotheses:

\begin{enumerate}
\item[$(h_{1})$] $\mu\in C^{1}(\mathbb{R}^{+})\cap L^{1}(\mathbb{R}^{+}%
),\quad\mu(s)\geq0,\quad\mu^{\prime}(s)\leq0$,\quad$\forall s\in\mathbb{R}%
^{+}$;

\item[$(h_{2})$] $\mu^{\prime}(s)+\varpi\mu(s)\leq0$,\quad$\forall
s\in\mathbb{R}^{+}$, ~~for some~ $\varpi>0$.
\end{enumerate}

Then $L^{2}_{\mu}(\mathbb{R}^{+};H^{1}_{0}(\mathcal{O}))$ is a {Hilbert space}
of functions $w:\mathbb{R}^{+}\rightarrow H^{1}_{0}(\mathcal{O})$ with inner
product,
\[
((w_{1},w_{2}))_{\mu}=\int_{0}^{\infty}\mu(s)(\nabla w_{1}(s),\nabla
w_{2}(s))ds.
\]
The solutions of \eqref{eq3.3a} are proved to exist in \cite{C1} and permit to
construct a dynamical system $S(t):L^{2}(\mathcal{O})\times L^{2}_{\mu
}(\mathbb{R}^{+};H^{1}_{0}(\mathcal{O}))\to L^{2}(\mathcal{O})\times
L^{2}_{\mu}(\mathbb{R}^{+};H^{1}_{0}(\mathcal{O}))$ via,
\[
S(t)(u_{0},\eta_{0})=(u(t;0,{(u_{0},\eta_{0})}),\eta^{t}(\cdot;0,{(u_{0}%
,\eta_{0})})),
\]
which possesses a global attractor in this phase space. However, this global
attractor does not reflect the complete asymptotic dynamics of the original
problem \eqref{Eq1-1} since the latter problem is not equivalent to
\eqref{Eq1-1}. In other words, not for every $\eta_{0}\in L^{2}_{\mu
}(\mathbb{R}^{+};H^{1}_{0}(\mathcal{O}))$, there exists $u_{0}:(-\infty,0]\to
H^{1}_{0}(\mathcal{O})$ such that,
\[
\eta_{0}(s)=\int_{-s}^{0}u_{0}(r)dr, \qquad\forall s\geq0.
\]
In fact, both problems are equivalent (cf. \cite{G2}) if and only if the
initial value $\eta_{0}$ belongs to a proper subspace of $L^{2}_{\mu
}(\mathbb{R}^{+};H^{1}_{0}(\mathcal{O}))$. Precisely, the domain of the
distributional derivative with respect to $s$, denoted by $D(\mathbf{T})$,
\[
D(\mathbf{T})=\left\{  \eta(\cdot)\in L^{2}_{\mu}(\mathbb{R}^{+};H^{1}%
_{0}(\mathcal{O}))\ | \ \eta_{s}(\cdot)\in L^{2}_{\mu}(\mathbb{R}^{+};
H^{1}_{0}(\mathcal{O})), ~\eta(0)=0\right\}  ,
\]
and $\mathbf{T}$ is defined by $\mathbf{T}\eta=-\eta_{s}$, $\eta\in
D(\mathbf{T}).$

Hence, it seems natural to construct a dynamical system generated by
\eqref{eq3.3a} in $L^{2}(\mathcal{O})\times D(\mathbf{T})$ and to prove the
existence of attractors to problem \eqref{Eq1-1} via the above relationship.
Up to our knowledge, it is not possible to prove the existence of attractors
in $L^{2}(\mathcal{O})\times D(\mathbf{T})$ unless the solutions admit more
regularity.
%These facts motivate the second alternative.

$\bullet$ Alternative 2: The idea comes from a simpler case in \cite{C10} when
the kernel is the so called non-singular one and has the expression
{$k(t)=e^{-d_{0}t}$, $d_{0}>0$, considering the phase space $L_{H_{0}^{1}}%
^{2}$ formed by the functions $\varphi:(-\infty,0]\rightarrow H_{0}%
^{1}(\mathcal{O})$ with $\int_{-\infty}^{0}{e^{\gamma s}}\Vert\varphi
(s)\Vert_{H_{0}^{1}}^{2}ds<+\infty$ for {certain $\gamma>0$}. The authors in
\cite{C10} proved that solutions of problem \eqref{Eq1-1} with initial value
$u_{0}$ generate a dynamical system which possesses a global attractor in
$L_{H_{0}^{1}}^{2}$. However, when working with delay problems, it is natural
(see e.g., \cite{delfour, CR} and the references therein) to consider the
phase space $L^{2}(\mathcal{O})\times L_{H_{0}^{1}}^{2}$ and set up the
problem as,
\begin{equation}
\hskip-1cm%
\begin{cases}
\displaystyle c_{0}\partial_{t}u-k_{0}\Delta
u-{\textcolor{red}\displaystyle\int_{-\infty}^{t}k(t-s)\Delta uds}+f(u)=h,\\
u(x,t)=0,\\[0.8ex]%
u(x,\tau)=u_{0}(x),\\
u(x,\tau+{t})=\varphi(x,{t}),\\
\end{cases}
\begin{aligned} &\hskip-.3 cm\mbox{in}~~\mathcal{O}\times (\tau,+\infty),\\ &\hskip-.3 cm\mbox{on}~\partial\mathcal{O}\times(\tau,+\infty),\\[0.8ex] &\hskip-.3 cm\mbox{in}~\mathcal{O},\\ &\hskip-.3 cm\mbox{in}~~\mathcal{O}\times{{(-\infty,0)}}.\\ \end{aligned} \label{eq1-1a}%
\end{equation}
Thanks to the results in \cite{C10}, we are able to construct a dynamical
system $S(t):L^{2}(\mathcal{O})\times L_{H_{0}^{1}}^{2}\rightarrow
L^{2}(\mathcal{O})\times L_{H_{0}^{1}}^{2}$ via the relation,
\begin{equation}
{S(t)(u_{0},\varphi):=(u(t;0,(u_{0},\varphi)),u_{t}(\cdot;0,(u_{0}%
,\varphi))),} \label{semigroup}%
\end{equation}
where $u(\cdot;0,(u_{0},\varphi))$ denotes the solution of problem
\eqref{eq1-1a} (see \cite{CR}), and $u_{t}$ the history up to time $t$,
\[
u_{t}(s;0,(u_{0},\varphi))=u(t+s;0,(u_{0},\varphi)),\quad s\leq0.
\]
The existence of global attractors in the space $L^{2}(\mathcal{O})\times
L_{H_{0}^{1}}^{2}$ is proved in \cite{C10}. In fact, it was proved for a
non-autonomous version which is much more general than the one explained here.
Nevertheless, as we mentioned before, the technique applied in this case
(essentially the Galerkin approach) requires the kernel to be non-singular
($k(t)=e^{-d_{0}t}$, $d_{0}>0$). This is a strong restriction on the kernel
{$k$} (and consequently, on $\mu$}) because in applied science singularities
appear very often, e.g., ${k(t)=e^{-d_{0}t}t^{-\alpha}},\alpha\in(0,1)$.
Motivated by this fact, recently we have proved in \cite{X2} the existence of
global attractors in $L^{2}(\mathcal{O})\times L_{H_{0}^{1}}^{2}$ for the
general singular case, even for a more general model containing nonlocal
diffusion coefficients thanks to a combination of Galerkin's method and
Dafermos' transformation. More precisely, the following nonlocal PDE
associated with singular memory {was} considered in \cite{X2},
\begin{equation}
{\small \hskip-1.3cm%
\begin{cases}
\partial_{t}u-{a(l(u))}\Delta u-\displaystyle{\int_{-\infty}^{t}}%
{k(t-s)}\Delta uds+f(u)=g,\\
u(x,t)=0,\\
u(x,\tau)={u_{0}(x)},\\
u(x,t+\tau)={\phi(x,t)},
\end{cases}
\begin{aligned} &\mbox{in }\mathcal{O}\times (\tau,\infty),\\ &\mbox{on
}\partial\mathcal{O}\times(\tau,\infty),\\ &\mbox{in }\mathcal{O}\\ &\mbox{in
}\mathcal{O}\times (-\infty,0),\\ \end{aligned}} \label{eq0}%
\end{equation}
where the function $a\in C(\mathbb{R};\mathbb{R}^{+})$ satisfies,
\begin{equation}
0<m\leq a(r),\qquad\forall r\in\mathbb{R}, \label{eq2.2}%
\end{equation}
with initial value $u_{0}\in L^{2}(\mathcal{O})$ and initial function $\phi\in
L_{H_{0}^{1}(\mathcal{O})}^{2}$. Then, the semigroup defined as in
\eqref{semigroup} for the solutions of \eqref{eq0} possesses a global
attractor in the phase space $L^{2}(\mathcal{O})\times L_{H_{0}^{1}}^{2}$.
Also it is worth mentioning that the same results were proved in the phase
space $L^{2}(\mathcal{O})\times L_{\mu}^{2}(\mathbb{R}^{+};{H_{0}%
^{1}(\mathcal{O})})$ (see \cite{X3}) by using the classical Dafermos' method.

Our interest in this paper is to analyze the behavior of the nonlocal problem
with memory when some stochastic disturbance appears in the model. Assume that
this perturbation appears as an additive noise, namely, our objective is to
study the following stochastic nonlocal PDEs with long time memory,
\begin{equation}%
\begin{cases}
\partial_{t}u-a(l(u))\Delta u-\displaystyle\int_{-\infty}^{t}k(t-s)\Delta
uds+f(u)=h+\phi\frac{dW(t)}{dt},\\
u(x,t)=0,\\
u(x,\tau)=u_{0}(x),\\
u(x,t+\tau)=\varphi(x,t),\\
\end{cases}
\begin{aligned} &\mbox{in}~~\mathcal{O}\times (\tau,\infty),\\ &\mbox{on}~\partial\mathcal{O}\times(\tau,\infty),\\ &\mbox{in}~~\mathcal{O},\\ &\mbox{in}~~\mathcal{O}\times (-\infty,0),\\ \end{aligned} \label{eq1-1}%
\end{equation}
where $\tau\in\mathbb{R}$, $l\in\mathcal{L}(L^{2}(\mathcal{O});\mathbb{R})$,
$\mathcal{O}\subset\mathbb{R}^{N}$ is a fixed bounded domain with regular
boundary, $h\in L^{2}(\mathcal{O})$, $W(t)$ is a two-sided standard Brownian
motion on a probability space $(\Omega,\mathcal{F},\mathbb{P})$ and $f$ is a
polynomial of odd degree $2p-1,\ p\in\mathbb{N}$. Suppose that there are two
constants $m$ and $M$ such that the function $a\in C(\mathbb{R};\mathbb{R}%
^{+})$ satisfies,
\begin{equation}
0<m\leq a(r)\leq M,\qquad\forall r\in\mathbb{R}. \label{eq1-2}%
\end{equation}
Here, $k:\mathbb{R}^{+}\rightarrow\mathbb{R}$ is the memory kernel whose
properties will be specified later. The initial value $u_{0}$ belongs to
$L^{2}(\mathcal{O})$, while the initial function $\varphi$ belongs to the
space $L_{H_{0}^{1}}^{2},$ which is given by the measurable functions
$\varphi:(-\infty,0)\rightarrow{H_{0}^{1}(\mathcal{O})}$, such that
\[
\int_{-\infty}^{0}e^{\gamma s}\Vert\varphi(s)\Vert_{H_{0}^{1}}^{2}ds<\infty,
\]
for certain $\gamma>0$. Furthermore, assume that $\phi\in H_{0}^{1}%
(\mathcal{O})\cap H^{2}(\mathcal{O})\cap L^{2p}(\mathcal{O})$ is such that
$\Delta\phi\in L^{2p}(\mathcal{O})$.
%by the classical Sobolev embedding result, it is obvious that $\phi\in L^r(\mathcal{O})$ for $r\leq\frac{2N}{N-2}$.
A local version of this problem has been analyzed in \cite{C2} when the noise
is a Hilbert valued Wiener process by using Dafermos' transformation,
obtaining the existence of random attractors in the corresponding phase space
for the Dafermos set-up. In this case, our intention now is to construct a
framework to solve problems with delay and memory in an appropriate phase
space, like in the deterministic model (see \cite{X2}).

{In general, the Wiener process $W$ can be chosen as a stochastic process to
represent the position of the Brownian particle, but the velocity of the
particle cannot be obtained from the Wiener process because of the nowhere
differentiability of the sample paths of $W$ \cite{G6}. Therefore, it is
natural to approximate Brownian motion by more regular stochastic process,
which is the so-called colored noise. In recent decades, the Wong-Zakai
approximations to reaction-diffusion differential equations have been
extensively studied in the literature, see, e.g., \cite{B1, G5, L3, N1, T3,
T2, W4, WZ, Y1} and the references therein.} One of the goals of this paper is
to derive the relations between the solutions of problem \eqref{eq1-1} and the
corresponding limiting problem.

To this end, in Section \ref{s2}, we will first set-up problem \eqref{eq1-1}
in an appropriate form and will do the transformation to obtain a random
partial differential equation with delay. Then, a random partial differential
system is obtained thanks to Dafermos' transformation. We will also include in
Section \ref{s2} some necessary preliminaries and notation to tackle our
problem. The well-posedness of the transformed system is proved in Section
\ref{s3}. Next, we prove in Section \ref{s4} the existence of random
attractors to problem \eqref{eq1-1} in the phase space $L^{2}(\mathcal{O}%
)\times L_{H_{0}^{1}}^{2}$. In Section \ref{s5}, we consider the approximation
of the original problem by a parameterized family of problems containing
colored noise which possesses a parameterized family of corresponding random
attractors. Finally, in Section \ref{s6}, we prove the upper-semicontinuity
property of this family of parameterized random attractors with respect to the
random attractors to problem \eqref{eq1-1}.

\section{ Preliminaries}

\label{s2}

\subsection{Set-up of the problem}

The standard probability space $(\Omega,\mathcal{F},\mathbb{P})$ will be used
throughout this paper, where $\Omega=\{\omega\in C(\mathbb{R};\mathbb{R}%
):\omega(0)=0\}$, $\mathcal{F}$ is the Borel $\sigma$-algebra induced by the
compact-open topology of $\Omega$ and $\mathbb{P}$ is the Wiener measure on
$(\Omega,\mathcal{F})$. Given $t\in\mathbb{R}$, define $\theta_{t}%
:\Omega\rightarrow\Omega$ by
\[
\theta_{t}\omega(\cdot)=\omega(\cdot+t)-\omega(t),\quad\omega\in\Omega.
\]
Then $(\Omega,\mathcal{F},\mathbb{P},\{\theta_{t}\}_{t\in\mathbb{R}})$ is a
metric dynamical system and we identify $W(t,\omega)=\omega(t)$. Let
$z(t,\omega)$ be the unique stationary solution of the stochastic equation
$dz=-zdt+dW(t)$. This stationary solution is given by $z(t,\omega)=z_{\ast
}(\theta_{t}\omega)$, where the random variable $z_{\ast}(\omega)$ is defined
as,
\begin{equation}
z_{\ast}(\omega)=-\int_{-\infty}^{0}e^{s}W(s,\omega)ds. \label{eqz}%
\end{equation}
In addition, it follows from \cite{X1} that there exists a $\theta_{t}%
$-invariant set of full measure such that $z_{\ast}(\theta_{t}\omega)$ is
pathwise continuous for each fixed $\omega\in\Omega$ and satisfies,
\begin{equation}
\lim_{t\rightarrow\pm\infty}\frac{|z_{\ast}(\theta_{t}\omega)|}{|t|}%
=0\quad\mbox{and}\quad\lim_{t\rightarrow\pm\infty}\frac{1}{t}\int_{0}%
^{t}z_{\ast}(\theta_{s}\omega)ds=0. \label{zproperty}%
\end{equation}

We now transform the stochastic equation \eqref{eq1-1} into a pathwise
deterministic one by using the random variable $z_{\ast}$. Given $\tau
\in\mathbb{R}$, $t\geq\tau$ and $\omega\in\Omega$, if $u(t,\omega)$ is a
solution of \eqref{eq1-1}, then we introduce a new variable $v(t,\omega)$, by
\begin{equation}
v(t,\omega)=u(t,\omega)-\phi z_{\ast}(\theta_{t}\omega). \label{eq1-3}%
\end{equation}
In this way, problem \eqref{eq1-1} can be rewritten as,
\begin{equation}%
\begin{cases}
v_{t}-a(l(v+\phi z_{\ast}(\theta_{t}\omega)))\Delta v-a(l(v+\phi z_{\ast
}(\theta_{t}\omega)))z_{\ast}(\theta_{t}\omega)\Delta\phi\\
\quad\displaystyle-\int_{-\infty}^{t}k(t-s)\Delta vds+f(v+\phi z_{\ast}%
(\theta_{t}\omega))=h+\phi z_{\ast}(\theta_{t}\omega)+{z_{k}^{\phi}(\theta
_{t}\omega)},\\
v(x,t)=0,\\
v(x,\tau)=v_{0}(x):=u_{0}(x)-\phi z_{\ast}(\theta_{\tau}\omega),\\
v(x,t+\tau):=u(x,t+\tau)-\phi z_{\ast}(\theta_{t+\tau}\omega)=\varphi
(x,t)-\phi z_{\ast}(\theta_{t+\tau}\omega):=\varphi_{v}(x,t),\\
\end{cases}
\begin{aligned} &\quad\\[2.2ex] &\mbox{in}~~\mathcal{O}\times (\tau,\infty),\\ &\mbox{on}~\partial\mathcal{O}\times(\tau,\infty),\\ &\mbox{in}~~\mathcal{O},\\ &\mbox{in}~~\mathcal{O}\times (-\infty,0),\\ \end{aligned} \label{eq1-4}%
\end{equation}
where $z_{k}^{\phi}(\omega)$ is a process defined by
\begin{equation}
{z_{k}^{\phi}(\omega)}=\Delta\phi\int_{-\infty}^{0}k(-s)z_{\ast}(\theta
_{s}\omega)ds. \label{eqphi}%
\end{equation}
Notice that a change of variable yields that%
\[
{z_{k}^{\phi}(\theta_{t}\omega)}=\Delta\phi\int_{-\infty}^{0}k(-s)z_{\ast
}(\theta_{s}\theta_{t}\omega)ds=\Delta\phi\int_{-\infty}^{t}k(t-s)z_{\ast
}(\theta_{s}\omega)ds.
\]

In order to use Dafermos' transform (see \cite{C1}) to establish the
well-posedness of problem \eqref{eq1-4}, let us define the new variables,
\[
v^{t}(x,s,\omega)=v(x,t-s,\omega),\qquad s\geq0,
\]%
\[
\eta^{t}(x,s,\omega)=\int_{0}^{s}v^{t}(x,r,\omega)dr=\int_{t-s}^{t}%
v(x,r,\omega)dr,\qquad s\geq0.
\]
Besides, assuming $k(\infty)=0$, a change of variable and a formal integration
by parts (see Lemma \ref{lemEquivInt} for a rigorous calculation) imply,
\[
\int_{-\infty}^{t}k(t-s)\Delta v(s)ds=-\int_{0}^{\infty}k^{\prime}%
(s)\Delta\eta^{t}(s)ds.
\]
Setting $\mu(s)=-k^{\prime}(s)$, problem \eqref{eq1-4} turns into the
following system without delay,
\begin{equation}
\left\{
\begin{aligned} &v_t-a(l(v+\phi z_*(\theta_t\omega)))\Delta v-a(l(v+\phi z_*(\theta_t\omega)))z_*(\theta_t\omega)\Delta \phi\\ &\quad-\displaystyle\int_{0}^{\infty}\mu(s)\Delta \eta^t(s)ds +f(v+\phi z_*(\theta_t\omega))=h+\phi z_*(\theta_t\omega)+z_k^\phi(\theta_t\omega),&&\mbox{in}~~\mathcal{O}\times (\tau,\infty),\\ &\eta^t_t(s)=-\eta^t_s(s)+v(t),&&\mbox{in}~~\mathcal{O}\times(\tau,\infty), s>0,\\ &v(x,t)=\eta^t(x,s)=0,&&\mbox{on}~\partial\mathcal{O}\times(\tau,\infty), s>0,\\ &v(x,\tau)=v_0(x):=u_0(x)-\phi z_*(\theta_{\tau}\omega),&&\mbox{in}~~\mathcal{O},\\ &\eta^{\tau}(x,s)=\eta_0(x,s),&&\mbox{in}~~\mathcal{O}\times\mathbb{R}^+,\\ \end{aligned}\right.
\label{eq1-5}%
\end{equation}
where
\begin{equation}
{\eta_{0}(s)(\omega)}=\int_{\tau-s}^{\tau}v(x,r,\omega)dr=\int_{-s}%
^{0}(\varphi(r)-\phi z_{\ast}(\theta_{r+\tau}\omega))dr:=\int_{-s}^{0}%
\varphi_{v}(r)dr, \label{eq1-9}%
\end{equation}
which {contains} the initial integrated past history of $\varphi$ with
vanishing boundary {and a piece of value of $z_{\ast}(\theta_{\cdot+\tau
}\omega)$ on $(-s,0]$}. Moreover, $\eta_{s}^{t}$ denotes the distributional
derivative of $\eta^{t}(s)$ with respect to the internal variable $s$.

\subsection{Assumptions}

We will enumerate the assumptions on the nonlinear term $f$ and the variable
$\mu$. In our analysis, suppose that $h\in L^{2}(\mathcal{O})$, $f:\mathbb{R}%
\rightarrow\mathbb{R}$ is a polynomial of odd degree with positive leading
coefficient,
\begin{equation}
f(u)=\sum_{k=1}^{2p}f_{2p-k}u^{k-1},\qquad p\in\mathbb{N}. \label{eq2-1}%
\end{equation}
The variable $\mu$ is required to verify the following hypotheses:\\[1ex]%
$(h_{1})$ $\mu\in C^{1}(\mathbb{R}^{+})\cap L^{1}(\mathbb{R}^{+}),~~\mu
(s)\geq0,~~\forall s\in\mathbb{R}^{+}$;\\[1ex]$(h_{2})$ $\mu^{\prime
}(s)+\varpi\mu(s)\leq0$,~~$\forall s\in\mathbb{R}^{+}$,~~for some $\varpi>0$.

\begin{remark}
\label{rem2-1} (i) Recall that from $\mu(s)=-k^{\prime}(s)$, together with
condition $(h_{1})$, we infer immediately that there exists a constant
$M_{1}>0$ such that,
\[
k(t)=-\int_{t}^{\infty}k^{\prime}(s)ds=\int_{t}^{\infty}\mu(s)ds<M_{1}%
,\qquad\forall\ 0\leq t<\infty.
\]
(ii) In terms of assumption $(h_{2})$ imposed on $\mu$, it is easy to see
that,
\[
\mu(s_{2})\leq\mu(s_{1})e^{-\varpi(s_{2}-s_{1})},\qquad\forall\ 0<s_{1}%
<s_{2}.
\]
(iii) Combining the results of $(i)$ and $(ii)$, we have
\[
k(t)=\int_{t}^{\infty}\mu(s)ds\leq\mu(t)\int_{t}^{\infty}e^{-\varpi
(s-t)}ds:=\frac{\mu(t)}{\varpi},\qquad\forall\ 0<t<\infty.
\]

\end{remark}

\subsection{Notation}

Let $\mathcal{O}$ be a fixed bounded domain in $\mathbb{R}^{N}$ with regular
boundary. On this set, we introduce the Lebesgue space $L^{p}(\mathcal{O})$
with the natural norm $\Vert\cdot\Vert_{p}$, where $1\leq p\leq\infty$.
Besides, $W^{1,p}(\mathcal{O})$ is the subspace of $L^{p}(\mathcal{O})$
consisting of functions such that the first order weak derivative belongs to
$L^{p}(\mathcal{O})$. For convenience, $L^{2}(\mathcal{O})$ is denoted by $H$,
$H_{0}^{1}(\mathcal{O})$ is denoted by $V$ and $H^{-1}(\mathcal{O})$, the dual
space of $H_{0}^{1}(\mathcal{O})$, is denoted by $V^{\ast}$. We will use the
norms and inner products of $H$ and $V$ as $|\cdot|$, $\Vert\cdot\Vert$, and
$(\cdot,\cdot)$, $((\cdot,\cdot))$, respectively. Moreover, $<\cdot,\cdot>$
will denote the duality pairing between $V$ and $V^{\ast}$.

Taking into account system \eqref{eq1-5} and $(h_{1})$, we need to modify
slightly the notation before showing main results. Let $L_{\mu}^{2}%
(\mathbb{R}^{+};H)$ be a Hilbert space of functions $w:\mathbb{R}%
^{+}\rightarrow H$ endowed with the inner product,
\[
(w_{1},w_{2})_{\mu}=\int_{0}^{\infty}\mu(s)(w_{1}(s),w_{2}(s))ds,
\]
and let $|\cdot|_{\mu}$ denote the corresponding norm. In a similar way, we
introduce the inner products $((\cdot,\cdot))_{\mu}$, $(((\cdot,\cdot)))_{\mu
}$ and relative norms $\Vert\cdot\Vert_{\mu}$, $|||\cdot|||_{\mu}$ on $L_{\mu
}^{2}(\mathbb{R}^{+};V)$, $L_{\mu}^{2}(\mathbb{R}^{+};V\cap H^{2}%
(\mathcal{O}))$, respectively. It follows then that
\[
((\cdot,\cdot))_{\mu}=(\nabla\cdot,\nabla\cdot)_{\mu},~~\mbox{and}~~(((\cdot
,\cdot)))_{\mu}=(\Delta\cdot,\Delta\cdot)_{\mu}.
\]
We also define the Hilbert spaces,
\[
\mathcal{H}=H\times L_{\mu}^{2}(\mathbb{R}^{+};V),
\]
and
\[
\mathcal{V}=V\times L_{\mu}^{2}(\mathbb{R}^{+};V\cap H^{2}(\mathcal{O})),
\]
which are respectively endowed with the inner products,
\[
(w_{1},w_{2})_{\mathcal{H}}=(w_{1},w_{2})+((w_{1},w_{2}))_{\mu},
\]
and
\[
(w_{1},w_{2})_{\mathcal{V}}=((w_{1},w_{2}))+(((w_{1},w_{2})))_{\mu},
\]
where $w_{i}\in\mathcal{H}$ or $\mathcal{V}$ $(i=1,2)$. The norms induced on
$\mathcal{H}$ and $\mathcal{V}$ are the so-called energy ones, which read
\[
\Vert(w_{1},w_{2})\Vert_{\mathcal{H}}^{2}=|w_{1}|^{2}+\int_{0}^{\infty}%
\mu(s)\Vert w_{2}(s)\Vert^{2}ds,
\]
and
\[
\Vert(w_{1},w_{2})\Vert_{\mathcal{V}}^{2}=\Vert w_{1}\Vert^{2}+\int%
_{0}^{\infty}\mu(s)\Vert\nabla w_{2}(s)\Vert^{2}ds.
\]

At last, with standard notation, $\mathcal{D}(I;X)$ is the space of infinitely
differentiable $X$-valued functions with compact support in $I\subset
\mathbb{R}$, whose dual space is the distribution one on $I$ with values in
$X^{\ast}$ (dual of $X$), denoted by $\mathcal{D}^{\prime}(I;X^{\ast})$. For
convenience, we define $L_{V}^{2}$ as the space of functions $u\left(
\text{\textperiodcentered}\right)  $ such that,%
\[
\int_{-\infty}^{0}e^{\gamma s}\left\Vert u\left(  s\right)  \right\Vert
^{2}ds<\infty,
\]
where $0<\gamma<\min\{\frac{m\lambda_{1}}{2},\varpi\}$, $\lambda_{1}$ is the
first eigenvalue of $-\Delta$ with zero Dirichlet boundary condition and
$\varpi$ comes from $(h_{2})$.

\section{Well-posedness of problem \eqref{eq1-1}}

\label{s3}

In this section, before presenting the well-posedness of problem
\eqref{eq1-1}, we first state an auxiliary result for the regularity of
initial value $\eta_{0}$ (cf. \eqref{eq1-9}), which is the essential point to
prove the existence and uniqueness of solutions to problem \eqref{eq1-5}.

\subsection{An auxiliary result}

Let us first recall a crucial technical lemma in \cite[Lemma 3.1]{X2}.

\begin{lemma}
\label{lem3-1} Assume $(h_{1})$-$(h_{2})$ hold. Then, the operator
$\mathcal{J}:L_{V}^{2}\rightarrow L_{\mu}^{2}(\mathbb{R}^{+};V)$ defined by,
\[
(\mathcal{J}\varphi)(s)=\int_{-s}^{0}\varphi(r)dr,\qquad s\in\mathbb{R}^{+},
\]
is a linear and continuous operator. In particular, there exists a positive
constant $K_{\mu}=e^{\gamma}\int_{0}^{1}\mu(s)ds+\mu(1)e^{\varpi}%
(\gamma-\varpi)^{-2}$ such that for any $\varphi\in L_{V}^{2}$, it holds
\[
\Vert\mathcal{J}\varphi\Vert_{L_{\mu}^{2}(\mathbb{R}^{+};V)}^{2}\leq K_{\mu
}\Vert\varphi\Vert_{L_{V}^{2}}^{2}.
\]

\end{lemma}

For the sake of simplicity, define $(\mathcal{J}_{\omega,\tau}{\varphi}
)(s):=\mathcal{J}(\varphi-\phi z_{*}(\theta_{\cdot+\tau}\omega))(s)$. By
slightly modifying the proof of Lemma \ref{lem3-1} in \cite{X2}, we have the
following corollary.

\begin{corollary}
\label{cor3-2} Assume $(h_{1})$-$(h_{2})$ hold and $\phi\in V\cap
H^{2}(\mathcal{O})\cap L^{2p}(\mathcal{O})$. Then, for every $\omega\in\Omega$
and $\tau\in\mathbb{R}$, the operator $\mathcal{J}_{\omega,\tau}:L_{V}%
^{2}\rightarrow L_{\mu}^{2}(\mathbb{R}^{+};V)$ defined by%
\begin{equation}
(\mathcal{J}_{\omega,\tau}\varphi)(s):=\int_{-s}^{0}\varphi(r)dr-\int_{-s}%
^{0}z_{\ast}(\theta_{r+\tau}\omega)\phi dr=\mathcal{J}(\varphi-\phi z_{\ast
}(\theta_{\text{\textperiodcentered}+\tau}\omega))(s), \label{eq3-0}%
\end{equation}
is continuous. Additionally, there exists a positive constant $K_{\mu}$ which
is the same as in Lemma \ref{lem3-1}, such that for any $\varphi\in L_{V}^{2}%
$, we have
\[
\Vert\mathcal{J}_{\omega,\tau}\varphi\Vert_{L_{\mu}^{2}(\mathbb{R}^{+};V)}%
^{2}\leq K_{\mu}\left(  \Vert\varphi-z_{\ast}(\theta_{\cdot+\tau}\omega
)\phi\Vert_{L_{V}^{2}}^{2}\right)  \leq2K_{\mu}\left(  \Vert\varphi
\Vert_{L_{V}^{2}}^{2}+\Vert z_{\ast}(\theta_{\cdot+\tau}\omega)\phi
\Vert_{L_{V}^{2}}^{2}\right)  .
\]

\end{corollary}

\textbf{Proof.} First of all, we show that the operator $\mathcal{J}%
_{\omega,\tau}$ is well-defined. As $\phi\in V\cap H^{2}(\mathcal{O})\cap
L^{2p}(\mathcal{O})$, it is easy to check that $z_{\ast}(\theta_{\cdot+\tau
}\omega)\phi\in L_{V}^{2}$. Indeed, for every $\omega\in\Omega$,
\[%
\begin{split}
\Vert z_{\ast}(\theta_{\cdot+\tau}\omega)\phi\Vert_{L_{V}^{2}}^{2}  &
=\int_{-\infty}^{0}e^{\gamma t}|z_{\ast}(\theta_{t+\tau}\omega)|^{2}\Vert
\phi\Vert^{2}dt=\Vert\phi\Vert^{2}\int_{-\infty}^{\tau}e^{\gamma(t-\tau
)}|z_{\ast}(\theta_{t}\omega)|^{2}dt\\
&  \leq e^{-\gamma\tau}\Vert\phi\Vert^{2}\left(  \int_{-\infty}^{0}e^{\gamma
t}|z_{\ast}(\theta_{t}\omega)|^{2}dt+\int_{0}^{\tau}e^{\gamma t}|z_{\ast
}(\theta_{t}\omega)|^{2}dt\right)  <+\infty,
\end{split}
\]
where we have used the continuity property of $z_{\ast}(\theta_{t}\omega)$
with respect to $t$ and \eqref{zproperty}. Notice that, by Remark
\ref{rem2-1}$(ii)$, we find%
\[%
\begin{split}
&  \qquad\left\Vert \int_{-s}^{0}\varphi-z_{\ast}(\theta_{r+\tau}\omega)\phi
dr\right\Vert _{L_{\mu}^{2}(\mathbb{R}^{+};V)}^{2}\leq\int_{0}^{\infty}%
\mu(s)\left(  \int_{-s}^{0}\Vert\varphi(r)-z_{\ast}(\theta_{r+\tau}\omega
)\phi\Vert dr\right)  ^{2}ds\\[0.8ex]
&  \leq\int_{0}^{1}\mu(s)s\int_{-s}^{0}\Vert\varphi(r)-z_{\ast}(\theta
_{r+\tau}\omega)\phi\Vert^{2}drds+\int_{1}^{\infty}\mu(s)\left(  \int_{-s}%
^{0}\Vert\varphi(r)-z_{\ast}(\theta_{r+\tau}\omega)\phi\Vert dr\right)
^{2}ds\\[0.8ex]
&  \leq\int_{-1}^{0}\Vert\varphi(r)-z_{\ast}(\theta_{r+\tau}\omega)\phi
\Vert^{2}e^{\gamma r}e^{-\gamma r}\int_{0}^{1}s\mu(s)dsdr\\[0.8ex]
&  ~~+\mu(1)e^{\varpi}\int_{-\infty}^{0}e^{\gamma r}\Vert\varphi(r)-z_{\ast
}(\theta_{r+\tau}\omega)\phi\Vert^{2}\int_{-r}^{\infty}se^{-\varpi s}e^{\gamma
s}dsdr\\[0.8ex]
&  \leq\Vert\varphi-z_{\ast}(\theta_{\cdot+\tau}\omega)\phi\Vert_{L_{V}^{2}%
}^{2}e^{\gamma}\int_{0}^{1}s\mu(s)ds+\mu(1)e^{\varpi}(\varpi-\gamma)^{-2}%
\Vert\varphi-z_{\ast}(\theta_{\cdot+\tau}\omega)\phi\Vert_{L_{V}^{2}}%
^{2}\\[0.8ex]
&  \leq\left(  e^{\gamma}\int_{0}^{1}\mu(s)ds+\mu(1)e^{\varpi}(\varpi
-\gamma)^{-2}\right)  \Vert\varphi-z_{\ast}(\theta_{\cdot+\tau}\omega
)\phi\Vert_{L_{V}^{2}}^{2}=K_{\mu}\Vert\varphi-z_{\ast}(\theta_{\cdot+\tau
}\omega)\phi\Vert_{L_{V}^{2}}^{2}.
\end{split}
\]
The proof of this corollary is complete. $\Box$

\begin{remark}
\label{rem3-2} Once we fix an initial function $\varphi\in L_{V}^{2}$ to
problem \eqref{eq1-1}, then for every $\omega\in\Omega$, the initial function
$\varphi_{v}:=\varphi-\phi z_{\ast}(\theta_{\cdot+\tau}\omega)$ of
\eqref{eq1-4} belongs to $L_{V}^{2}$ as $\phi\in V\cap H^{2}(\mathcal{O})\cap
L^{2p}(\mathcal{O})$. Also, the initial value for the second component
$\eta_{0}(s,\omega)=(\mathcal{J}_{\omega,\tau}\varphi)(s)$ of problem
\eqref{eq1-5} belongs to $L_{\mu}^{2}(\mathbb{R}^{+};V)$ thanks to Corollary
\ref{cor3-2}. Analogously, if we pick up the initial function $\varphi\in
L_{V\cap H^{2}({\mathcal{O}})}^{2}$ to problem \eqref{eq1-1}, then the initial
function $\varphi_{v}:=\varphi-\phi z_{\ast}(\theta_{\cdot+\tau}\omega)$ of
\eqref{eq1-4} also belongs to $L_{V\cap H^{2}({\mathcal{O}})}^{2}$. Thus,
making use of a similar proof as in Corollary \ref{cor3-2}, it is easy to
check that $\eta_{0}(s,\omega)$ defined by \eqref{eq3-0} belongs to $L_{\mu
}^{2}(\mathbb{R}^{+};V\cap H^{2}(\mathcal{O}))$.
\end{remark}

%$\varphi-\phi z_*(\theta_{\cdot+\tau}\omega)\in L_{V\cap H^2({\mathcal{O}})}^2$

%We will first prove well-posedness of problem \eqref{eq1-5} and then obtain the existence and uniqueness of solutions to \eqref{eq1-1} via the transform \eqref{eq1-3}.

\subsection{Well-posedness of problem \eqref{eq1-5}}

\begin{theorem}
\label{thm3-3} Assume that \eqref{eq1-2}, \eqref{eq2-1} and $(h_{1})$%
-$(h_{2})$ hold. Let $\phi\in V\cap H^{2}(\mathcal{O})\cap L^{2p}%
(\mathcal{O})$ be such that $\Delta\phi\in L^{2p}(\mathcal{O})$, let $h\in H$
and $a$ be a locally Lipschitz function. Then:

\begin{enumerate}
\item[$(i)$] For every $\omega\in\Omega$, any initial value $v_{0}\in H$ and
initial function $\varphi\in L_{V}^{2}$, there exists a unique solution
$(v,\eta)$ to problem \eqref{eq1-5} in the weak sense {with initial value
$(v_{0},\eta_{0})$,} where {\ $\eta_{0}(s,\omega)=(\mathcal{J}_{\omega,\tau
}\varphi)(s),$} fulfilling
\[
v\in L^{\infty}(\tau,T;H)\cap L^{2}(\tau,T;V)\cap L^{2p}(\tau,T;L^{2p}%
(\mathcal{O})),\qquad\forall T>\tau,
\]%
\[
\eta\in L^{\infty}(\tau,T;L_{\mu}^{2}(\mathbb{R}^{+};V)),\qquad\forall
T>\tau.
\]
Furthermore, the solution $(v,\eta)$ of \eqref{eq1-5} is continuous with
respect to the initial value $(v_{0},\eta_{0})$ for all $t\in\lbrack\tau,T]$
in $\mathcal{H}$;

\item[$(ii)$] For any initial value $(v_{0},\eta_{0})\in\mathcal{V}$, the
unique solution $(v,\eta)$ to problem \eqref{eq1-5} satisfies,
\[
v\in L^{\infty}(\tau,T;V)\cap L^{2}(\tau,T;V\cap H^{2}(\mathcal{O}%
)),\qquad\forall T>\tau,
\]%
\[
\eta\in L^{\infty}(\tau,T;L_{\mu}^{2}(\mathbb{R}^{+};V\cap H^{2}%
(\mathcal{O}))),\qquad\forall T>\tau.
\]
In addition, the solution $(v,\eta)$ of \eqref{eq1-5} is continuous with
respect to the initial value $(v_{0},\eta_{0})$ for all $t\in\lbrack\tau,T]$
in $\mathcal{V}$.
\end{enumerate}
\end{theorem}

\begin{remark}
By a similar argument as in \cite[p.443]{X2} one can prove that%
\[
v\in C([\tau,T],H),\ \eta\in C([\tau,T],L_{\mu}^{2}(\mathbb{R}^{+};V)).
\]

\end{remark}

\textbf{Proof.} On the one hand, for every $\omega\in\Omega$, it follows from
Corollary \ref{cor3-2} that $(\mathcal{J}_{\omega,\tau}{\varphi})(s)\in
L_{\mu}^{2}(\mathbb{R}^{+};V)$ thanks to the facts $\varphi\in L_{V}^{2}$ and
$\phi\in V\cap H^{2}(\mathcal{O})\cap L^{2p}(\mathcal{O})$. {We will prove
$(i)$ in five steps. }

\textbf{Step 1. (Faedo-Galerkin scheme)} Let $\{w_{j}\}_{j=1}^{\infty}$ be the
eigenfunctions of the operator $-\Delta$, which is an orthonormal basis in
$H$. For regular domains these functions belong to $V\cap L^{2p}(\mathcal{O})$
\cite[Chapter 8]{R1}. We select also an orthonormal basis $\{\zeta_{j}%
\}_{j=1}^{\infty}$ of $L_{\mu}^{2}(\mathbb{R}^{+};V)$ which {belongs} to
$\mathcal{D}(\mathbb{R}^{+};V)$ as well. Fix $T>\tau$. For a given integer $n$
and each $\omega$, denote by $P_{n}$ and $Q_{n}$ the projectors on the
subspaces,
\[
\mbox{span}\{w_{1},\cdots,w_{n}\}\subset V,\qquad\mbox{and}\qquad
\mbox{span}\{\zeta_{1},\cdots,\zeta_{n}\}\subset L_{\mu}^{2}(\mathbb{R}%
^{+};V),
\]
respectively. We will look for a function $(v_{n},\eta_{n})$ of the form,
\[
v_{n}(t)=\sum_{j=1}^{n}b_{j}(t)w_{j}\qquad\mbox{and}\qquad\eta_{n}^{t}%
(s)=\sum_{j=1}^{n}c_{j}(t)\zeta_{j}(s),
\]
satisfying for all $t\geq\tau$,
\begin{equation}
\left\{
\begin{aligned} &\frac{d}{dt}b_k(t)+\lambda_ka(l(\sum_{j=1}^nb_j(t)w_j) +l(\phi)z_*(\theta_t\omega))b_k(t) +a(l(\sum_{j=1}^nb_j(t)w_j)+l(\phi)z_*(\theta_t\omega))z_*(\theta_t\omega) \\ &\qquad ~\quad\times \lambda_k<\phi,w_k>+\sum_{j=1}^nc_j(t)((\zeta_j,w_k))_\mu+(f(\sum_{j=1}^nb_j(t)w_j+\phi z_*(\theta_t\omega)),w_k)\\ &\qquad~\quad=z_*(\theta_t\omega)(\phi,w_k)+(z_k^\phi(\theta_t\omega),w_k)+(h,w_k),\\ &\frac{d}{dt}c_k(t)=-\sum_{j=1}^nc_j(t)((\zeta_j^\prime,\zeta_k))_\mu+\sum_{j=1}^nb_j(t)((w_j,\zeta_k))_\mu,~~k=1,2,\cdots,n, \end{aligned}\right.
\label{eq3-1}%
\end{equation}
where $\lambda_{j}$ is the eigenvalue associated to the eigenfunction $w_{j}$,
subject to the initial conditions
\begin{equation}
b_{k}(\tau)=(v_{0},w_{k}),\qquad c_{k}(\tau)=((\eta_{0},\zeta_{k}))_{\mu}.
\label{eq3-2}%
\end{equation}
According to the standard existence theory for ordinary differential systems,
there exists a continuous solution of \eqref{eq3-1}-\eqref{eq3-2} on some
interval $(\tau,t_{n})$. Proceeding as usual, by establishing some a priori
estimates below, we can ensure that $t_{n}=\infty$.

\textbf{Step 2. (Energy estimates)} Multiplying the first equation of
\eqref{eq3-1} by $b_{k}$ and the second one by $c_{k}$, respectively, summing
over $k$ $(k=1,2,\cdots,n)$ and adding the results, we have
\[
\begin{aligned} &\frac{1}{2}\frac{d}{dt}|v_n(t)|^2+a(l( v_n(t))+l(\phi)z_*(\theta_t\omega))\|v_n(t)\|^2+a(l(v_n(t))+l(\phi)z_*(\theta_t\omega))z_*(\theta_t\omega) ((\phi,v_n))\\[1.0ex] &\qquad\qquad ~~~+((\eta_n^t(s),v_n(t)))_\mu+(f(v_n(t)+\phi z_*(\theta_t\omega)),v_n(t))\\[1.0ex] &\qquad\qquad~~=z_*(\theta_t\omega)(\phi,v_n(t))+(z_k^\phi(\theta_t\omega),v_n(t))+(h,v_n(t)),\\[1.0ex] &\frac{1}{2}\frac{d}{dt}\|\eta^t_n(s)\|_\mu^2=-((\eta^t_n(s), (\eta^t_n(s))^\prime))_\mu+((\eta_n^t(s),v_n(t)))_\mu. \end{aligned}
\]
Combining the two equations above, for every fixed $\omega\in\Omega$ and
$t\in\lbrack\tau,T]$, by condition \eqref{eq1-2}, we obtain
\[%
\begin{split}
\frac{1}{2}\frac{d}{dt}|v_{n}(t)|^{2}  &  +\frac{1}{2}\frac{d}{dt}\Vert
\eta_{n}^{t}(s)\Vert_{\mu}^{2}+m\Vert v_{n}(t)\Vert^{2}+a(l(v_{n}%
(t))+l(\phi)z_{\ast}(\theta_{t}\omega))z_{\ast}(\theta_{t}\omega)((\phi
,v_{n}(t)))\\[0.8ex]
&  +(f(v_{n}(t)+\phi z_{\ast}(\theta_{t}\omega)),v_{n}(t))\\[0.8ex]
&  \leq-((\eta_{n}^{t}(s),(\eta_{n}^{t}(s))^{\prime}))_{\mu}+|z_{\ast}%
(\theta_{t}\omega)||\phi||v_{n}(t)|+(z_{k}^{\phi}(\theta_{t}\omega
),v_{n}(t))+(h,v_{n}\left(  t\right)  ),
\end{split}
\]
which is equivalent to
\begin{equation}%
\begin{split}
\frac{d}{dt}|v_{n}(t)|^{2}  &  +\frac{d}{dt}\Vert\eta_{n}^{t}(s)\Vert_{\mu
}^{2}+2m\Vert v_{n}(t)\Vert^{2}+2a(l(v_{n}(t))+l(\phi)z_{\ast}(\theta
_{t}\omega))z_{\ast}(\theta_{t}\omega)((\phi,v_{n}(t)))\\[0.8ex]
&  +2(f(v_{n}(t)+\phi z_{\ast}(\theta_{t}\omega)),v_{n}(t))\\[0.8ex]
&  \leq-2((\eta_{n}^{t}(s),(\eta_{n}^{t}(s))^{\prime}))_{\mu}+2|z_{\ast
}(\theta_{t}\omega)||\phi||v_{n}(t)|+2(z_{k}^{\phi}(\theta_{t}\omega
),v_{n}(t))+2(h,v_{n}\left(  t\right)  ).
\end{split}
\label{eq3-3}%
\end{equation}
Let us do estimates for \eqref{eq3-3} one by one. First of all, by integration
by parts, we have
\begin{equation}
-2((\eta_{n}^{t}(s),(\eta_{n}^{t}(s))^{\prime}))_{\mu}=\int_{0}^{\infty}%
\mu^{\prime}(s)|\nabla\eta_{n}^{t}(s)|^{2}ds\leq0. \label{eq3-4}%
\end{equation}
Second, by \eqref{eq1-2} and the Young inequality, we obtain
\begin{equation}%
\begin{split}
-2a(l(v_{n}(t))+l(\phi)z_{\ast}(\theta_{t}\omega))z_{\ast}(\theta_{t}%
\omega)((\phi,v_{n}(t)))  &  \leq2M|z_{\ast}(\theta_{t}\omega)|\Vert\phi
\Vert\Vert v_{n}(t)\Vert\\
&  \leq\frac{m}{4}\Vert v_{n}(t)\Vert^{2}+\frac{4M^{2}}{m}|z_{\ast}(\theta
_{t}\omega)|^{2}\Vert\phi\Vert^{2}.
\end{split}
\label{eq3-5}%
\end{equation}
Third, by the Poincar\'{e} and Young inequalities, we derive%
\begin{equation}
2|z_{\ast}(\theta_{t}\omega)||\phi||v_{n}(t)|\leq2|z_{\ast}(\theta_{t}%
\omega)||\phi|\frac{\Vert v_{n}(t)\Vert}{\sqrt{\lambda_{1}}}\leq\frac{m}%
{4}\Vert v_{n}(t)\Vert^{2}+\frac{4}{m\lambda_{1}}|z_{\ast}(\theta_{t}%
\omega)|^{2}|\phi|^{2}. \label{eq3-6}%
\end{equation}
Fourth, for the nonlinear term $f$, making use of the Young inequality, we
infer that there exist constants $\alpha,\beta>0$, such that
\begin{equation}
f(u)u\geq\frac{1}{2}f_{0}u^{2p}-\alpha,\qquad\mbox{and}\qquad|f(u)|\leq
\beta(1+|u|^{2p-1}). \label{fcondition}%
\end{equation}
Therefore, by the continuity of $z_{\ast}(\theta_{t}\omega)$,
\eqref{fcondition} and the Young inequality, we deduce there exists a constant
$\tilde{C}_{1}$ such that
\begin{equation}%
\begin{split}
&  ~\quad2(f(v_{n}(t)+\phi z_{\ast}(\theta_{t}\omega)),v_{n}(t))\\[0.8ex]
&  =2\int_{\mathcal{O}}f(v_{n}(t)+\phi z_{\ast}(\theta_{t}\omega))\left(
v_{n}(t)+\phi z_{\ast}(\theta_{t}\omega)\right)  dx-2\int_{\mathcal{O}}%
f(v_{n}(t)+\phi z_{\ast}(\theta_{t}\omega))\phi z_{\ast}(\theta_{t}\omega)dx\\
&  \geq f_{0}\int_{\mathcal{O}}\left\vert v_{n}(t)+\phi z_{\ast}(\theta
_{t}\omega)\right\vert _{2p}^{2p}dx-\widetilde{C}_{1}\int_{\mathcal{O}}\left(
1+|v_{n}(t)|^{2p-1}+|z_{\ast}(\theta_{t}\omega)|^{2p-1}|\phi|^{2p-1}\right)
|z_{\ast}(\theta_{t}\omega)||\phi|dx\\[0.8ex]
&  -2\alpha|\mathcal{O}|.
\end{split}
\label{eq3-9}%
\end{equation}
Since%
\[
\left\vert v\right\vert ^{2p}=\left\vert v+r-r\right\vert ^{2p}\leq D\left(
\left\vert v+r\right\vert ^{2p}+\left\vert r\right\vert ^{2p}\right)  ,
\]
for some $D=D(p)>0$, we obtain%
\begin{align*}
&  ~~\quad2(f(v_{n}(t)+\phi z_{\ast}(\theta_{t}\omega)),v_{n}(t))\\[0.8ex]
&  \geq\frac{f_{0}}{D}\int_{\mathcal{O}}\left\vert v_{n}(t)\right\vert
^{2p}dx-\frac{f_{0}}{D}|z_{\ast}(\theta_{t}\omega)|^{2p}\left\Vert
\phi\right\Vert _{2p}^{2p}\\
&  -\widetilde{C}_{1}\max\{|z_{\ast}(\theta_{t}\omega)|,|z_{\ast}(\theta
_{t}\omega)|^{2p}\}\int_{\mathcal{O}}\left(  1+|v_{n}(t)|^{2p-1}+|\phi
|^{2p-1}\right)  |\phi|dx-2\alpha|\mathcal{O}|\\
&  \geq\frac{f_{0}}{2D}\Vert v_{n}(t)\Vert_{2p}^{2p}-{C_{1}(\theta_{t}%
\omega)(1+}\Vert\phi\Vert_{2p}^{2p})-2\alpha|\mathcal{O}|.
\end{align*}
Here, {$C_{1}(\omega):=C_{1}(|z_{\ast}(\omega)|,p,\left\vert \mathcal{O}%
\right\vert )=\widetilde{C}_{2}({1+}$}$|z_{\ast}(\omega)|^{4p^{2}}),$ for some
$\widetilde{C}_{2}=\widetilde{C}_{2}(p,\left\vert \mathcal{O}\right\vert )>0.$
As for the last term, by the Young inequality, the properties of $z_{\ast
}(\theta_{t}\omega)$ (cf. \eqref{zproperty}) and Remark \ref{rem2-1}, we
deduce that there exists a random variable $C_{2}(\omega)$, such that
\begin{equation}%
\begin{split}
2(z_{k}^{\phi}(\theta_{t}\omega),v_{n}(t))  &  \leq2\left(  \int_{0}^{\infty
}k(s)|z_{\ast}(\theta_{t-s}\omega)|ds\right)  \Vert\phi\Vert\Vert
v_{n}(t)\Vert\\[0.8ex]
&  \leq2\left(  \int_{0}^{1}k(s)|z_{\ast}(\theta_{t-s}\omega)|ds+\frac
{1}{\varpi}\int_{1}^{\infty}\mu(s)|z_{\ast}(\theta_{t-s}\omega)|ds\right)
\Vert\phi\Vert\Vert v_{n}(t)\Vert\\[0.8ex]
&  \leq2\left(  M_{1}\int_{0}^{1}|z_{\ast}(\theta_{t-s}\omega)|ds+\frac
{\mu(1)e^{\varpi}}{\varpi}\int_{1}^{\infty}e^{-\varpi s}|z_{\ast}(\theta
_{t-s}\omega)|ds\right)  \Vert\phi\Vert\Vert v_{n}(t)\Vert\\
&  =C_{2}(\theta_{t}\omega)\Vert\phi\Vert\Vert v_{n}(t)\Vert\leq\frac
{C_{2}(\theta_{t}\omega)^{2}}{m}\Vert\phi\Vert^{2}+\frac{m}{4}\Vert
v_{n}(t)\Vert^{2}.
\end{split}
\label{eq3-10}%
\end{equation}
Also, by the Young inequality, we have
\begin{equation}
{2}\left(  h,v_{n}\right)  \leq\frac{{4}}{m\lambda_{1}}|h|^{2}+\frac{m}%
{4}\Vert v_{n}(t)\Vert^{2}. \label{eq3-10b}%
\end{equation}
Substituting \eqref{eq3-4}-\eqref{eq3-10b} into \eqref{eq3-3}, we obtain
\begin{equation}%
\begin{split}
&  ~\quad\frac{d}{dt}|v_{n}(t)|^{2}+\frac{d}{dt}\Vert\eta_{n}^{t}(s)\Vert
_{\mu}^{2}+m\Vert v_{n}(t)\Vert^{2}+\frac{f_{0}}{2D}\Vert v_{n}(t)\Vert
_{2p}^{2p}\leq\frac{4M^{2}}{m}|z_{\ast}(\theta_{t}\omega)|^{2}\Vert\phi
\Vert^{2}\\[0.8ex]
&  +\frac{4}{m\lambda_{1}}|z_{\ast}(\theta_{t}\omega)|^{2}|\phi|^{2}%
+2\alpha|\mathcal{O}|+{C_{1}(\theta_{t}\omega)(1+}\Vert\phi\Vert_{2p}%
^{2p})+\frac{{C_{2}(\theta_{t}\omega)^{2}}}{m}\Vert\phi\Vert^{2}+\frac{{4}%
}{m\lambda_{1}}|h|^{2}.
\end{split}
\label{eq3-11}%
\end{equation}
Denote
\begin{equation}
\Theta_{1}(\omega)=\frac{4M^{2}}{m}|z_{\ast}(\omega)|^{2}\Vert\phi\Vert
^{2}+\frac{4}{m\lambda_{1}}|z_{\ast}(\omega)|^{2}|\phi|^{2}+2\alpha
|\mathcal{O}|+{C_{1}(\omega)(1+}\Vert\phi\Vert_{2p}^{2p})+\frac{{C_{2}%
(\omega)}^{2}}{m}\Vert\phi\Vert^{2}+\frac{4}{m\lambda_{1}}|h|^{2}.
\label{theta1}%
\end{equation}
Subsequently, it follows from \eqref{eq3-11} that
\[
\frac{d}{dt}|v_{n}(t)|^{2}+\frac{d}{dt}\Vert\eta_{n}^{t}(s)\Vert_{\mu}%
^{2}+m\Vert v_{n}(t)\Vert^{2}+\frac{f_{0}}{{2D}}\Vert v_{n}(t)\Vert_{2p}%
^{2p}\leq\Theta_{1}(\theta_{t}\omega).
\]
Denote $y_{n}(t):=(v_{n}(t),\eta_{n}^{t}(s))$, then $y_{n}(\tau):=y_{0}%
^{n}=(P_{n}v_{0},Q_{n}\eta_{0})$. Integrating the above inequality over
$(\tau,T)$ for every $T>\tau$ and $\omega\in\Omega$, we have
\begin{equation}
\Vert y_{n}(T)\Vert_{\mathcal{H}}^{2}+m\int_{\tau}^{T}\Vert v_{n}(t)\Vert
^{2}dt+\frac{f_{0}}{{2D}}\int_{\tau}^{T}\Vert v_{n}(t)\Vert_{2p}^{2p}%
dt\leq\Vert y_{0}^{n}\Vert_{\mathcal{H}}^{2}+\int_{\tau}^{T}\Theta_{1}%
(\theta_{t}\omega)dt. \label{eq3-12}%
\end{equation}
Thanks to $\phi\in V\cap H^{2}(\mathcal{O})\cap L^{2p}(\mathcal{O})$, together
with the fact that for every $\omega\in\Omega$, $z_{\ast}(\theta_{t}\omega)$
is continuous with respect to $t$, we infer that for every $T>\tau$,
$\Theta_{1}(\theta_{\cdot}\omega)\in L^{1}(\tau,T)$. Hence, for each
$\omega\in\Omega$, there exists a constant $C_{3}(\omega,T)>0$, such that
\begin{equation}
\Vert y_{n}(T)\Vert_{\mathcal{H}}^{2}+m\int_{\tau}^{T}\Vert v_{n}(t)\Vert
^{2}dt+\frac{f_{0}}{2{D}}\int_{\tau}^{T}\Vert v_{n}(t)\Vert_{2p}^{2p}dt\leq
C_{3}(\omega,T). \label{eq3-13}%
\end{equation}
Making use of (\ref{fcondition}), a compactness argument and the Aubin-Lions
lemma, for every $\omega\in\Omega$, there exist subsequences $\{v_{n}\}$ and
$\{\eta_{n}\}$ (relabeled the same), $v\in L^{\infty}(\tau,T;H)\cap L^{2}%
(\tau,T;V)\cap L^{2p}(\tau,T;L^{2p}(\mathcal{O}))$ and $\eta\in L^{\infty
}(\tau,T;L_{\mu}^{2}(\mathbb{R}^{+};V))$, such that%
\begin{align}
v_{n}  &  \rightarrow v~~\mbox{weak-$*$ in}~~L^{\infty}(\tau,T;H);\nonumber\\
v_{n}  &  \rightarrow v~~\mbox{weakly in}~~L^{2}(\tau,T;V);\nonumber\\
v_{n}  &  \rightarrow v~~\mbox{weakly in}~~L^{2p}(\tau,T;L^{2p}(\mathcal{O}%
));\nonumber\\
\eta_{n}  &  \rightarrow\eta~~\mbox{weak-$*$ in}~~L^{\infty}(\tau,T;L_{\mu
}^{2}(\mathbb{R}^{+};V));\nonumber\\
\frac{dv_{n}}{dt}  &  \rightarrow\frac{dv}{dt}~~\mbox {weakly in}~~L^{2}%
(\tau,T;V^{\ast})+L^{q}(\tau,T;L^{q}(\mathcal{O}));\label{eq3-14}\\
v_{n}  &  \rightarrow v~~\mbox{strongly in}~~L^{2}(\tau,T;H);\nonumber\\
v_{n}(t,\omega)  &  \rightarrow v(t,\omega
)~~\mbox{strongly in}~~H,~~\mbox{a.e.}~t\in(\tau,T];\nonumber\\
v_{n}(x,t,\omega)  &  \rightarrow v(x,t,\omega)~~\mbox{a.e.}~(x,t)\in
\mathcal{O}\times(\tau,T];\nonumber\\
f(v_{n}+\phi z_{\ast}(\theta_{\text{\textperiodcentered}}\omega))  &
\rightarrow f(v+\phi z_{\ast}(\theta_{\text{\textperiodcentered}}%
\omega))\text{ weakly in \ }L^{q}(\tau,T;L^{q}(\mathcal{O})),\nonumber
\end{align}
for all $T>\tau$, where $q=\frac{2p}{2p-1}$ is the conjugate number of $2p$.

\textbf{Step 3. (Passage to limit)} For a fixed integer $m$ and each
$\omega\in\Omega$, choose a function
\[
l=(\sigma,\pi)\in\mathcal{D}((\tau,T);V\cap L^{2p}(\mathcal{O}))\times
\mathcal{D}((\tau,T);\mathcal{D}(\mathbb{R}^{+};V)),
\]
of the form
\[
\sigma(t,\omega)=\sum_{j=1}^{m}\tilde{b}_{j}(t)w_{j}\qquad\mbox{and}\qquad
\pi^{t}(s,\omega)=\sum_{j=1}^{m}\tilde{c}_{j}(t)\zeta_{j}(s),
\]
where $\{\tilde{b}_{j}\}_{j=1}^{m}$ and $\{\tilde{c}_{j}\}_{j=1}^{m}$ are
given functions in $\mathcal{D}((\tau,T))$.

Our main target is to prove that problem \eqref{eq1-5} has a solution in the
weak sense, i.e., for arbitrary $l\in\mathcal{D}((\tau,T);V\cap L^{2p}%
(\mathcal{O}))\times\mathcal{D}((\tau,T);\mathcal{D}(\mathbb{R}^{+};V))$, the
following equality%
\begin{align*}
&  \qquad\int_{\tau}^{T}(\partial_{t}v_{n},\sigma)dt+\int_{\tau}^{T}%
((\partial_{t}\eta_{n}^{t},\pi))_{\mu}dt\\
&  =-\int_{\tau}^{T}\big(a(l(v_{n})+l(\phi)z_{\ast}(\theta_{t}\omega
))((v_{n},\sigma))+a(l(v_{n})+l(\phi)z_{\ast}(\theta_{t}\omega))z_{\ast
}(\theta_{t}\omega)((\phi,\sigma))\big)dt\\
&  -\int_{\tau}^{T}((\eta_{n}^{t},\sigma))_{\mu}dt-\int_{\tau}^{T}%
(f(v_{n}+\phi z_{\ast}(\theta_{t}\omega)),\sigma)dt+\int_{\tau}^{T}z_{\ast
}(\theta_{t}\omega)(\phi,\sigma)dt\\
&  +\int_{\tau}^{T}(z_{k}^{\phi}(\theta_{t}\omega),\sigma)dt {+\int_{\tau}%
^{T}(h,\sigma)dt}-\int_{\tau}^{T}\left\langle \left\langle (\eta_{n}%
^{t})^{\prime},\pi\right\rangle \right\rangle dt+\int_{\tau}^{T}((v_{n}%
,\pi))_{\mu}dt,
\end{align*}
holds in the space $\mathcal{D}^{\prime}((\tau,T))$. Here, we denote by
$\left\langle \left\langle \text{\textperiodcentered}%
,\text{\textperiodcentered}\right\rangle \right\rangle $ the duality map
between $H_{\mu}^{1}(\mathbb{R}^{+};V)$ and its dual space. With the help of
\eqref{eq3-14} and the continuity property of function $a$, we proceed
likewise as in the proof of \cite[p.344]{G1} (see also \cite[Step 3. Appendix
A]{X2}) to finish these arguments.

\textbf{Step 4. (Continuity of solution)} By means of similar arguments as in
\cite[Step 4. Appendix A]{X2}, it is immediate to see that for every
$\omega\in\Omega$, $(\frac{dv}{dt},\frac{d\eta}{dt})$ fulfills
\[%
\begin{split}
&  \frac{dv}{dt}\in L^{2}(\tau,T;V^{\ast})+L^{q}(\tau,T;L^{q}(\mathcal{O}%
));\\[1ex]
&  \frac{d\eta}{dt}\in L^{2}(\tau,T;H_{\mu}^{-1}(\mathbb{R}^{+};V)).
\end{split}
\]
Using a slightly modified version of \cite[Lemma III.1.2]{T1}, together with
\eqref{eq3-14}, we deduce that $v\in C([\tau,T];H)$. As for the second
component, by applying the same arguments as for the theorem in \cite[Section
2]{G1}, we obtain that $\eta\in C([\tau,T];L_{\mu}^{2}(\mathbb{R}^{+};V))$.
Thus, $(v(\tau),\eta^{\tau})$ makes sense and the equality $(v(\tau
),\eta^{\tau})=(v_{0},\eta_{0})$ follows from the fact that $(P_{n}v_{0}%
,Q_{n}\eta_{0})$ converges to $(v_{0},\eta_{0})$ strongly for each $\omega
\in\Omega$.

\textbf{Step 5. (Continuity with respect to initial value and uniqueness)} Let
$y_{1}=(v_{1},\eta_{1})$ and $y_{2}=(v_{2},\eta_{2})$ be two solutions of
\eqref{eq1-5} with initial data $y_{10}$ and $y_{20}$, respectively. Due to
the fact that $v\in C([\tau,T];H)$, for every $\omega\in\Omega$, there exists
a bounded set $S\subset H$ such that $v_{i}(t)\in S$ for all $t\in\lbrack
\tau,T]$ and $i=1,2$. In addition, taking into account that $l\in
\mathcal{L}(H;\mathbb{R})$, $\phi\in V\cap H^{2}(\mathcal{O})\cap
L^{2p}(\mathcal{O})$ and $z_{\ast}(\theta_{t}\omega)$ is uniformly bounded for
each $\omega\in\Omega$ and all $t\in\lbrack\tau,T]$, we have $\{l(v_{i}%
(t)+l(\phi)z_{\ast}(\theta_{t}\omega))\}_{t\in\lbrack\tau,T]}\subset
\lbrack-R,R]$ for $i=1,2$ and some $R>0$. Hence, let $\bar{y}=y_{1}%
-y_{2}=(\bar{v},\bar{\eta})=(v_{1}-v_{2},\eta_{1}-\eta_{2})$ and $\bar{y}%
_{0}=y_{10}-y_{20}$, we have
\[%
\begin{split}
\frac{1}{2}\frac{d}{dt}|\bar{v}(t)|^{2}+\frac{1}{2}\frac{d}{dt}\Vert\bar{\eta
}^{t}\Vert_{\mu}^{2}  &  \leq|a(l(v_{1})+l(\phi)z_{\ast}(\theta_{t}%
\omega))-a(l(v_{2})+l(\phi)z_{\ast}(\theta_{t}\omega))||z_{\ast}(\theta
_{t}\omega)|\Vert\phi\Vert\Vert\bar{v}\Vert\\[0.8ex]
&  ~~-(((\bar{\eta}^{t})^{\prime},\bar{\eta}^{t}))_{\mu}+|a(l(v_{1}%
)+l(\phi)z_{\ast}(\theta_{t}\omega))-a(l(v_{2})+l(\phi)z_{\ast}(\theta
_{t}\omega))|\Vert v_{2}\Vert\Vert\bar{v}\Vert\\[0.8ex]
&  ~~-a(l(v_{1})+l(\phi)z_{\ast}(\theta_{t}\omega))\Vert\bar{v}\Vert
^{2}-(f(v_{1}+\phi z_{\ast}(\theta_{t}\omega))-f(v_{2}+\phi z_{\ast}%
(\theta_{t}\omega)),\bar{v})_{L^{p,q}},
\end{split}
\]
where $(\cdot,\cdot)_{L^{p,q}}$ is the duality between $L^{2p}$ and $L^{q}$.
Observe that it follows from integration by parts and the fact $\mu^{\prime
}\leq0$ that,
\[
2(((\bar{\eta}^{t})^{\prime},\bar{\eta}^{t}))_{\mu}=-\lim_{s\rightarrow0}%
\mu(s)|\nabla\bar{\eta}^{t}(s)|^{2}-\int_{0}^{\infty}\mu^{\prime}%
(s)|\nabla\bar{\eta}^{t}(s)|^{2}ds\geq0.
\]
All these operations are formal but can be justified using mollifiers (see
\cite[Section 2]{G1}). Applying the Poincar\'{e} and Young inequalities,
together with \eqref{eq1-2} and the above results, we obtain
\begin{equation}%
\begin{split}
\frac{d}{dt}\Vert\bar{y}\Vert_{\mathcal{H}}^{2}  &  \leq-2m\Vert\bar{v}%
\Vert^{2}+2L_{a}(R)|l||\bar{v}||z_{\ast}(\theta_{t}\omega)|\Vert\phi\Vert
\Vert\bar{v}\Vert+2L_{a}(R)|l||\bar{v}|\Vert v_{2}\Vert\Vert\bar{v}%
\Vert\\[0.8ex]
&  ~~-2(f(v_{1}+\phi z_{\ast}(\theta_{t}\omega))-f(v_{2}+\phi z_{\ast}%
(\theta_{t}\omega)),\bar{v})_{L^{p,q}}\\[0.8ex]
&  \leq-2m\Vert\bar{v}\Vert^{2}+2m\Vert\bar{v}\Vert^{2}+\frac{1}{m}L_{a}%
^{2}(R)|l|^{2}|\bar{v}|^{2}\left(  \Vert\phi\Vert^{2}|z_{\ast}(\theta
_{t}\omega)|^{2}+\Vert v_{2}\Vert^{2}\right) \\[0.8ex]
&  ~~-2(f(v_{1}+\phi z_{\ast}(\theta_{t}\omega))-f(v_{2}+\phi z_{\ast}%
(\theta_{t}\omega)),\bar{v})_{L^{p,q}}.
\end{split}
\label{eq3-16}%
\end{equation}
{Since $f$ is a polynomial of odd degree with positive leading coefficient, we
find that there exists a positive constant $\sigma$, such that}
\begin{equation}
{f^{\prime}(s)\geq-\frac{\sigma}{2},\qquad\forall s\in\mathbb{R}.}
\label{eq630}%
\end{equation}
{With help of the mean value theorem, we deduce}
\begin{equation}%
\begin{split}
&  ~\quad-2(f(v_{1}+\phi z_{\ast}(\theta_{t}\omega))-f(v_{2}+\phi z_{\ast
}(\theta_{t}\omega)),\bar{v})_{L^{p,q}}\\[0.8ex]
&  =-2\int_{\mathcal{O}}(f(v_{1}+\phi z_{\ast}(\theta_{t}\omega))-f(v_{2}+\phi
z_{\ast}(\theta_{t}\omega)))\bar{v}dx\\[0.8ex]
&  =-2\int_{\mathcal{O}}f^{\prime}(s_{x})|\bar{v}|^{2}dx\leq\sigma|\bar
{v}|^{2}\leq\sigma\Vert\bar{y}\Vert_{\mathcal{H}}^{2},
\end{split}
\label{eq3-17}%
\end{equation}
where $s_{x}$ is an intermediate point between $v_{1}(x)+\phi(x)z_{\ast
}(\theta_{t}\omega)$ and $v_{2}(x)+\phi(x)z_{\ast}(\theta_{t}\omega)$.

Subsequently, \eqref{eq3-16}-\eqref{eq3-17} imply that
\[
\frac{d}{dt}\Vert\bar{y}\Vert_{\mathcal{H}}^{2}\leq\left(  \frac{1}{m}%
L_{a}^{2}(R)|l|^{2}\big(\Vert\phi\Vert^{2}|z_{\ast}(\theta_{t}\omega
)|^{2}+\Vert v_{2}\Vert^{2}\big)+{\sigma}\right)  \Vert\bar{y}\Vert
_{\mathcal{H}}^{2}.
\]
The uniqueness and continuous dependence on initial data of solutions to
problem \eqref{eq1-5} follow from the Gronwall lemma. Till now, the proof of
the first assertion is finished.

$(ii)$ \textbf{(Further regularity)} We are going to study further regularity
of $(v,\eta)$. To this end, for every $\omega\in\Omega$ and $\tau\in
\mathbb{R}$, let us first consider the operator $\mathcal{I}_{\tau,\omega
}:L_{V\cap H^{2}(\mathcal{O})}^{2}\rightarrow L_{\mu}^{2}(\mathbb{R}^{+};V\cap
H^{2}(\mathcal{O}))$ defined by
\[
(\mathcal{I}_{\tau,\omega}{\varphi})(s)=\int_{-s}^{0}\varphi(r)dr-\int%
_{-s}^{0}z_{\ast}(\theta_{r+\tau}\omega)\phi dr.
\]
Thus, similar to \cite{X2}, we know that the operator $\mathcal{I}%
_{\tau,\omega}$ introduced above is a continuous mapping. Particularly, there
exists a positive constant $K_{\mu}$, which is the same as in Corollary
\ref{cor3-2}, such that for any $\varphi\in L_{V\cap H^{2}(\mathcal{O})}^{2}$
and $\phi\in V\cap H^{2}(\mathcal{O})\cap L^{2p}(\mathcal{O})$, it holds
\[
\Vert\mathcal{I}_{\tau,\omega}\varphi\Vert_{L_{\mu}^{2}(\mathbb{R}^{+};V\cap
H^{2}(\mathcal{O}))}^{2}\leq2K_{\mu}\left(  \Vert\varphi\Vert_{L_{V\cap
H^{2}(\mathcal{O})}^{2}}^{2}+\Vert z_{\ast}(\theta_{\cdot+\tau}\omega
)\phi\Vert_{L_{V\cap H^{2}(\mathcal{O})}^{2}}^{2}\right)  .
\]
Next, multiplying $\eqref{eq1-5}_{1}$ by $-\Delta v$ with respect to the inner
product of $H$, the Laplacian of $\eqref{eq1-5}_{2}$ by $\Delta\eta^{t}$ with
respect to the inner product of $L_{\mu}^{2}(\mathbb{R}^{+};H)$ and adding the
two terms, for every $\omega\in\Omega$, it follows from \eqref{eq1-2} that,%
\begin{equation}%
\begin{split}
\frac{d}{dt}\Vert v\Vert^{2}  &  +\frac{d}{dt}\Vert\Delta\eta^{t}\Vert_{\mu
}^{2}+2m|\Delta v|^{2}+2((((\eta^{t})^{\prime},\eta^{t})))_{\mu}\leq
2M|z_{\ast}(\theta_{t}\omega)||\Delta\phi||\Delta v|\\[0.8ex]
&  ~~+2(f(v+\phi z_{\ast}(\theta_{t}\omega)),\Delta v)+2|z_{\ast}(\theta
_{t}\omega)|\Vert\phi\Vert\Vert v\Vert+2(z_{k}^{\phi}(t,\omega),-\Delta
v)+{2}(h,-\Delta v).
\end{split}
\label{eq3-18}%
\end{equation}
Applying the same arguments as in \cite[Appendix A]{X2}, we know that%
\begin{equation}
2((((\eta^{t})^{\prime},\eta^{t})))_{\mu}=-2\int_{0}^{\infty}\mu^{\prime
}(s)|\Delta\eta^{t}(s)|^{2}ds>0. \label{eq3-19}%
\end{equation}
By means of the Young inequality, we infer that
\begin{equation}
2M|z_{\ast}(\theta_{t}\omega)||\Delta\phi||\Delta v|\leq\frac{2M^{2}}%
{m}|z_{\ast}(\theta_{t}\omega)|^{2}|\Delta\phi|^{2}+\frac{m}{2}|\Delta v|^{2},
\label{eq3-20}%
\end{equation}
and
\begin{equation}
2|z_{\ast}(\theta_{t}\omega)|\Vert\phi\Vert\Vert v\Vert\leq|z_{\ast}%
(\theta_{t}\omega)|^{2}\Vert\phi\Vert^{2}+\Vert v\Vert^{2}. \label{eq3-21}%
\end{equation}
{Next, taking into account of \eqref{eq630} and assumption $\Delta\phi\in
L^{2p}(\mathcal{O})$}, together with the Young inequality and the Green
formula, it yields
\begin{equation}%
\begin{split}
&  ~\quad2(f(v+\phi z_{\ast}(\theta_{t}\omega)),\Delta v)\\[0.8ex]
&  =2(f(v+\phi z_{\ast}(\theta_{t}\omega)),\Delta(v+\phi z_{\ast}(\theta
_{t}\omega)))-2(f(v+\phi z_{\ast}(\theta_{t}\omega)),z_{\ast}(\theta_{t}%
\omega)\Delta\phi)\\[0.8ex]
&  \leq2\int_{\mathcal{O}}f_{2p-1}(\Delta v+z_{\ast}(\theta_{t}\omega
)\Delta\phi)dx+2|z_{\ast}(\theta_{t}\omega)|\int_{\mathcal{O}}|f(v+\phi
z_{\ast}(\theta_{t}\omega))||\Delta\phi|dx\\[0.8ex]
&  ~~-2\int_{\mathcal{O}}\nabla(v+\phi z_{\ast}(\theta_{t}\omega))\cdot
f^{\prime}(v+z_{\ast}(\theta_{t}\omega)\phi)\nabla(v+\phi z_{\ast}(\theta
_{t}\omega))dx\\[0.8ex]
&  \leq\frac{4}{m}f_{2p-1}^{2}|\mathcal{O}|+\frac{m}{2}|\Delta v|^{2}+\frac
{m}{2}|z_{\ast}(\theta_{t}\omega)|^{2}|\Delta\phi|^{2}+2\sigma\Vert v\Vert
^{2}+2\sigma|z_{\ast}(\theta_{t}\omega)|^{2}\Vert\phi\Vert^{2}\\[0.8ex]
&  ~~+\frac{2|z_{\ast}(\theta_{t}\omega)|}{q}\Vert f(v+\phi z_{\ast}%
(\theta_{t}\omega))\Vert_{q}^{q}+\frac{|z_{\ast}(\theta_{t}\omega)|}{p}%
\Vert\Delta\phi\Vert_{2p}^{2p},
\end{split}
\label{eq3-22}%
\end{equation}
where $f(v+\phi z_{\ast}(\theta_{t}\omega))\in L^{q}(\mathcal{O})$ since $v\in
L^{2p}(\mathcal{O})$ and $\phi\in L^{2p}(\mathcal{O})$. In the end, it follows
from the Young inequality that,
\begin{equation}
2(z_{k}^{\phi}(\theta_{t}\omega),-\Delta v)\leq\frac{C_{2}(\theta_{t}%
\omega)^{2}}{m}|\Delta\phi|^{2}+\frac{m}{4}|\Delta v|^{2}, \label{eq3-23}%
\end{equation}%
\[
{2}(h,-\Delta v)\leq\frac{{4}}{m}\left\vert h\right\vert ^{2}+\frac{m}%
{4}\left\vert \Delta v\right\vert ^{2},
\]
where $C_{2}$ is the same as in \eqref{eq3-10}.

Substituting \eqref{eq3-19}-\eqref{eq3-23} into \eqref{eq3-18}, we have
\[%
\begin{split}
\frac{d}{dt}\Vert v\Vert^{2}+\frac{d}{dt}\Vert\Delta\eta^{t}\Vert_{\mu}%
^{2}+\frac{m}{2}|\Delta v|^{2}  &  \leq\left(  \frac{2M^{2}}{m}+\frac{m}%
{2}\right)  |z_{\ast}(\theta_{t}\omega)|^{2}|\Delta\phi|^{2}+(2\sigma
+1)|z_{\ast}(\theta_{t}\omega)|^{2}\Vert\phi\Vert^{2}\\[0.8ex]
&  ~~+\frac{4}{m}f_{2p-1}^{2}|\mathcal{O}|+\left(  2\sigma+1\right)  \Vert
v\Vert^{2}+\frac{C_{2}(\theta_{t}\omega)^{2}}{m}|\Delta\phi|^{2}\\[0.8ex]
&  ~~+\frac{2|z_{\ast}(\theta_{t}\omega)|}{q}\Vert f(v+\phi z_{\ast}%
(\theta_{t}\omega))\Vert_{q}^{q}+\frac{|z_{\ast}(\theta_{t}\omega)|}{p}%
\Vert\Delta\phi\Vert_{2p}^{2p}+\frac{4}{m}\left\vert h\right\vert ^{2}.
\end{split}
\]
Denote
\[%
\begin{split}
\Theta_{2}(t,\omega):  &  =\left(  \frac{2M^{2}}{m}+\frac{m}{2}\right)
|z_{\ast}(\omega)|^{2}|\Delta\phi|^{2}+(2\sigma+1)|z_{\ast}(\omega)|^{2}%
\Vert\phi\Vert^{2}+\frac{4}{m}f_{2p-1}^{2}|\mathcal{O}|+\frac{C_{2}%
(\omega)^{2}}{m}|\Delta\phi|^{2}\\[0.8ex]
&  ~~+\frac{2|z_{\ast}(\omega)|}{q}\Vert f(v(t)+\phi z_{\ast}(\omega
))\Vert_{q}^{q}+\frac{|z_{\ast}(\omega)|}{p}\Vert\Delta\phi\Vert_{2p}%
^{2p}{+\frac{4}{m}|h|^{2}}\in L^{1}(\tau,T).
\end{split}
\]
Then, for every $\omega\in\Omega$ and $t\in(\tau,T]$, we obtain
\begin{equation}
\frac{d}{dt}\Vert y\Vert_{\mathcal{V}}^{2}+\frac{m}{2}|\Delta v|^{2}\leq
\Theta_{2}(t,\theta_{t}\omega)+\left(  2\sigma+1\right)  \Vert v\Vert^{2}.
\end{equation}
By the continuity of $z_{\ast}(\theta_{t}\omega)$ on $(\tau,T]$ and
integrating the above inequality between $\tau$ and $t$ with $\tau\leq t\leq
T$, we have
\[
\Vert y(t)\Vert_{\mathcal{V}}^{2}+\frac{m}{2}\int_{\tau}^{t}|\Delta
v(s)|^{2}ds\leq\Vert y_{0}\Vert_{\mathcal{V}}^{2}+\int_{\tau}^{t}\Theta
_{2}(s,\theta_{s}\omega)ds+(2\sigma+1)\int_{\tau}^{t}\Vert v(s)\Vert^{2}ds.
\]
Thus, we conclude that
\[%
\begin{split}
&  v\in L^{\infty}(\tau,T;V)\cap L^{2}(\tau,T;V\cap H^{2}(\mathcal{O}));\\
&  \eta\in L^{\infty}(\tau,T;L_{\mu}^{2}(\mathbb{R}^{+};V\cap H^{2}%
(\mathcal{O}))).
\end{split}
\]
Furthermore, the continuity of $v$ follows again using a slightly modified
version of \cite[Lemma III.1.2.]{T1} and the continuity of $\eta$ can be
proved mimicking the idea of the proof of Step 4 of $(i)$, with $V\cap
H^{2}(\mathcal{O})$ in place of $V$. The proof of this theorem is complete.
$\Box$

\begin{lemma}
\label{lemEquivInt}Let conditions (h1)-(h2) hold. If $u\in L_{V}^{2}$, then
$\eta(s)=\int_{-s}^{0}u\left(  r\right)  dr$ belongs to $L_{\mu}%
^{2}(\mathbb{R}^{+};V)$ and%
\begin{equation}
\int_{0}^{\infty}\mu(s)\Delta\eta(s)ds=\int_{-\infty}^{0}k(-s)\Delta u(s)ds.
\label{EquivIntegrals}%
\end{equation}

\end{lemma}

\textbf{Proof. }The fact that $\eta\in L_{\mu}^{2}(\mathbb{R}^{+};V)$ is given
by Lemma \ref{lem3-1}. From the arguments in \cite[pp-174-175]{Gajewski}, it
follows the existence of a sequence of functions $u_{n}\left(
\text{\textperiodcentered}\right)  \in C^{1}((-\infty,0],V)\cap L_{V}^{2}$
such that
\[
u_{n}\rightarrow u\text{ in }L_{V}^{2}.
\]
First, we will show that $u_{n},\eta_{n}$, where $\eta_{n}(s)=\int_{-s}%
^{0}u_{n}\left(  r\right)  dr$, satisfy (\ref{EquivIntegrals}). For any $w\in
V$, we have%
\begin{align*}
\left\langle \int_{0}^{\infty}\mu(s)\Delta\eta_{n}(s)ds,w\right\rangle  &
=\int_{0}^{\infty}\mu(s)\left\langle \Delta\eta_{n}(s),w\right\rangle
ds=\int_{0}^{\infty}k^{\prime}\left(  s\right)  \left(  \nabla\eta
_{n}(s),\nabla w\right)  ds\\
&  =\int_{0}^{\infty}k^{\prime}\left(  s\right)  \left(  \nabla\int_{-s}%
^{0}u_{n}(r)dr,\nabla w\right)  ds=\int_{0}^{\infty}k^{\prime}\left(
s\right)  \int_{-s}^{0}\left(  \nabla u_{n}(r),\nabla w\right)  drds\\
&  =-\int_{0}^{\infty}k(s)\left(  \nabla u_{n}(-s),\nabla w\right)
ds+\lim_{s\rightarrow\infty}\ k(s)\int_{-s}^{0}\left(  \nabla u_{n}(r),\nabla
w\right)  dr\\
&  ~~ -\lim_{s\rightarrow0}\ k(s)\int_{-s}^{0}\left(  \nabla u_{n}(r),\nabla
w\right)  dr.
\end{align*}
Let us check that the last two limits of the above equality are equal to $0$.
By Remark \ref{rem2-1}, we derive%
\[
k(s)e^{\gamma s}\leq\frac{\mu(1)}{\varpi}e^{\varpi}e^{(\gamma-\varpi
)s},~~~\text{ for any }s\geq1.
\]
Hence, $\gamma<\varpi$ implies%
\begin{align*}
\left\vert k(s)\int_{-s}^{0}\left(  \nabla u_{n}(r),\nabla w\right)
dr\right\vert  &  \leq k(s)e^{\gamma s}\left\Vert w\right\Vert \int_{-s}%
^{0}e^{\gamma r}\left\Vert u_{n}(r)\right\Vert dr\\
&  \leq\frac{k(s)e^{\gamma s}\left\Vert w\right\Vert }{2}\left(  \int%
_{-\infty}^{0}e^{\gamma r}\left\Vert u_{n}(r)\right\Vert ^{2}dr+\frac
{1}{\gamma}\right)  \leq C_{1}e^{(\gamma-\varpi)s}\underset{s\rightarrow
\infty}{\rightarrow}0.
\end{align*}
Also, from $k\left(  s\right)  \underset{s\rightarrow0}{\rightarrow}\int%
_{0}^{\infty}\mu\left(  r\right)  dr$ and $u_{n}\in L_{V}^{2}$, it follows
that the second limit is $0$ as well. Hence,%
\[
\left\langle \int_{0}^{\infty}\mu(s)\Delta\eta_{n}(s)ds,w\right\rangle
=-\int_{0}^{\infty}k(s)\left(  \nabla u_{n}(-s),\nabla w\right)
ds=\left\langle \int_{-\infty}^{0}k(-s)\Delta u_{n}(s)ds,w\right\rangle ,
\]
this proves (\ref{EquivIntegrals}) for $u_{n}$.

Furthermore, for any $w\in V$, we infer%
\begin{align*}
&  \left\vert \left\langle \int_{-\infty}^{0}k(-s)(\Delta u_{n}(s)-\Delta
u(s))ds,w\right\rangle \right\vert \\
&  =\left\vert \int_{-\infty}^{0}k(-s)\left\langle \Delta u_{n}(s)-\Delta
u(s),w\right\rangle ds\right\vert \\
&  \leq\left\Vert w\right\Vert \left(  C_{2}\int_{-1}^{0}\left\Vert
u_{n}(s)-u(s)\right\Vert ds+\frac{\mu(1)}{\varpi}e^{\varpi}\int_{-\infty}%
^{-1}e^{\varpi s}\left\Vert u_{n}(s)-u(s)\right\Vert ds\right) \\
&  \leq C_{3}\left(  \left(  \int_{-1}^{0}\left\Vert u_{n}(s)-u(s)\right\Vert
^{2}ds\right)  ^{\frac{1}{2}}+\left(  \int_{-\infty}^{-1}e^{\gamma
s}\left\Vert u_{n}(s)-u(s)\right\Vert ^{2}ds\right)  ^{\frac{1}{2}}\right)
\underset{n\rightarrow\infty}{\rightarrow}0,
\end{align*}
and Lemma \ref{lem3-1} implies%
\begin{align*}
&  \left\vert \left\langle \int_{0}^{\infty}\mu(s)(\Delta\eta_{n}%
(s)-\Delta\eta(s))ds,w\right\rangle \right\vert \\
&  =\left\vert \int_{0}^{\infty}\mu(s)\left\langle \Delta\eta_{n}%
(s)-\Delta\eta(s),w\right\rangle ds\right\vert \leq\left\Vert w\right\Vert
\int_{0}^{\infty}\mu(s)\left\Vert \eta_{n}(s)-\eta(s)\right\Vert ds\\
&  \leq\left\Vert w\right\Vert \left(  \int_{0}^{\infty}\mu(s)ds\right)
^{\frac{1}{2}}\left(  \int_{0}^{\infty}\mu(s)\left\Vert \eta_{n}%
(s)-\eta(s)\right\Vert ^{2}ds\right)  ^{\frac{1}{2}}\leq C_{4}\left\Vert
u_{n}-u\right\Vert _{L_{V}^{2}}\underset{n\rightarrow\infty}{\rightarrow}0.
\end{align*}
By these convergences we deduce (\ref{EquivIntegrals}). The proof of this
lemma is complete. $\Box$

Lemma \ref{lemEquivInt} implies that the solution given in Theorem
\ref{thm3-3} is in fact the unique weak solution to problem (\ref{eq1-4}).

\begin{corollary}
\label{CorExistSol}Assume that \eqref{eq1-2}, \eqref{eq2-1} and $(h_{1}%
)$-$(h_{2})$ hold, and that $\phi\in V\cap H^{2}(\mathcal{O})\cap
L^{2p}(\mathcal{O})$ is such that $\Delta\phi\in L^{2p}(\mathcal{O})$. Let
{$h\in H$} and $a$ be a locally Lipschitz function. If for fixed $\tau
\in\mathbb{R}$ and $\omega\in\Omega$, the function $(v,\eta)$ {is} the unique
weak solution to problem \eqref{eq1-5} corresponding to the initial values
$v_{0}\in H$ and $\varphi\in L_{V}^{2}$, then $v$ is the unique weak solution
to problem (\ref{eq1-4}).
%%
%%\red{Then} for every $\tau
%%\in\mathbb{R}$ and $\omega\in\Omega$, the function $(v,\eta)$ \red{is} be the unique
%%weak solution to problem \eqref{eq1-5}, \red{and} $v$ is the unique weak solution
%%to problem (\ref{eq1-4}).

\end{corollary}

Now by the transform \eqref{eq1-3}, we derive the well-posedness of problem \eqref{eq1-1}.

\begin{theorem}
\label{thm3-5} Assume that \eqref{eq1-2}, \eqref{eq2-1} and $(h_{1})$%
-$(h_{2})$ hold, and that $\phi\in V\cap H^{2}(\mathcal{O})\cap L^{2p}%
(\mathcal{O})$ is such that $\Delta\phi\in L^{2p}(\mathcal{O})$. Let {$h\in
H$} and $a$ be a locally Lipschitz function. Then, for every $\tau
\in\mathbb{R}$ and $\omega\in\Omega$, it holds:

\begin{enumerate}
\item[$(i)$] For any initial vale $u_{0}\in H$ and initial function
$\varphi\in L_{V}^{2}$, there exists a unique solution $u$ to problem
\eqref{eq1-1} in the weak sense, fulfilling
\[
u\in L^{\infty}(\tau,T;H)\cap L^{2}(\tau,T;V)\cap L^{2p}(\tau,T;L^{2p}%
(\mathcal{O})),\qquad\forall T>\tau.
\]
Furthermore, the solution $u$ of \eqref{eq1-1} is continuous with respect to
the initial values {$(u_{0},\varphi)$} for all $t\in[\tau,T]$ in $H$;

\item[$(ii)$] For any initial value $u_{0}\in V$ and initial function
$\varphi\in L_{V\cap H^{2}(\mathcal{O})}^{2}$, the unique solution $u$ to
problem \eqref{eq1-1} satisfies,
\[
u\in L^{\infty}(\tau,T;V)\cap L^{2}(\tau,T;V\cap H^{2}(\mathcal{O}%
)),\qquad\forall T>\tau.
\]
In addition, the solution $u$ of \eqref{eq1-1} is continuous with respect to
the initial values {$(u_{0},\varphi)$} for all $t\in\lbrack\tau,T]$ in $V $.
\end{enumerate}
\end{theorem}

\begin{remark}
The proof of Theorem \ref{thm3-3} is correct for a general function $f\in
C^{1}(\mathbb{R})$ satisfying (\ref{fcondition}) and (\ref{eq630}). The same
applies to the results in Sections \ref{s4}-\ref{s5}.
\end{remark}

\section{Existence of random attractors}

\label{s4}

This section is devoted to studying the long time behavior of \eqref{eq1-1} in
the natural phase space,
\[
X=H\times L_{V}^{2},
\]
endowed with the norm
\[
\Vert(w_{1},w_{2})\Vert_{X}^{2}=|w_{1}|^{2}+\Vert w_{2}\Vert_{L_{V}^{2}}^{2}.
\]
%It is worth emphasizing that we will take $\tau\equiv 0$ throughout this section since problem \eqref{eq1-1} is autonomous.
{It is worth emphasizing that we will take $\tau=0$ in this section since
problem \eqref{eq1-1} is autonomous.} Taking into account the results in the
previous section, problem \eqref{eq1-1} generates a random dynamical system in
$X$. Let us denote by $u(\cdot;0,\omega,(u_{0},\varphi))$ the unique solution
to \eqref{eq1-1}. Then, the random dynamical system (RDS) generated by
\eqref{eq1-1}, denoted by $\Xi:\mathbb{R}^{+}\times\Omega\times X\rightarrow
X,$ is defined, for every $t\in\mathbb{R}^{+}$, $\omega\in\Omega$ and
$(u_{0},\varphi)\in X$, as
\[
\Xi(t,\omega,(u_{0},\varphi))=(u(t;0,\omega,(u_{0},\varphi)),u_{t}%
(\cdot;0,\omega,(u_{0},\varphi))).
\]
Moreover, problem \eqref{eq1-5} also generates a random dynamical system
$\Phi$ on the phase space $H\times L_{\mu}^{2}(\mathbb{R}^{+};V)$, which is
defined, for every $t\in\mathbb{R}^{+}$, $\omega\in\Omega$ and $(v_{0}%
,\eta_{0})\in H\times L_{\mu}^{2}(\mathbb{R}^{+};V)$, by
\[
\Phi(t,\omega,(v_{0},\eta_{0}))=(v(\cdot;0,\omega,(v_{0},\eta_{0}%
)),\eta^{\cdot}(\cdot;0,\omega,(v_{0},\eta_{0}))),
\]
where the right-hand side of the above equality denotes the solution to
\eqref{eq1-5} for $\tau=0$, the initial values $(v_{0},\eta_{0})\in H\times
L_{\mu}^{2}(\mathbb{R}^{+};V)$ and $\eta_{0}$ is given in \eqref{eq1-9}.
Thanks to Dafermos' transformation, we can obtain a random dynamical system
$\Psi:\mathbb{R}^{+}\times\Omega\times X\rightarrow X$ generated by
\eqref{eq1-4} which is given, for every $t\in\mathbb{R}^{+}$, $\omega\in
\Omega$ and $(v_{0},\psi)\in X$, by
\[
\Psi(t,\omega,(v_{0},\psi))=(v(t;0,\omega,(v_{0},(\mathcal{J}\psi
))),v_{t}(\cdot;0,\omega,(v_{0},(\mathcal{J}\psi)))).
\]
Then, on account of the random transformation \eqref{eq1-3}, for
$(u_{0},\varphi)\in X$, we deduce
\begin{align}
&  \quad\Xi(t,\omega,(u_{0},\varphi))\label{cocycle}\\
&  =(u(t;0,\omega,(u_{0},\varphi)),u_{t}(\cdot;0,\omega,(u_{0},\varphi
)))\nonumber\\
&  =(v(t;0,\omega,(u_{0}-\phi z_{\ast}(\omega),(\mathcal{J}_{\omega,0}%
{\varphi})))+\phi z_{\ast}(\theta_{t}\omega),v_{t}(\cdot;0,\omega,(u_{0}-\phi
z_{\ast}(\omega),(\mathcal{J}_{\omega,0}{\varphi})))+\phi z_{\ast}%
(\theta_{t+\cdot}\omega))\nonumber\\
&  =\Psi(t,\omega,(u_{0}-\phi z_{\ast}(\omega),{\varphi_{v}}))+(\phi z_{\ast
}(\theta_{t}\omega),\phi z_{\ast}(\theta_{t+\cdot}\omega)).\nonumber
\end{align}
It is straightforward to check that the cocycles $\Xi$ and $\Psi$ are
conjugated. Indeed, consider the mapping $T:\Omega\times X\rightarrow X$
defined by,
\[
T(\omega,(u_{0},\varphi))=(u_{0}-\phi z_{\ast}(\omega),\varphi-\phi z_{\ast
}(\theta_{\cdot}\omega)).
\]
Then, it holds that
\[
T^{-1}(\omega,(u_{0},\varphi))=(u_{0}+\phi z_{\ast}(\omega),\varphi+\phi
z_{\ast}(\theta_{\cdot}\omega)).
\]
In addition, by \eqref{cocycle}, it is clear that
\begin{equation}%
\begin{split}
\Xi(t,\omega,(u_{0},\varphi))  &  =\Psi(t,\omega,(u_{0}-\phi z_{\ast}%
(\omega),\varphi-\phi z_{\ast}(\theta_{\cdot}\omega)))+(\phi z_{\ast}%
(\theta_{t}\omega),\phi z_{\ast}(\theta_{t+\cdot}\omega))\\[0.6ex]
&  =T^{-1}(\theta_{t}\omega,\Psi(t,\omega,T(\omega,(u_{0},\varphi))).
\end{split}
\label{conjugation}%
\end{equation}

Let $D=\{D(\omega):\omega\in\Omega\}$ be a family of bounded nonempty subsets
of $X$. Such a family $D$ is called tempered if for every $c>0$ and $\omega
\in\Omega$,
\[
\lim_{t\rightarrow\infty}e^{-ct}\Vert D(\theta_{-t}\omega)\Vert=0,
\]
where the norm $\Vert D\Vert$ of a set $D$ in $X$ is defined by $\Vert
D\Vert=\sup_{u\in D}\Vert u\Vert_{X}$. From now on, we will use $\mathcal{D}$
to denote the collection of all tempered families of bounded nonempty subsets
of $X$:
\[
\mathcal{D}=\{D=\{D(\omega):\omega\in\Omega\}:D~\mbox{is tempered in}~X\}.
\]
This family will be adopted to prove the existence of random pullback
attractors for the RDS $\Xi$. Notice that, for $D\in\mathcal{D}$, the set
$\widetilde{D}$ whose fibers are given by,
\[
\widetilde{D}(\omega)=\{(u_{0}-\phi z_{\ast}(\omega),\varphi-\phi z_{\ast
}(\theta_{\cdot}\omega)):(u_{0},\varphi)\in D(\omega)\},
\]
also belongs to $\mathcal{D}$ thanks to the arguments in the proof of
Corollary \ref{cor3-2} and the properties of the random variable $z_{\ast
}(\omega)$ (cf. \eqref{zproperty}).
%which means that $z_*(\theta_t\omega)$ has subexponential growth as $t\to\pm\infty$, for each $\omega\in\Omega$.

%%We will first prove the cocycle $\Phi$ is well defined.
%%
%%
%%
%%\begin{lemma}\label{lem4-1}
%%Under assumptions of Theorem \ref{thm3-3}, if $(v_0,\varphi-\phi z_*(\theta_t\omega))\in X$, then for all  $t\in\mathbb{R}^+$ and $\omega\in\Omega$, $\Phi(t,\omega,(v_0,\varphi-\phi z_*(\theta_t\omega)))\in X$.
%%\end{lemma}
%%%=========detail of proof ======
%%{\bf Proof.}  The results follows the lines of the proof \cite[Lemma 3.6]{X2}, we omit the details here. $\Box$

%Let $(v_0,\rho)\in X$, it follows from Theorem
%\ref{thm3-3} that the first component $v(t)$ belongs to $H$, thus it only remains to prove the segment of solution $v_t(\cdot)$ belongs to $L_V^2$. Indeed,
%\begin{equation*}
%\begin{split}
%\int_{-\infty}^0e^{\gamma s}\|v_t(s)\|^2ds&=\int_{-\infty}^0e^{\gamma s}\|v(t+s)\|^2ds=\int_{-\infty}^te^{\gamma(\sigma-t)}\|v(\sigma)\|^2d\sigma\\
%&=e^{-\gamma t}\int_{-\infty}^te^{\gamma\sigma}\|v(\sigma)\|^2d\sigma\\
%&=e^{-\gamma t}\int_{-\infty}^0e^{\gamma\sigma}\|\rho(\sigma)\|^2d\sigma+%\int_0^te^{\gamma(\sigma-t)}\|v(\sigma)\|^2d\sigma<\infty,
%\end{split}
%\end{equation*}
%where the above estimates holds since $\rho\in L_V^2$ and $v\in L^2(0,T;V)$ for every $\omega\in\Omega$ and $T>0$.

%The proof of this lemma is finished. $\Box$

\begin{lemma}
\label{absorbing} Under assumptions of Theorem \ref{thm3-5}, there exists
$B\in\mathcal{D}$ which is $\mathcal{D}$-pullback absorbing for the RDS $\Xi$.
In other words, for any given $\omega\in\Omega$ and $D\in\mathcal{D}$, there
exists $t_{0}:=t_{0}(\omega,D)\geq0$, such that
\[
\Xi(t,\theta_{-t}\omega,D(\theta_{-t}\omega))\subset B(\omega),\ \quad
\text{for all}\ t\geq t_{0}(\omega,D),
\]
where $B(\omega)$ is the ball in $X$ centered at $0$ with radius {$\rho
(\omega)$} and
\[
\rho^{2}(\omega)=1+{2}K_{2}\int_{-\infty}^{0} e^{\gamma s}{\Theta_{1}%
(\theta_{s}\omega)}ds+{2}|\phi|^{2}|z_{*}(\omega)|^{2}+{2}\|\phi z_{*}%
(\theta_{\cdot}\omega)\|^{2}_{L_{V}^{2}},
\]
where $\Theta_{1}(\omega)$ is given in \eqref{theta1} and $K_{2}>0$ is a constant.
\end{lemma}

\textbf{Proof.} Let us first pick $(u_{0},\varphi)\in D$. Thanks to
\eqref{cocycle}, we have
\begin{align*}
&  \quad\Xi(t,\omega,(u_{0},\varphi))\\[0.4ex]
&  =(v(t;0,\omega,(u_{0}-\phi z_{\ast}(\omega),(\mathcal{J}_{\omega,0}%
{\varphi}))),v_{t}(\cdot;0,\omega,(u_{0}-\phi z_{\ast}(\omega),(\mathcal{J}%
_{\omega,0}{\varphi})))\\[0.4ex]
&  ~~+(\phi z_{\ast}(\theta_{t}\omega),\phi z_{\ast}(\theta_{t+\cdot}%
\omega))\\[0.4ex]
&  =\Psi(t,\omega,(u_{0}-\phi z_{\ast}(\omega),\varphi-\phi z_{\ast}%
(\theta_{\cdot}\omega)))+(\phi z_{\ast}(\theta_{t}\omega),\phi z_{\ast}%
(\theta_{t+\cdot}\omega)).
\end{align*}
For the sake of simplicity, denote by $y(t,\omega):=(v(t,\omega),\eta
^{t}(s,\omega))$ the solution to \eqref{eq1-5} with initial value $(v_{0}%
,\eta_{0})=(u_{0}-\phi z_{\ast}(\omega),\mathcal{J}_{\omega,0}{\varphi})$.
Now, for every $\omega\in\Omega$, we multiply the first equation of
\eqref{eq1-5} by $v(t)$ in $H$ and the second equation of \eqref{eq1-5} by
$\eta^{t}$ in $L_{\mu}^{2}(\mathbb{R}^{+};V)$, respectively. Then, by means of
the same estimates as in the proof of Theorem \ref{thm3-3} (cf.
\eqref{eq3-11}) and the Poincar\'{e} inequality, we obtain
\begin{equation}
\frac{d}{dt}\Vert y(t)\Vert_{\mathcal{H}}^{2}+\frac{m\lambda_{1}}{2}%
|v(t)|^{2}+\frac{m}{2}\Vert v(t)\Vert^{2}+2((\eta^{t}(s),(\eta^{t}%
(s))^{\prime}))_{\mu}+\frac{f_{0}}{2D}\Vert v(t)\Vert_{2p}^{2p}\leq\Theta
_{1}(\theta_{t}\omega), \label{eq4-1}%
\end{equation}
(see \eqref{theta1} for the expression of $\Theta_{1}(\cdot)$). With the help
of condition $(h_{2})$, we infer
\begin{equation}
2(((\eta^{t})^{\prime},\eta^{t}))_{\mu}=-\int_{0}^{\infty}\mu^{\prime
}(s)|\nabla\eta^{t}(s)|^{2}ds\geq\varpi\int_{0}^{\infty}\mu(s)|\nabla\eta
^{t}(s)|^{2}ds:=\varpi\Vert\eta^{t}\Vert_{\mu}^{2}. \label{eq44}%
\end{equation}
Recall that $0<\gamma<\min\{\frac{m\lambda_{1}}{2},\varpi\}$, which, together
with \eqref{eq4-1} and \eqref{eq44}, implies that
\begin{equation}
\frac{d}{dt}\Vert y(t)\Vert_{\mathcal{H}}^{2}+\gamma\Vert y(t)\Vert
_{\mathcal{H}}^{2}+\frac{m}{2}\Vert v(t)\Vert^{2}+\frac{f_{0}}{2D}\Vert
v(t)\Vert_{2p}^{2p}\leq\Theta_{1}(\theta_{t}\omega). \label{eq4-12}%
\end{equation}
Next, multiplying the above inequality by $e^{\gamma t}$ and integrating over
$(0,t)$, neglecting the last term on the left hand side of \eqref{eq4-12}, we
find%
\begin{equation}
\Vert y(t)\Vert_{\mathcal{H}}^{2}+\frac{m}{2}\int_{0}^{t}e^{-\gamma(t-s)}\Vert
v(s)\Vert^{2}ds\leq\Vert y_{0}\Vert_{\mathcal{H}}^{2}e^{-\gamma t}+\int%
_{0}^{t}e^{-\gamma(t-s)}\Theta_{1}(\theta_{s}\omega)ds. \label{eq4-2}%
\end{equation}
Then,
\begin{equation}%
\begin{split}
\frac{m}{2}\Vert v_{t}\Vert_{L_{V}^{2}}^{2}  &  =\frac{m}{2}\int_{-\infty}%
^{0}e^{-\gamma(t-s)}\Vert\varphi(s)-\phi z_{\ast}(\theta_{s}\omega)\Vert
^{2}ds+\frac{m}{2}\int_{0}^{t}e^{-\gamma(t-s)}\Vert v(s)\Vert^{2}ds\\[0.8ex]
&  \leq\frac{m}{2}e^{-\gamma t}\left(  \Vert\varphi-\phi z_{\ast}%
(\theta_{\cdot}\omega)\Vert_{L_{V}^{2}}^{2}\right)  +\Vert y_{0}%
\Vert_{\mathcal{H}}^{2}e^{-\gamma t}+\int_{0}^{t}e^{-\gamma(t-s)}\Theta
_{1}(\theta_{s}\omega)ds.
\end{split}
\label{eq4-3}%
\end{equation}
On account of Corollary \ref{cor3-2}, we infer
\begin{equation}
\Vert y_{0}\Vert_{\mathcal{H}}^{2}=|v_{0}|^{2}+\Vert\mathcal{J}_{\omega
,0}{\varphi}\Vert_{L_{\mu}^{2}(\mathbb{R}^{+};V)}^{2}\leq|u_{0}-\phi z_{\ast
}(\omega)|^{2}+2K_{\mu}\left(  \Vert\varphi-\phi z_{\ast}(\theta_{\cdot}%
\omega)\Vert_{L_{V}^{2}}^{2}\right)  . \label{eq4-4}%
\end{equation}
Hence, collecting \eqref{eq4-2}-\eqref{eq4-4}, we arrive at
\begin{equation}%
\begin{split}
&  ~\quad\Vert\Psi(t,\omega,(u_{0}-\phi z_{\ast}(\omega),\varphi_{v}%
))\Vert_{X}^{2}=|v(t)|^{2}+\Vert v_{t}\Vert_{L_{V}^{2}}^{2}\leq\Vert
y(t)\Vert_{\mathcal{H}}^{2}+\Vert v_{t}\Vert_{L_{V}^{2}}^{2}\\[0.8ex]
&  \leq K_{1}e^{-\gamma t}\left(  |u_{0}-\phi z_{\ast}(\omega)|^{2}%
+\Vert\varphi-\phi z_{\ast}(\theta_{\cdot}\omega)\Vert_{L_{V}^{2}}^{2}\right)
+K_{2}\int_{0}^{t}e^{-\gamma(t-s)}\Theta_{1}(\theta_{s}\omega)ds,
\end{split}
\label{eq48}%
\end{equation}
where $K_{1},K_{2}>0$ are constants which neither depend on $\omega$ nor on
the initial functions. Now, replacing $\omega$ by $\theta_{-t}\omega$ in
\eqref{eq48}, we obtain
\begin{equation}%
\begin{split}
&  ~\quad\Vert\Psi(t,\theta_{-t}\omega,(u_{0}-\phi z_{\ast}(\theta_{-t}%
\omega),\varphi-\phi z_{\ast}(\theta_{-t+\cdot}\omega)))\Vert_{X}^{2}\\[0.4ex]
&  \leq K_{1}e^{-\gamma t}\left(  |u_{0}-\phi z_{\ast}(\theta_{-t}\omega
)|^{2}+\Vert\varphi-\phi z_{\ast}(\theta_{-t+\cdot}\omega)\Vert_{L_{V}^{2}%
}^{2}\right)  +K_{2}\int_{0}^{t}e^{-\gamma(t-s)}\Theta_{1}(\theta_{-t+s}%
\omega)ds\\[0.4ex]
&  \leq K_{1}e^{-\gamma t}\left(  |u_{0}-\phi z_{\ast}(\theta_{-t}\omega
)|^{2}+\Vert\varphi-\phi z_{\ast}(\theta_{-t+\cdot}\omega)\Vert_{L_{V}^{2}%
}^{2}\right)  +K_{2}\int_{-\infty}^{0}e^{\gamma s}{\Theta_{1}(\theta_{s}%
\omega)}ds.
\end{split}
\label{eq49}%
\end{equation}
Therefore, for any $(u_{0},\varphi)\in D(\theta_{-t}\omega)$, we have
\begin{equation}%
\begin{split}
&  \quad~\Vert\Xi(t,\theta_{-t}\omega,(u_{0},\varphi))\Vert_{X}^{2}\\[0.8ex]
&  =\Vert\Psi(t,\theta_{-t}\omega,(u_{0}-\phi z_{\ast}(\theta_{-t}%
\omega),\varphi-\phi z_{\ast}(\theta_{-t+\cdot}\omega)))+(\phi z_{\ast}%
(\omega),\phi z_{\ast}(\theta_{\cdot}\omega))\Vert_{X}^{2}\\[0.4ex]
&  \leq2K_{1}e^{-\gamma t}\left(  |u_{0}-\phi z_{\ast}(\theta_{-t}\omega
)|^{2}+\Vert\varphi-\phi z_{\ast}(\theta_{-t+\cdot}\omega)\Vert_{L_{V}^{2}%
}^{2}\right) \\
&  \quad+2K_{2}\int_{-\infty}^{0}e^{\gamma s}\Theta_{1}(\omega)ds+2|\phi
|^{2}|z_{\ast}(\omega)|^{2}+2\Vert\phi z_{\ast}(\theta_{\cdot}\omega
)\Vert_{L_{V}^{2}}^{2}\\[0.4ex]
&  \leq2K_{1}e^{-\gamma t}\Vert\widetilde{D}(\theta_{-t}\omega)\Vert
^{2}+2K_{2}\int_{-\infty}^{0}e^{\gamma s}{\Theta_{1}(\theta_{s}\omega
)ds}+2|\phi|^{2}|z_{\ast}(\omega)|^{2}+2\Vert\phi z_{\ast}(\theta_{\cdot
}\omega)\Vert_{L_{V}^{2}}^{2}.
\end{split}
\end{equation}
Consequently, let
\[
\rho^{2}(\omega)=1+2K_{2}\int_{-\infty}^{0}e^{\gamma s}{\Theta_{1}(\theta
_{s}\omega)}ds+2|\phi|^{2}|z_{\ast}(\omega)|^{2}+2\Vert\phi z_{\ast}%
(\theta_{\cdot}\omega)\Vert_{L_{V}^{2}}^{2}.
\]
Taking into account the temperedness of $\widetilde{D}$ and the property of
Ornstein-Uhlenbeck process, it is straightforward to check that the set,
\[
B_{\rho}(\omega)=\{(u_{0},\varphi)\in X:\Vert(u_{0},\varphi)\Vert_{X}\leq
\rho(\omega)\},
\]
is tempered, i.e., $B_{\rho}(\omega)$ belongs to $\mathcal{D}$ and is pullback
absorbing for the universe $\mathcal{D}$. The proof of this lemma is finished.
$\Box$

\begin{remark}
\label{absorbingv} It follows from the proof of Lemma \ref{absorbing} that,
under the assumptions of Theorem \ref{thm3-3}, there exists $\tilde{B}%
\in\mathcal{D}$ which is $\mathcal{D}$-pullback absorbing for the RDS $\Psi$.
In other words, for any given $\omega\in\Omega$ and $D\in\mathcal{D}$, there
exists $\tilde{t}_{0}:=\tilde{t}_{0}(\omega,D)\geq0$, such that
\[
\Psi(t,\theta_{-t}\omega,D(\theta_{-t}\omega))\subset\tilde{B}(\omega
),\ \quad\text{for all}\ t\geq\tilde{t}_{0}(\omega,D).
\]

\end{remark}

Next, we will prove the asymptotic compactness of the cocycle $\Xi$. Namely,
we will show that for any $\omega\in\Omega$, $D\in\mathcal{D}$ and for any
sequence $t_{n}\to+\infty$, $(u_{0}^{n},\varphi^{n})\in D(\theta_{-t_{n}%
}\omega)$, the sequence $\{\Xi(t_{n},\theta_{-t_{n}}\omega, (u_{0}^{n}%
,\varphi^{n}))\}$ possesses a convergent subsequence in $X$. To this end, let
us first prove an auxiliary result.

\begin{lemma}
\label{lem4-5} Assume the hypotheses in Theorem \ref{thm3-3} hold. Let
$\{v_{0}^{n},\varphi_{v}^{n}\}$ be a sequence such that $(v_{0}^{n}%
,\varphi_{v}^{n})\rightarrow(v_{0},\varphi_{v})$ weakly in $X$ as
$n\rightarrow\infty$. Then, for every $\omega\in\Omega$, $\Psi(t,\omega
,(v_{0}^{n},\varphi_{v}^{n}))=(v^{n}(t),v_{t}^{n}(\cdot))$ fulfills:
\begin{equation}
v^{n}\rightarrow v\quad\mbox{in}\quad C([r,T];H)~~\mbox{for all}~~0<r<T;
\label{eq4-6}%
\end{equation}%
\begin{equation}
v^{n}\rightarrow v\quad\mbox{weakly in}\quad L^{2}%
(0,T;V)~~\mbox{for all}~~T>0; \label{eq4-7}%
\end{equation}%
\begin{equation}
v^{n}\rightarrow v\quad\mbox{in}\quad L^{2}(0,T;H)~~\mbox{for all}~~T>0;
\label{eq4-8}%
\end{equation}%
\begin{equation}
\limsup_{n\rightarrow\infty}\Vert v_{t}^{n}-v_{t}\Vert_{L_{V}^{2}}^{2}\leq
Ke^{-\gamma t}\limsup_{n\rightarrow\infty}\left(  |v_{0}^{n}-v_{0}|^{2}%
+\Vert\varphi_{v}^{n}-\varphi_{v}\Vert_{L_{V}^{2}}^{2}\right)
~\mbox{for all}~~t\geq0, \label{eq4-9}%
\end{equation}
where $K=(1+\frac{2K_{\mu}}{m}+\frac{1}{m})$ and $(v(t),v_{t}($%
\textperiodcentered$))=\Psi(t,\omega,(v_{0},\varphi_{v}))$. Moreover, if
$(v_{0}^{n},\varphi_{v}^{n})\rightarrow(v_{0},\varphi_{v})$ strongly in $X$ as
$n\rightarrow\infty$, then
\begin{equation}
v^{n}\rightarrow v\quad\mbox{in}\quad L^{2}(0,T;V)~\mbox{for all}~~T>0;
\label{eq4-10}%
\end{equation}%
\begin{equation}
v_{t}^{n}\rightarrow v_{t}\quad\mbox{in}\quad L_{V}^{2}~~\mbox{for all}~~t\geq
0. \label{eq4-11}%
\end{equation}

\end{lemma}

\textbf{Proof.} Let $T>0$ be arbitrary. Integrating in \eqref{eq4-12}, we
deduce that $v^{n}$ is bounded in $L^{\infty}(0,T;H)$, $L^{2}(0,T;V)$ and
$L^{2p}(0,T;L^{2p}(\mathcal{O}))$, $\eta_{n}$ is bounded in $L^{\infty
}(0,T;L_{\mu}^{2}(\mathbb{R}^{+};V))$. Hence, passing to a subsequence, for
every $\omega\in\Omega$, we have
\[
\left\{
\begin{array}
[c]{rcl}%
\begin{aligned} &v^{n} \rightarrow v ~~\mbox{weak-$*$ in}~~ L^{\infty}(0,T;H); \\ &v^n \rightarrow v ~~ \mbox{weakly in}~~ L^{2}(0,T;V);\\ &v^n\rightarrow v~~ \mbox{weakly in}~~ L^{2p}(0,T;L^{2p}(\mathcal{O}));\\ & \eta_n\rightarrow \eta~~\mbox{weak-star in} ~~L^\infty(0,T;L_\mu^2(\mathbb{R}^+;V)), \end{aligned} &
&
\end{array}
\right.
\]
thus \eqref{eq4-7} holds. By the same arguments in the proof of Theorem
\ref{thm3-3}, we deduce
\begin{equation}
\label{eq41811}\left\{
\begin{array}
[c]{rcl}%
\begin{aligned} &\frac{dv^n}{dt}\rightarrow \frac{dv}{dt}~\mbox{weakly in} ~~L^2(0,T;V^*)+L^q(0,T;L^q(\mathcal{O}));\\ &f(v^n+\phi z_*(\theta_\cdot\omega))\rightarrow \chi~~\mbox{weakly in} ~~L^q(0,T;L^q(\mathcal{O})). \end{aligned} &
&
\end{array}
\right.
\end{equation}
In view of \eqref{eq4-7} and the above results, making use of the Compactness
Theorem \cite{R1}, we infer that \eqref{eq4-8} is true. Thus, $v^{n}%
(t,x,\omega)\rightarrow v(t,x,\omega)$, $f(v^{n}(t,x,\omega)+\phi z_{\ast
}(\theta_{t}\omega))\rightarrow f(v(t,x,\omega)+\phi z_{\ast}(\theta_{t}%
\omega))$ for a.a. $(t,x)\in(0,T)\times\mathcal{O}$. Also it follows from
\cite[Lemma 1.3]{L1} that $\chi=f(v+\phi z_{\ast}(\theta_{\cdot}\omega))$.

By proceeding as in the proof of Theorem \ref{thm3-3}, we obtain that
$y=(v,\eta)$ is a solution to problem \eqref{eq1-5} with initial value
$y(0)=(v_{0},\eta_{0}):=(v_{0},\mathcal{J}_{\omega,0}{\varphi})$. Thanks to
the uniqueness of solution, a standard argument implies that the above
convergences are true for the whole sequence. Furthermore, we know that
$v^{n}\in C([0,T];H)$ and $v\in C([0,T];H)$ for each $\omega\in\Omega$.

Since $\{(v^{n})^{\prime}\}$ is bounded in $L^{q}(0,T;V^{\ast}+L^{q}%
(\mathcal{O}))$, we have that $\{v^{n}\}$ is equicontinuous in $V^{\ast}%
+L^{q}(\mathcal{O})$ on $[0,T]$. Indeed,%
\begin{equation}
\Vert v^{n}(s_{2})-v^{n}(s_{1})\Vert_{V^{\ast}+L^{q}}\leq\int_{s_{1}}^{s_{2}%
}\left\Vert (v^{n})^{\prime}\right\Vert _{V^{\ast}+L^{q}}ds\leq\left\vert
s_{2}-s_{1}\right\vert ^{\frac{1}{2p}}\left\Vert (v^{n})^{\prime}\right\Vert
_{L^{q}(0,T;V^{\ast}+L^{q}(\mathcal{O}))}. \label{eq2}%
\end{equation}
In addition, as $\{v^{n}\}$ is bounded in $C([0,T];H)$ and the embedding
$H\subset V^{\ast}$ is compact, by the Arzel\`{a}-Ascoli theorem, we obtain
(relabeled the same) that
\begin{equation}
v^{n}\rightarrow v~~\mbox{strongly in}~~C([0,T];V^{\ast}+L^{q}(\mathcal{O})).
\label{eq411}%
\end{equation}
Now, consider a sequence $\{s_{n}\}\in\lbrack0,T]$ which converges to
$s_{\ast}\in(0,T]$. Since $\{v^{n}\}$ is bounded in $C([0,T];H)$, there exist
a subsequence of $\{v^{n}(s_{n})\}$ and $\hat{u}\in H$ such that%
\[
v^{n}(s_{n})\rightharpoonup\hat{u}~~\mbox{weakly in}~~H.
\]
Using (\ref{eq411}) we deduce that $\hat{u}=u(s_{\ast})$ and that the whole
sequence converges. Therefore,
\begin{equation}
|v(s_{\ast})|\leq\liminf_{n\rightarrow\infty}|v^{n}(s_{n})|. \label{eq711}%
\end{equation}
We will prove that $v^{n}(s_{n})\rightarrow v(s_{\ast})$ strongly in $H$,
which implies (\ref{eq4-6}). By Corollary \ref{CorExistSol}, $v^{n}$ and $v$
are weak solutions to problem (\ref{eq1-4}), so multiplying the equation by
$v^{n}$, we obtain that
\[%
\begin{split}
&  \quad~~\frac{1}{2}\frac{d}{dt}|v^{n}(t)|^{2}+m\Vert v^{n}\Vert^{2}%
+(f(v^{n}+\phi z_{\ast}(\theta_{t}\omega)),v^{n})\\[0.8ex]
&  \leq M|z_{\ast}(\theta_{t}\omega)|\Vert\phi\Vert\Vert v^{n}\Vert+\left(
\int_{-\infty}^{t}k(t-s)\Delta v^{n}(s)ds,v^{n}(t)\right)  +(h,v^{n}%
(t))\\[0.8ex]
&  ~~+|z_{\ast}(\theta_{t}\omega)||\phi||v^{n}|+(z_{k}^{\phi}(\theta_{t}%
\omega),v^{n}).
\end{split}
\]
By using similar arguments as in Theorem \ref{thm3-3} and the Young
inequality, we deduce
\begin{equation}%
\begin{split}
&  ~~\quad\frac{d}{dt}|v^{n}|^{2}+m\Vert v^{n}\Vert^{2}+\frac{f_{0}}{2D}\Vert
v^{n}\Vert_{2p}^{2p}\leq2\alpha|\mathcal{O}|+C_{1}(\theta_{t}\omega
)(1+\Vert\phi\Vert_{2p}^{2p})\\[0.8ex]
&  ~~+\frac{4M^{2}}{m}|z_{\ast}(\theta_{t}\omega)|^{2}\Vert\phi\Vert^{2}%
+\frac{4}{m\lambda_{1}}|h|^{2}+\frac{4}{m\lambda_{1}}|z_{\ast}(\theta
_{t}\omega)|^{2}|\phi|^{2}\\[0.8ex]
&  ~~+\frac{C_{2}(\theta_{t}\omega)}{m}\left\Vert \phi\right\Vert ^{2}%
+2\int_{-\infty}^{t}k(t-s)\Vert v^{n}(s)\Vert ds\Vert v^{n}(t)\Vert.
\end{split}
\label{eq27}%
\end{equation}
Now we estimate the last term in the above inequality. Notice that
\[
\int_{-\infty}^{t}k(t-s)\Vert v^{n}(s)\Vert ds=\int_{-\infty}^{0}k(t-s)\Vert
v^{n}(s)\Vert ds+\int_{0}^{t}k(t-s)\Vert v^{n}(s)\Vert ds:=I_{1}+I_{2}.
\]
For $I_{1}$, by the fact that $\varphi_{v}\in L_{V}^{2}$, $\gamma<\min
\{\frac{m\lambda_{1}}{2},\varpi\}$ and $k(t)\leq M_{1},\ \forall t\in
\lbrack0,\infty)$ (see Remark \ref{rem2-1}), we find
\[%
\begin{split}
I_{1}  &  =\int_{-\infty}^{0}k(t-s)e^{-\frac{\gamma s}{2}}e^{\frac{\gamma
s}{2}}\Vert\varphi_{v}(s)\Vert ds\leq\left(  \int_{-\infty}^{0}k^{2}%
(t-s)e^{-\gamma s}ds\right)  ^{\frac{1}{2}}\\[0.8ex]
&  ~~\times\left(  \int_{-\infty}^{0}e^{\gamma s}\Vert\varphi_{v}\Vert
^{2}ds\right)  ^{\frac{1}{2}}\leq\Vert\varphi_{v}\Vert_{L_{V}^{2}}\left(
\int_{t}^{\infty}k(s)e^{-\gamma(t-s)}ds\right)  ^{\frac{1}{2}}M_{1}^{\frac
{1}{2}}\\[0.8ex]
&  \leq\frac{\Vert\varphi_{v}\Vert_{L_{V}^{2}}M_{1}^{\frac{1}{2}}}%
{\varpi^{\frac{1}{2}}}\left(  \int_{t}^{\infty}\mu(s)e^{-\gamma(t-s)}%
ds\right)  ^{\frac{1}{2}}\\[0.8ex]
&  \leq\frac{|\varphi_{v}\Vert_{L_{V}^{2}}M_{1}^{\frac{1}{2}}}{\varpi
^{\frac{1}{2}}}\left(  \int_{t}^{\infty}\mu(t)e^{-\varpi(s-t)}e^{-\gamma
(t-s)}ds\right)  ^{\frac{1}{2}}\\[0.8ex]
&  \leq\frac{M_{1}^{\frac{1}{2}}\mu^{\frac{1}{2}}(t)\Vert\varphi_{v}%
\Vert_{L_{V}^{2}}}{\varpi^{\frac{1}{2}}(\varpi-\gamma)^{\frac{1}{2}}}.
\end{split}
\]
For $I_{2}$, by means of the property $k(t)\leq M_{1}$ and the boundedness of
$v^{n}$ in $L^{2}(0,T;V)$, there exists a constant $M^{\prime\prime}$ such
that
\[
I_{2}\leq M_{1}\int_{0}^{t}\Vert v^{n}(s)\Vert ds\leq M_{1}M^{\prime\prime
}\sqrt{t}.
\]
Therefore, it follows from the above inequalities and \eqref{eq27} that for
every $t<T$,%
\begin{equation}%
\begin{split}
\frac{d}{dt}|v^{n}(t)|^{2}  &  +m\Vert v^{n}(t)\Vert^{2}+\frac{f_{0}}{2D}\Vert
v^{n}(t)\Vert_{2p}^{2p}\leq2\alpha|\mathcal{O}|+C_{1}(\theta_{t}%
\omega)(1+\Vert\phi\Vert_{2p}^{2p})\\[0.8ex]
&  ~~+\frac{4M^{2}}{m}|z_{\ast}(\theta_{t}\omega)|^{2}\Vert\phi\Vert^{2}%
+\frac{4}{m\lambda_{1}}|h|^{2}+\frac{4}{m\lambda_{1}}|z_{\ast}(\theta
_{t}\omega)|^{2}|\phi|^{2}\\[0.8ex]
&  ~~+\frac{C_{2}(\theta_{t}\omega)}{m}\left\Vert \phi\right\Vert
^{2}+2\left(  \frac{M_{1}^{\frac{1}{2}}\mu^{\frac{1}{2}}(t)\Vert\varphi
_{v}\Vert_{L_{V}^{2}}}{\varpi^{\frac{1}{2}}(\varpi-\gamma)^{\frac{1}{2}}%
}+M_{1}M^{\prime\prime}\sqrt{T}\right)  \Vert v^{n}(t)\Vert.
\end{split}
\label{eq271}%
\end{equation}
We will estimate the last term of the above inequality. By the Young
inequality, we have%
\begin{equation}
2\frac{M_{1}^{\frac{1}{2}}\mu^{\frac{1}{2}}(t)\Vert\varphi_{v}\Vert_{L_{V}%
^{2}}}{\varpi^{\frac{1}{2}}(\varpi-\gamma)^{\frac{1}{2}}}\Vert v(t)\Vert
\leq\frac{4M_{1}\mu(t)\Vert\varphi_{v}\Vert_{L_{V}^{2}}^{2}}{\varpi
m(\varpi-\gamma)}+\frac{m}{4}\Vert v(t)\Vert^{2}, \label{eq272}%
\end{equation}
and%
\begin{equation}
\quad2M_{1}M^{\prime\prime}\sqrt{T}\Vert v(t)\Vert\leq\frac{4(M_{1}%
)^{2}(M^{\prime\prime})^{2}T}{m}+\frac{m}{4}\Vert v(t)\Vert^{2}. \label{eq273}%
\end{equation}
Collecting \eqref{eq271}-\eqref{eq273}, we obtain%
\begin{equation}%
\begin{split}
\frac{d}{dt}|v^{n}(t)|^{2}  &  +\frac{m}{2}\Vert v^{n}(t)\Vert^{2}+\frac
{f_{0}}{2D}\Vert v^{n}(t)\Vert_{2p}^{2p}\leq2\alpha|\mathcal{O}|+C_{1}%
(\theta_{t}\omega)(1+\Vert\phi\Vert_{2p}^{2p})\\[0.8ex]
&  +\frac{4M^{2}}{m}|z_{\ast}(\theta_{t}\omega)|^{2}\Vert\phi\Vert^{2}%
+\frac{4}{m\lambda_{1}}|h|^{2}+\frac{4}{m\lambda_{1}}|z_{\ast}(\theta
_{t}\omega)|^{2}|\phi|^{2}\\[0.8ex]
&  +\frac{4M_{1}\mu(t)\Vert\varphi_{v}\Vert_{L_{V}^{2}}^{2}}{\varpi
m(\varpi-\gamma)}+\frac{4(M_{1})^{2}(M^{\prime\prime})^{2}T}{m}+\frac
{C_{2}(\theta_{t}\omega)}{m}\left\Vert \phi\right\Vert ^{2},
\end{split}
\label{eq274}%
\end{equation}
and the same is true for the function $v$. Hence, we define the functions%
\[%
\begin{split}
J_{n}(t)  &  =|v^{n}(t)|^{2}-2\alpha|\mathcal{O}|t-\frac{4(M_{1}%
)^{2}(M^{\prime\prime})^{2}T}{m}t-\int_{0}^{t}C_{1}(\theta_{r}\omega
)(1+\Vert\phi\Vert_{2p}^{2p})dr\\[0.8ex]
&  ~~-\int_{0}^{t}\left(  \frac{4M^{2}}{m}|z_{\ast}(\theta_{r}\omega
)|^{2}\Vert\phi\Vert^{2}+\frac{4}{m\lambda_{1}}|z_{\ast}(\theta_{r}%
\omega)|^{2}|\phi|^{2}\right)  dr\\[0.8ex]
&  ~~-\frac{4\left\vert h\right\vert ^{2}}{m\lambda_{1}}t-\frac{4M_{1}%
\Vert\varphi_{v}\Vert_{L_{V}^{2}}^{2}}{\varpi m(\varpi-\gamma)}\int_{0}^{t}%
\mu(r)dr-\int_{0}^{t}\frac{C_{2}(\theta_{r}\omega)}{m}\left\Vert
\phi\right\Vert ^{2}dr,
\end{split}
\]%
\[%
\begin{split}
J(t)  &  =|v(t)|^{2}-2\alpha|\mathcal{O}|t-\frac{4(M_{1})^{2}(M^{\prime\prime
})^{2}T}{m}t-\int_{0}^{t}C_{1}(\theta_{r}\omega)(1+\Vert\phi\Vert_{2p}%
^{2p})dr\\[0.8ex]
&  ~~-\int_{0}^{t}\left(  \frac{4M^{2}}{m}|z_{\ast}(\theta_{r}\omega
)|^{2}\Vert\phi\Vert^{2}+\frac{4}{m\lambda_{1}}|z_{\ast}(\theta_{r}%
\omega)|^{2}|\phi|^{2}\right)  dr\\[0.8ex]
&  ~~-\frac{4\left\vert h\right\vert ^{2}}{m\lambda_{1}}t-\frac{4M_{1}%
\Vert\varphi_{v}\Vert_{L_{V}^{2}}^{2}}{\varpi m(\varpi-\gamma)}\int_{0}^{t}%
\mu(r)dr-\int_{0}^{t}\frac{C_{2}(\theta_{r}\omega)}{m}\left\Vert
\phi\right\Vert ^{2}dr.
\end{split}
\]
From the regularity of $v$ and all $v^{n}$, together with (\ref{eq274}), it
holds that these functions $J$ and $J_{n}$ are continuous and non-increasing
on $[0,T]$, and
\[
J_{n}(s)\rightarrow J(s)~~\mbox{a.e.}~~s\in\lbrack0,T]~\mbox{as}~n\rightarrow
\infty.
\]
Hence, there exists a sequence $\{\tilde{t}_{k}\}\in(0,s_{\ast})$ such that
$\tilde{t}_{k}\rightarrow s_{\ast}$ when $k\rightarrow\infty$, and
\[
\lim_{n\rightarrow\infty}J_{n}(\tilde{t}_{k})=J(\tilde{t}_{k}),\qquad\forall
k\geq1.
\]
Fix an arbitrary value $\epsilon>0$. From the continuity of $J$ on $[0,T]$,
there exists $k(\epsilon)\geq1$ such that
\[
|J(\tilde{t}_{k})-J(s_{\ast})|\leq\frac{\epsilon}{2},~~\forall k\geq
k(\epsilon).
\]
Now consider $n(\epsilon)\geq1$ such that
\[
t_{n}\geq\tilde{t}_{k(\epsilon)}~~\mbox{and}~~|J_{n}(\tilde{t}_{k(\epsilon
)})-J(\tilde{t}_{k(\epsilon)})|\leq\frac{\epsilon}{2},\quad\forall n\geq
n(\epsilon).
\]
Then, since all $J_{n}$ are non-increasing, we deduce that
\[%
\begin{split}
J_{n}(t_{n})-J(s_{\ast})  &  \leq J_{n}(\tilde{t}_{k(\epsilon)})-J(s_{\ast
})\leq|J_{n}(\tilde{t}_{k(\epsilon)})-J(s_{\ast})|\\
&  \leq|J_{n}(\tilde{t}_{k(\epsilon)})-J(\tilde{t}_{k(\epsilon)}%
)|+|J(\tilde{t}_{k(\epsilon)})-J(s_{\ast})|\leq\epsilon,\qquad\forall n\geq
n(\epsilon).
\end{split}
\]
As $\epsilon>0$ is arbitrary, we obtain
\[
\limsup_{n\rightarrow\infty}J_{n}(t_{n})\leq J(s_{\ast}).
\]
Thus,
\begin{equation}
\limsup_{n\rightarrow\infty}\ |u^{n}(t_{n})|\leq|u(s_{\ast})|. \label{eq911}%
\end{equation}
Therefore, \eqref{eq711} and \eqref{eq911} imply that $v^{n}(s_{n})\rightarrow
v(s_{\ast})$ strongly in $H$, and (\ref{eq4-6}) holds true.

Define the functions $\bar{y}^{n}=y^{n}-y$ and $\bar{\eta}_{n}^{t}=\eta
_{n}^{t}-\eta^{t}$ with $\bar{y}_{0}^{n}=\bar{y}_{0}^{n}-y_{0}$, where
$y_{0}=(v_{0},\eta_{0})$. Similar to the uniqueness part in the proof of
Theorem \ref{thm3-3}, for every $\omega\in\Omega$, we have
\[%
\begin{split}
\frac{d}{dt}\Vert\bar{y}^{n}\Vert_{\mathcal{H}}^{2}  &  +2(((\bar{\eta}%
_{n}^{t})^{\prime},\eta_{n}^{t}))_{\mu}\leq-2\int_{\mathcal{O}}(f(v^{n}+\phi
z_{\ast}(\theta_{t}\omega))-f(v+\phi z_{\ast}(\theta_{t}\omega)))(v^{n}%
-v)dx\\[0.8ex]
&  -2\int_{\mathcal{O}}\left(  a(l(v^{n})+l(\phi)z_{\ast}(\theta_{t}%
\omega))\nabla v^{n}-a(l(v)+l(\phi)z_{\ast}(\theta_{t}\omega))\nabla v\right)
\cdot\nabla(v^{n}-v)dx\\[0.8ex]
&  +2\int_{\mathcal{O}}(a(l(v^{n})+l(\phi)z_{\ast}(\theta_{t}\omega
))-a(l(v)+l(\phi)z_{\ast}(\theta_{t}\omega)))z_{\ast}(\theta_{t}\omega
)\Delta\phi(v^{n}-v)dx.
\end{split}
\]
Since $a$ is a locally Lipschitz function, by \eqref{eq1-2} and the Young
inequality, we find
\[%
\begin{split}
&  \quad-2\int_{\mathcal{O}}\left(  a(l(v^{n})+l(\phi)z_{\ast}(\theta
_{t}\omega))\nabla v^{n}-a(l(v)+l(\phi)z_{\ast}(\theta_{t}\omega))\nabla
v\right)  \cdot\nabla(v^{n}-v)dx\\[0.8ex]
&  \leq-2m\Vert v^{n}-v\Vert^{2}+2L_{a}(R)|l||v^{n}-v|\Vert v\Vert\Vert
v^{n}-v\Vert\\[0.8ex]
&  \leq(\alpha-2m)\Vert v^{n}-v\Vert^{2}+\frac{L_{a}^{2}(R)|l|^{2}}{\alpha
}|v^{n}-v|^{2}\Vert v\Vert^{2},
\end{split}
\]
where $\alpha\leq(m\lambda_{1}-\gamma)/\lambda_{1}$ and for all $n\geq1$,
$t\geq0$, and we have chosen $R>0$ such that $\{l(v^{n}(t)+l(\phi)z_{\ast
}(\theta_{t}\omega))\}_{t\in\lbrack\tau,T]}\subset\lbrack
-R,R],\ \{l(v(t)+l(\phi)z_{\ast}(\theta_{t}\omega))\}_{t\in\lbrack\tau
,T]}\subset\lbrack-R,R]$, which can be done because $\left\vert v^{n}\left(
t\right)  \right\vert $ are uniformly bounded in $[\tau,T]$. Then, by the
above estimates, we deduce
\[%
\begin{split}
&  ~\quad\frac{d}{dt}\Vert\bar{y}^{n}\Vert_{\mathcal{H}}^{2}+\gamma\Vert
\bar{y}^{n}\Vert_{\mathcal{H}}^{2}+m\Vert v^{n}-v\Vert^{2}\\[0.8ex]
&  \leq\frac{d}{dt}\Vert\bar{y}^{n}\Vert_{\mathcal{H}}^{2}+(2m-\alpha)\Vert
v^{n}-v\Vert^{2}+\varpi\int_{0}^{\infty}\mu(s)|\nabla\bar{\eta}_{n}%
^{t}(s)|^{2}ds\\[0.8ex]
&  \leq\frac{L_{a}^{2}(R)|l|^{2}}{\alpha}|v^{n}-v|^{2}\Vert v\Vert^{2}%
-2\int_{\mathcal{O}}(f(v^{n}+\phi z_{\ast}(\theta_{t}\omega))-f(v+\phi
z_{\ast}(\theta_{t}\omega)))(v^{n}-v)dx\\[0.8ex]
&  ~~+2\int_{\mathcal{O}}(a(l(v^{n})+l(\phi)z_{\ast}(\theta_{t}\omega
))-a(l(v)+l(\phi)z_{\ast}(\theta_{t}\omega)))z_{\ast}(\theta_{t}\omega
)\Delta\phi(v^{n}-v)dx,
\end{split}
\]
where we have used that $0<\gamma\leq\min\{(m-\alpha)\lambda_{1},\varpi\}$ by
the choice of $\alpha$. Multiplying by $e^{\gamma t}$ on both sides of the
above inequality and integrating over $(0,t)$, we obtain
\[%
\begin{split}
&  ~\quad\Vert\bar{y}^{n}(t)\Vert_{\mathcal{H}}^{2}+m\int_{0}^{t}%
e^{-\gamma(t-s)}\Vert v^{n}(s)-v(s)\Vert^{2}ds\\
&  \leq e^{-\gamma t}\Vert\bar{y}_{0}^{n}\Vert_{\mathcal{H}}^{2}+\frac
{L_{a}^{2}(R)|l|^{2}}{\alpha}\int_{0}^{t}e^{-\gamma(t-s)}|v^{n}-v|^{2}\Vert
v\Vert^{2}ds\\
&  -2\int_{0}^{t}e^{-\gamma(t-s)}\int_{\mathcal{O}}(f(v^{n}+\phi z_{\ast
}(\theta_{s}\omega))-f(v+\phi z_{\ast}(\theta_{s}\omega)))(v^{n}-v)dxds\\
&  +2\int_{0}^{t}e^{-\gamma(t-s)}\int_{\mathcal{O}}(a(l(v^{n})+l(\phi)z_{\ast
}(\theta_{s}\omega))-a(l(v)+l(\phi)z_{\ast}(\theta_{s}\omega)))|z_{\ast
}(\theta_{s}\omega)||\Delta\phi||v^{n}-v|dxds.
\end{split}
\]
On the one hand, by \eqref{eq4-6}, we know that $|v^{n}(s)-v(s)|^{2}\Vert
v(s)\Vert^{2}\rightarrow0$ and $|v^{n}(s)-v(s)|\rightarrow0$ for a.e.
$s\in(0,t)$. On the other hand, $e^{-\gamma(t-s)}|v^{n}(s)-v(s)|^{2}\Vert
v(s)\Vert^{2}$ and $e^{-\gamma(t-s)}(a(l(v^{n})+l(\phi)z_{\ast}(\theta
_{s}\omega))-a(l(v)+l(\phi)z_{\ast}(\theta_{s}\omega)))|z_{\ast}(\theta
_{s}\omega)||\Delta\phi||v^{n}-v|$ can be bounded by $4R^{2}e^{-\gamma(t-s)}$
$\Vert v(s)\Vert^{2}$ and $4MRe^{-\gamma(t-s)}\sup_{s\in\lbrack0,t]}|z_{\ast
}(\theta_{s}\omega)||\Delta\phi|$, respectively. Hence, the Lebesgue theorem
implies that
\[
\int_{0}^{t}e^{-\gamma(t-s)}|v^{n}(s)-v(s)|^{2}\Vert v(s)\Vert^{2}%
ds\rightarrow0\quad\mbox{as}\quad n\rightarrow0,
\]
and
\[%
\begin{split}
&  \int_{0}^{t}e^{-\gamma(t-s)}\int_{\mathcal{O}}(a(l(v^{n})+l(\phi)z_{\ast
}(\theta_{s}\omega))-a(l(v)+l(\phi)z_{\ast}(\theta_{s}\omega)))|z_{\ast
}(\theta_{s}\omega)||\Delta\phi||v^{n}-v|dxds\\[1ex]
&  \qquad\longrightarrow0\quad\mbox{as}\quad n\rightarrow0,
\end{split}
\]
respectively. Moreover, it follows from the argument after \eqref{eq41811}
that $f(v^{n}+\phi z_{\ast}(\theta_{s}\omega))\rightarrow f(v+\phi z_{\ast
}(\theta_{s}\omega))$ weakly in $L^{q}(0,T;L^{q}(\mathcal{O}))$ as
$n\rightarrow\infty$, therefore
\[
\int_{0}^{t}e^{-\gamma(t-s)}\int_{\mathcal{O}}\left(  f(v^{n}+\phi z_{\ast
}(\theta_{s}\omega))-f(v+\phi z_{\ast}(\theta_{s}\omega))\right)
vdxds\rightarrow0\quad\mbox{as}\quad n\rightarrow\infty.
\]
By \eqref{fcondition} and the Young inequality, we deduce that there are two
positive constants {$\kappa_{1}(\omega,\phi,T)>0$ and $\kappa_{2}>0$ such
that}
\[%
\begin{split}
&  ~\quad f(v^{n}(s,x,\omega)+\phi z_{\ast}(\theta_{s}\omega))v^{n}%
(s,x,\omega)\\
&  =f(v^{n}(s,x,\omega)+\phi z_{\ast}(\theta_{s}\omega))\left(  v^{n}%
(s,x,\omega\right)  +\phi z_{\ast}(\theta_{s}\omega))-f(v^{n}(s,x,\omega)+\phi
z_{\ast}(\theta_{s}\omega))\phi z_{\ast}(\theta_{s}\omega)\\[0.5ex]
&  \geq\frac{1}{2}f_{0}\left\vert v^{n}(s,x,\omega)+\phi z_{\ast}(\theta
_{s}\omega)\right\vert ^{2p}-\alpha-\beta\left(  1+\left\vert v^{n}%
(s,x,\omega)+\phi z_{\ast}(\theta_{s}\omega)\right\vert ^{2p-1}\right)
\left\vert \phi\right\vert |z_{\ast}(\theta_{s}\omega)|\\[0.5ex]
&  \geq-\kappa_{1}+\kappa_{2}\left\vert v^{n}(s,x,\omega)+\phi z_{\ast}%
(\theta_{s}\omega)\right\vert ^{2p}.
\end{split}
\]
%%
%%As $f(v^n(s,x,\omega)+\phi z_*(\theta_s\omega))v^n(s,x,\omega)=F(v^n,\omega)v^n=\sum_{k=1}^{2p}F_{2p-k}(v^n(s,x,\omega))^{k-1}v^n(s,x,\omega)\geq -\kappa_1+\kappa_2|v^n(t,x,\omega)|^{2p}$ (cf. \eqref{fcondition}) and $v^n(s,x,\omega)\rightarrow v(s,x,\omega)$, $f(v^n(s,x,\omega)+\phi z_*(\theta_s\omega))\rightarrow f(v(s,x,\omega)+\phi z_*(\theta_s\omega))$ for a.e. $(s,x,\omega)\in (0,T)\times\mathcal{O}\times\Omega$,
Thus, the Fatou-Lebesgue theorem implies that
\[%
\begin{split}
&  \qquad\limsup_{n\rightarrow\infty}\left(  -2\int_{0}^{t}e^{-\gamma
(t-s)}\int_{\mathcal{O}}f(v^{n}+\phi z_{\ast}(\theta_{s}\omega))v^{n}%
dxds\right) \\[0.8ex]
&  \leq-2\liminf_{n\rightarrow\infty}\int_{0}^{t}e^{-\gamma(t-s)}%
\int_{\mathcal{O}}f(v^{n}+\phi z_{\ast}(\theta_{s}\omega))v^{n}dxds\\[0.8ex]
&  \leq-2\int_{0}^{t}e^{-\gamma(t-s)}\int_{\mathcal{O}}\liminf_{n\rightarrow
\infty}f(v^{n}+\phi z_{\ast}(\theta_{s}\omega))v^{n}dxds\\[0.8ex]
&  =-2\int_{0}^{t}e^{-\gamma(t-s)}\int_{\mathcal{O}}f(v+\phi z_{\ast}%
(\theta_{s}\omega))vdxds.
\end{split}
\]
This inequality, together with
\begin{equation}
\int_{0}^{t}e^{-\gamma(t-s)}\int_{\mathcal{O}}f(v+\phi z_{\ast}(\theta
_{s}\omega))(v^{n}-v)dxds\rightarrow0\quad\mbox{as}\quad n\rightarrow\infty,
\label{eq4-17}%
\end{equation}
shows that
\[
\limsup_{n\rightarrow\infty}\left(  -2\int_{0}^{t}e^{-\gamma(t-s)}%
\int_{\mathcal{O}}(f(v^{n}+\phi z_{\ast}(\theta_{s}\omega))-f(v+\phi z_{\ast
}(\theta_{s}\omega)))(v^{n}-v)dxds\right)  {\leq0}.
\]
Notice that \eqref{eq4-17} follows from the facts that $f(v+\phi z_{\ast
}(\theta_{\cdot}\omega))\in L^{q}(0,T;L^{q}(\mathcal{O}))$ and $v^{n}%
\rightarrow v$ weakly in $L^{2p}(0,T;L^{2p}(\mathcal{O}))$ for every
$\omega\in\Omega$.

Collecting all inequalities derived above and using Corollary \ref{cor3-2}, we
find
\[%
\begin{split}
&  \qquad\limsup_{n\rightarrow\infty}\int_{0}^{t}e^{-\gamma(t-s)}\Vert
v^{n}(s)-v(s)\Vert^{2}ds\\
&  \leq\frac{1}{m}e^{-\gamma t}\limsup_{n\rightarrow\infty}\Vert\bar{y}%
_{0}^{n}\Vert_{\mathcal{H}}^{2}\\
&  \leq\frac{1}{m}e^{-\gamma t}\limsup_{n\rightarrow\infty}\left(
|v^{n}(0)-v(0)|^{2}+\int_{0}^{\infty}\mu(s)\Vert\bar{\eta}_{n}^{0}(s)\Vert
^{2}ds\right) \\
&  \leq\frac{1}{m}e^{-\gamma t}\limsup_{n\rightarrow\infty}\left(  |v_{0}%
^{n}-v_{0}|^{2}+2K_{\mu}\int_{-\infty}^{0}e^{\gamma s}\Vert\varphi_{v}%
^{n}(s)-\varphi_{v}(s)\Vert^{2}ds\right)  .
\end{split}
\]
Finally, \eqref{eq4-9} follows from
\[%
\begin{split}
\Vert v_{t}^{n}-v_{t}\Vert_{L_{V}^{2}}^{2}  &  =\int_{-t}^{0}e^{\gamma s}\Vert
v^{n}(t+s)-v(t+s)\Vert^{2}ds+\int_{-\infty}^{-t}e^{\gamma s}\Vert
v^{n}(t+s)-v(t+s)\Vert^{2}ds\\
&  =\int_{0}^{t}e^{-\gamma(t-s)}\Vert v^{n}(s)-v(s)\Vert^{2}ds+e^{-\gamma
t}\int_{-\infty}^{0}e^{\gamma s}\Vert\varphi_{v}^{n}(s)-\varphi_{v}%
(s)\Vert^{2}ds.
\end{split}
\]
If $(v_{0}^{n},\varphi_{v}^{n})\rightarrow(v_{0},\varphi_{v})$ in $X$, then
\eqref{eq4-9} implies \eqref{eq4-10}-\eqref{eq4-11}. The proof is complete.
$\Box$

\begin{remark}
The proof of (\ref{eq4-6}) also works in the deterministic case (see the Appendix).
\end{remark}

\begin{lemma}
\label{compact} Suppose that the conditions of Theorem \ref{thm3-3} hold. Then
the cocycle $\Xi$ is asymptotically compact.
\end{lemma}

\textbf{Proof.} Let $D\in\mathcal{D}$. It is sufficient to prove that for any
sequence $\{(u_{0}^{n},\varphi^{n})\}_{n\in\mathbb{N}}\subset D(\theta
_{-t_{n}}\omega)$, the sequence $\{\Xi(t_{n},\theta_{-t_{n}}\omega,(u_{0}%
^{n},\varphi^{n}))\}_{n\in\mathbb{N}}$ is relatively compact in $X$ as
$t_{n}\rightarrow\infty$. Recall that
\[%
\begin{split}
\Xi(t_{n},\theta_{-t_{n}}\omega,(u_{0}^{n},\varphi^{n}))  &  =\Psi
(t_{n},\theta_{-t_{n}}\omega,(u_{0}^{n}-\phi z_{\ast}(\theta_{-t_{n}}%
\omega),\mathcal{J}(\varphi^{n}-\phi z_{\ast}(\theta_{-t_{n}+\cdot}%
\omega))))\\
&  ~~+(\phi z_{\ast}(\omega),\phi z_{\ast}(\theta_{\cdot}\omega)).
\end{split}
\]
Hence, we only need to prove that the sequence $\{\Psi(t_{n},\theta_{-t_{n}%
}\omega,(u_{0}^{n}-\phi z_{\ast}(\theta_{-t_{n}}\omega),\mathcal{J}%
(\varphi^{n}-\phi z_{\ast}(\theta_{-t_{n}+\cdot}\omega))))\}$ possesses a
convergent subsequence in $X$. Observe that for $T>0$,
\[%
\begin{split}
&  ~\quad\Psi(t_{n},\theta_{-t_{n}}\omega,(u_{0}^{n}-\phi z_{\ast}%
(\theta_{-t_{n}}\omega),\varphi^{n}-\phi z_{\ast}(\theta_{-t_{n}+\cdot}%
\omega)))\\[0.4ex]
&  =\Psi(T,\theta_{-T}\omega,\Psi(t_{n}-T,\theta_{-t_{n}}\omega,(u_{0}%
^{n}-\phi z_{\ast}(\theta_{-t_{n}}\omega),\varphi^{n}-\phi z_{\ast}%
(\theta_{-t_{n}+\cdot}\omega))))\\[0.4ex]
&  =\Psi(T,\theta_{-T}\omega,\Psi(t_{n}-T,\theta_{-t_{n}+T}\theta_{-T}%
\omega,(u_{0}^{n}-\phi z_{\ast}(\theta_{-t_{n}+T}\theta_{-T}\omega
),\varphi^{n}-\phi z_{\ast}(\theta_{-t_{n}+T+\cdot}\theta_{-T}\omega
))))\\[0.4ex]
&  \subset\Psi(T,\theta_{-T}\omega,\tilde{B}(\theta_{-T}\omega)),
\end{split}
\]
for $t_{n}-T\geq\tilde{t}_{0}(\omega,D)$, where $\tilde{B}$ is the absorbing
ball of $\Psi$. Let
\[
\mathcal{Y}_{n}:=(\alpha^{n},\beta^{n})=\Psi(t_{n},\theta_{-t_{n}}%
\omega,(u_{0}^{n}-\phi z_{\ast}(\theta_{-t_{n}}\omega),\varphi^{n}-\phi
z_{\ast}(\theta_{-t_{n}+\cdot}\omega))).
\]
Then $(\alpha^{n},\beta^{n})=\Psi(T,\theta_{-T}\omega,\xi_{n}^{T})$, where
$\xi_{n}^{T}\in\tilde{B}(\theta_{-T}\omega)$. Let $(V^{n}(\cdot),V_{\cdot}%
^{n})$ be a sequence of solutions to problem \eqref{eq1-5} with initial
condition $\xi_{n}^{T}$ and $(V^{n}(T),V_{T}^{n})=(\alpha^{n},\beta^{n})$.
Since $\tilde{B}(\omega),\ \tilde{B}(\theta_{-T}\omega) $ are bounded in $X$,
by Lemma \ref{absorbing}, we can assume (up to a subsequence) that
$\mathcal{Y}_{n}\rightarrow\mathcal{Y=(\vartheta},\mathcal{\pi)}$, $\xi
_{n}^{T}\rightarrow\xi^{T}$ weakly in $X$ for every $\omega$.

It follows from Lemma \ref{lem4-5} that $(V^{n}(T),V_{T}^{n}(\cdot
))=\Psi(T,\theta_{-T}\omega,\xi_{n}^{T})$ satisfies
\eqref{eq4-6}-\eqref{eq4-8}. We deduce from the above convergence that
$\vartheta=V(T)$ in $H$ and $\pi=V_{T}$ in $L_{V}^{2}$, $\pi(s)=V_{T}(s)$ for
almost all $s\in(-\infty,0)$ and $\omega\in\Omega$. Also, in view of
\eqref{eq4-6}, we find that
\[
\alpha^{n}=V^{n}(T)\rightarrow V(T)=\vartheta\quad\mbox{in}~~H.
\]
Hence, in order to prove that $\mathcal{Y}_{n}\rightarrow\mathcal{Y}$ in $X$,
it remains to show that $\beta^{n}\rightarrow\pi$ in $L_{V}^{2}$ (up to a
subsequence). Notice that $\beta^{n}=V_{T}^{n}$ for all $T>0$ and $V_{T}=\pi$.
Since the family $\tilde{B}$ is tempered, we have that%
\begin{equation}
\lim_{T\rightarrow\infty}\ e^{-cT}\sup_{n}\left\Vert \xi_{n}^{T}\right\Vert
_{X}=0, \label{LimBT}%
\end{equation}
for any $c>0$. Thanks to \eqref{eq4-9}, we have, for each $T\in\mathbb{N}$,
\[%
\begin{split}
\limsup_{n\rightarrow\infty}\Vert\beta^{n}-\pi\Vert_{L_{V}^{2}}^{2}  &
=\limsup_{n\rightarrow\infty}\Vert V_{T}^{n}-V_{T}\Vert_{L_{V}^{2}}^{2}\\
&  \leq Ke^{-(\gamma-c)T}e^{-cT}\limsup_{n\rightarrow\infty}\left(  \Vert
\xi_{n}^{T}-\xi^{T}\Vert_{X}^{2}\right) \\
&  \leq\tilde{K}e^{-(\gamma-c)T},
\end{split}
\]
where $0<c<\gamma$ and the last inequality follows from (\ref{LimBT}). For
every $k>0$, there exists $T(k)$ such that for all $T\geq T(k)$,
\[
\limsup_{n\rightarrow\infty}\Vert\beta^{n}-\pi\Vert_{L_{V}^{2}}^{2}%
=\limsup_{n\rightarrow\infty}\Vert V_{T}^{n}-V_{T}\Vert_{L_{V}^{2}}^{2}%
\leq\frac{1}{k}.
\]
Taking $k\rightarrow\infty$ and using a diagonal argument, we obtain that
there exists a subsequence $\{\beta^{n_{k}}\}$ such that $\beta^{n_{k}%
}\rightarrow\pi$ in $L_{V}^{2}$ for all $\omega\in\Omega$. The proof of this
lemma is complete. $\Box$

A family of sets $K(\omega)$ is said to be measurable with respect to
$\mathcal{F}$, if for any $x\in X$ the mapping $\omega\mapsto dist(x,K(\omega
))$ is $\left(  \mathcal{F},\mathcal{B}(\mathbb{R})\right)  $-measurable.

We recall that the family of non-empty compact sets $\mathcal{A}%
=\mathcal{A}(\omega)\in\mathcal{D}$ is called a random attractor for the
cocycle $\Xi$, if $\mathcal{A}(\omega)$ is measurable with respect to
$\mathcal{F}$, it is invariant, that is,%
\[
\Xi(t,\omega,\mathcal{A}(\omega))=\mathcal{A}(\theta_{t}\omega)\text{ for all
}\omega\in\Omega,\ t\geq0,
\]
and it is pullback $\mathcal{D}$-attracting, that is, for any $D\in
\mathcal{D}$, it holds that%
\[
\lim_{t\rightarrow+\infty}dist(\Xi(t,\theta_{-t}\omega,D(\theta_{-t}%
\omega)),\mathcal{A}(\omega))=0,
\]
where $dist\left(  C_{1},C_{2}\right)  =\sup_{x\in C_{1}}\inf_{y\in C_{2}%
}\ \left\Vert x-y\right\Vert _{X}$ is the Hausdorff semidistance between sets
from $X$.

The existence and uniqueness of the random attractor $\mathcal{A}$ follow from
\cite[Proposition 2.10]{W2} (see also \cite{W1}, \cite{W3} for related
results) immediately based on Lemmas \ref{absorbing}, \ref{lem4-5} and
\ref{compact}. We observe that the radius $\rho\left(  \omega\right)  $ is
measurable, so it is easy to see that the family of closed balls $B(\omega)$
is measurable with respect to $\mathcal{F}$.

\begin{theorem}
Assume that \eqref{eq1-2}, \eqref{eq2-1} and $(h_{1})$-$(h_{2})$ hold, and
that $\phi\in V\cap H^{2}(\mathcal{O})\cap L^{2p}(\mathcal{O})$ is such that
$\Delta\phi\in L^{2p}(\mathcal{O})$. Let {$h\in H$} and $a(\cdot)$ be a
locally Lipschitz function. Then the cocycle $\Xi$ of problem \eqref{eq1-1}
has a unique random attractor $\mathcal{A}=\{\mathcal{A}(\omega):\omega
\in\Omega\}$ in $H$.
\end{theorem}

\section{Stochastic nonlocal PDEs with long time memory driven by colored
noise}

\label{s5} This section is devoted to discuss the approximations of stochastic
nonlocal PDEs with long time memory, namely, the following pathwise Wong-Zakai
approximated problem,
\begin{equation}%
\begin{cases}
\partial_{t}u_{\delta}-a(l(u_{\delta}))\Delta u_{\delta}-\displaystyle\int%
_{-\infty}^{t}k(t-s)\Delta u_{\delta}(s)ds+f(u_{\delta})=h+\phi\zeta_{\delta
}(\theta_{t}\omega),\\
u_{\delta}(t,x)=0,\\
u_{\delta}(\tau,x)=u_{0,\delta}(x),\\
u_{\delta}(t+\tau,x)=\varphi_{\delta}(t,x),\\
\end{cases}
\begin{aligned} &\mbox{in}~~\mathcal{O}\times (\tau,\infty),\\ &\mbox{on}~\partial\mathcal{O}\times(\tau,\infty),\\ &\mbox{in}~~\mathcal{O},\\ &\mbox{in}~~\mathcal{O}\times (-\infty,0),\\ \end{aligned} \label{eq5-1}%
\end{equation}
where $\zeta_{\delta}(\theta_{t}\omega)$ is the colored noise with correlation
time $\delta>0$, which is a stationary solution of the stochastic differential
equation
\[
d\zeta_{\delta}+\frac{1}{\delta}\zeta_{\delta}=\frac{1}{\delta}dW.
\]
This process satisfies%
\[
\lim_{t\rightarrow\pm\infty}\frac{\left\vert \zeta_{\delta}(\theta_{t}%
\omega)\right\vert }{t}=0\text{ for all }0<\delta\leq1,
\]%
\[
\lim_{\delta\rightarrow0^{+}}\sup_{t\in\lbrack\tau,\tau+T]}\left\vert \int%
_{0}^{t}\zeta_{\delta}(\theta_{s}\omega)ds-\omega(t)\right\vert =0,\ \forall
\tau\in\mathbb{R},\ T>0.
\]
For more details about colored noise see \cite{CW, CZ, G4} and the references
therein. For applications to stochastic Hamiltonian flows see \cite{Cui1},
\cite{Cui2}. The functions $a$, $f$ and $\phi$ fulfill the same assumptions as
in Section 2. Define a random variable,
\begin{equation}
v_{\delta}(t,\omega)=u_{\delta}(t,\omega)-\phi y_{\delta}(\theta_{t}\omega).
\label{eq5-2}%
\end{equation}
Recall that $y_{\delta}$ satisfies,
\begin{equation}
\frac{dy_{\delta}}{dt}=-\sigma y_{\delta}+\zeta_{\delta}(\theta_{t}\omega).
\label{eq5-3}%
\end{equation}
For almost all $\omega\in\Omega$, one special solution of \eqref{eq5-3} can be
represented by,
\[
Y_{\delta}(t,\omega)=e^{-\sigma t}\int_{-\infty}^{t}e^{\sigma s}\zeta
(\theta_{s}\omega)ds,
\]
which, in fact, can be rewritten as $Y_{\delta}(t,\omega)=y_{\delta}%
(\theta_{t}\omega)$. Here $y_{\delta}:\Omega\rightarrow\mathbb{R}$ is a
well-defined random variable given by $y_{\delta}(\omega):=\int_{-\infty}%
^{0}e^{\sigma s}\zeta_{\delta}(\theta_{s}\omega)ds$ and has the following properties.

\begin{lemma}
\cite[Lemma 3.2]{G4}\label{lem5-1} Let $y_{\delta}$ be the random variable
defined as above. Then the mapping
\begin{equation}
(t,\omega)\rightarrow y_{\delta}(\theta_{t}\omega)=e^{-\sigma t}\int_{-\infty
}^{t}e^{\sigma s}\zeta_{\delta}(\theta_{s}\omega)ds, \label{eq5-4}%
\end{equation}
is a stationary solution of \eqref{eq5-3} with continuous trajectories. In
addition, $\mathbb{E}(y_{\delta})=0$ and for every $\omega$,
\begin{equation}
\lim_{\delta\rightarrow0}y_{\delta}(\theta_{t}\omega)=z_{\ast}(\theta
_{t}\omega),\quad\mbox{uniformly on}~[\tau,\tau+T]~\mbox{with}~\tau
\in\mathbb{R},~T>0; \label{eq5-5}%
\end{equation}%
\begin{equation}
\lim_{t\rightarrow\pm\infty}\frac{|y_{\delta}(\theta_{t}\omega)|}{|t|}%
=0,\quad\mbox{uniformly for}~0<\delta\leq\tilde{\sigma}; \label{eq5-6}%
\end{equation}%
\begin{equation}
\lim_{t\rightarrow\pm\infty}\frac{1}{t}\int_{0}^{t}y_{\delta}(\theta_{s}%
\omega)ds=0,\quad\mbox{uniformly for}~0<\delta\leq\tilde{\sigma};
\label{eq5-7}%
\end{equation}
where $\tilde{\sigma}=\min\{1,\frac{1}{2\sigma}\}$ and $z_{\ast}(\theta
_{t}\omega)$ is given in Section 2.
\end{lemma}

\begin{remark}
\label{ConvergYdelta}It follows from (\ref{eq5-5})-(\ref{eq5-6}) that
$y_{\delta}(\theta_{\text{\textperiodcentered}}\omega)\rightarrow z_{\ast
}(\theta_{\text{\textperiodcentered}}\omega)$ in $L_{V}^{2}.$
\end{remark}

\begin{remark}
Throughout this paper, to simplify the computations, we take $\sigma=1$ in
\eqref{eq5-3}. Then the results of Lemma \ref{lem5-1} are true for $\sigma=1 $.
\end{remark}

Thus, it follows from \eqref{eq5-1}-\eqref{eq5-3} that, for $t>\tau$,
\begin{equation}%
\begin{cases}
\partial_{t}v_{\delta}-a(l(v_{\delta})+y_{\delta}(\theta_{t}\omega
)l(\phi))\Delta v_{\delta}-a(l(v_{\delta})+y_{\delta}(\theta_{t}\omega
)l(\phi))y_{\delta}(\theta_{t}\omega)\Delta\phi\\
\quad\displaystyle-\int_{-\infty}^{t}k(t-s)\Delta v_{\delta}(s)ds+f(v_{\delta
}+\phi y_{\delta}(\theta_{t}\omega))=h+\phi y_{\delta}(\theta_{t}%
\omega)+z_{k,\delta}^{\phi}(t,\omega),\\
v_{\delta}(t,x)=0,\\
v_{\delta}(\tau,x)=v_{0,\delta}(x):=u_{0,\delta}(x)-\phi y_{\delta}%
(\theta_{\tau}\omega),\\
v_{\delta}(t+\tau,x):=u_{\delta}(t+\tau,x)-\phi y_{\delta}(\theta_{t+\tau
}\omega)=\varphi_{\delta}(t,x)-\phi y_{\delta}(\theta_{t+\tau}\omega
):=\varphi_{v,\delta}(t,x),\\
\end{cases}
\begin{aligned} &\quad\\[2.2ex] &\mbox{in}~~\mathcal{O}\times (\tau,\infty),\\ &\mbox{on}~\partial\mathcal{O}\times(\tau,\infty),\\ &\mbox{in}~~\mathcal{O},\\ &\mbox{in}~~\mathcal{O}\times (-\infty,0),\\ \end{aligned} \label{eq5-8}%
\end{equation}
where $z_{k,\delta}^{\phi}(t)$ is a process defined by
\begin{equation}
z_{k,\delta}^{\phi}(t,\omega)=\int_{-\infty}^{t}k(t-s)y_{\delta}(\theta
_{s}\omega)\Delta\phi ds. \label{eq5-9}%
\end{equation}
To use Dafermos' transformation obtaining the well-posedness of problem
\eqref{eq5-8}, let us define the new variables,
\[
v_{\delta}^{t}(s,x,\omega)=v_{\delta}(t-s,x,\omega),\qquad s\geq0,
\]%
\[
\eta_{\delta}^{t}(s,x,\omega)=\int_{0}^{s}v_{\delta}^{t}(r,x,\omega
)dr=\int_{t-s}^{t}v_{\delta}(r,x,\omega)dr,\qquad s\geq0.
\]
Besides, assuming $k(\infty)=0$, a change of variable and a formal integration
by parts imply,
\[
\int_{-\infty}^{t}k(t-s)\nabla v_{\delta}(s)ds=-\int_{0}^{\infty}k^{\prime
}(s)\nabla\eta_{\delta}^{t}(s)ds.
\]
Setting $\mu(s)=-k^{\prime}(s)$, problem \eqref{eq5-8} turns into the
following system without delay,
\begin{equation}
\left\{
\begin{aligned} &\partial_tv_{\delta}-a(l(v_{\delta}+\phi y_\delta(\theta_t\omega)))\Delta v_{\delta}-a(l(v_{\delta}+\phi y_\delta(\theta_t\omega)))y_\delta(\theta_t\omega)\Delta \phi\\ &\quad-\displaystyle\int_{0}^{\infty}\mu(s)\Delta \eta_\delta^t(s)ds +f(v_{\delta}+\phi y_\delta(\theta_t\omega))=\phi y_\delta(\theta_t\omega)+z^\phi_{k,\delta}(t)+h,&&\mbox{in}~~\mathcal{O}\times (\tau,\infty),\\ &\partial_t\eta^t_\delta(s)=-\partial_s\eta^t_\delta(s)+v_{\delta}(t),&&\mbox{in}~~\mathcal{O}\times(\tau,\infty)\times\mathbb{R}^+,\\ &v_{\delta}(t,x)=\eta_\delta^t(x,s)=0,&&\mbox{on}~\partial\mathcal{O}\times(\tau,\infty)\times\mathbb{R}^+,\\ &v_{\delta}(\tau,x)=v_{0,\delta}(x):=u_{0,\delta}(x)-\phi y_\delta(\theta_{\tau}\omega),&&\mbox{in}~~\mathcal{O},\\ &\eta_\delta^{\tau}(s,x)=\eta_{0,\delta}(s,x),&&\mbox{in}~~\mathcal{O}\times\mathbb{R}^+,\\ \end{aligned}\right.
\label{eq5-10}%
\end{equation}
where
\[
\eta_{0,\delta}(s,x)(\omega)=\int_{\tau-s}^{\tau}v_{\delta}(r,x)dr=\int%
_{-s}^{0}(\varphi_{\delta}(r)-\phi y_{\delta}(\theta_{r+\tau}\omega
))dr=\int_{-s}^{0}\varphi_{v,\delta}(r)dr.
\]

The following result is proved exactly as Corollary \ref{cor3-2}.

\begin{corollary}
\label{cor3-2b} Assume that $(h_{1})$-$(h_{2})$ hold and $\phi\in V\cap
H^{2}(\mathcal{O})\cap L^{2p}(\mathcal{O})$. Then, for every $\omega\in\Omega$
and $\tau\in\mathbb{R}$, the operator $\mathcal{J}_{\omega,\tau}^{\delta
}:L_{V}^{2}\rightarrow L_{\mu}^{2}(\mathbb{R}^{+};V)$ defined by
\begin{equation}
(\mathcal{J}_{\omega,\tau}^{\delta}{\varphi})(s):=\int_{-s}^{0}\varphi
(r,\omega)dr-\int_{-s}^{0}y_{\delta}(\theta_{r+\tau}\omega)\phi dr=\mathcal{J}%
(\varphi-\phi y_{\delta}(\theta_{\tau+\cdot}\omega))(s), \label{eq3-0b}%
\end{equation}
is continuous. Additionally, there exists a positive constant $K_{\mu}$ which
is the same as in Lemma \ref{lem3-1} (which is also independent of $\delta$),
such that for any $\varphi\in L_{V}^{2}$, we have
\[
\Vert\mathcal{J}_{\omega,\tau}^{\delta}{\varphi}\Vert_{L_{\mu}^{2}%
(\mathbb{R}^{+};V)}^{2}\leq K_{\mu}\left\Vert \varphi-y_{\delta}%
(\theta_{\text{\textperiodcentered}+\tau}\omega)\phi\right\Vert _{L_{V}^{2}%
}^{2}\leq2K_{\mu}\left(  \Vert\varphi\Vert_{L_{V}^{2}}^{2}+\Vert y_{\delta
}(\theta_{\cdot+\tau}\omega)\phi\Vert_{L_{V}^{2}}^{2}\right)  .
\]

\end{corollary}

Since \eqref{eq5-10} can be viewed as a deterministic equation parameterized
by $\omega\in\Omega$, by the same procedures as in Theorems \ref{thm3-3} and
\ref{thm3-5}, we are able to prove the following results.

\begin{theorem}
\label{thm5-3} Assume that \eqref{eq1-2}, \eqref{eq2-1} and $(h_{1})$%
-$(h_{2})$ hold. Let $\phi\in V\cap H^{2}(\mathcal{O})\cap L^{2p}%
(\mathcal{O})$ be such that $\Delta\phi\in L^{2p}(\mathcal{O})$, let {$h\in
H$} and $a$ be a locally Lipschitz function. Then, for every $\tau
\in\mathbb{R}$ and $\omega\in\Omega$, it holds:

\begin{enumerate}
\item[$(i)$] For any initial value $v_{0,\delta}\in H$ and initial function
$\varphi_{\delta}\in L_{V}^{2}$, there exists a unique solution $(v_{\delta
},\eta_{\delta})$ to problem \eqref{eq5-10} in the weak sense {with initial
value $(v_{0,\delta},\eta_{0,\delta})$,} where $\eta_{0,\delta}=\mathcal{J}%
_{\omega,\tau}^{\delta}{\varphi_{\delta}}$, fulfilling
\[
v_{\delta}\in L^{\infty}(\tau,T;H)\cap L^{2}(\tau,T;V)\cap L^{2p}%
(\tau,T;L^{2p}(\mathcal{O})),\qquad\forall T>\tau;
\]%
\[
\eta_{\delta}\in L^{\infty}(\tau,T;L_{\mu}^{2}(\mathbb{R}^{+};V)),\qquad
\forall T>\tau.
\]
Furthermore, the solution $(v_{\delta},\eta_{\delta})$ of \eqref{eq5-10} is
continuous with respect to the initial value $(v_{0,\delta},\eta_{0,\delta})$
for all $t\in\lbrack\tau,T]$ in $\mathcal{H}$;

\item[$(ii)$] For any initial value $(v_{0,\delta},\eta_{0,\delta}%
)\in\mathcal{V}$, the unique solution $(v_{\delta},\eta_{\delta})$ to problem
\eqref{eq5-10} satisfies,
\[
v_{\delta}\in L^{\infty}(\tau,T;V)\cap L^{2}(\tau,T;V\cap H^{2}(\mathcal{O}%
)),\qquad\forall T>\tau;
\]
\[
\eta_{\delta}\in L^{\infty}(\tau,T;L_{\mu}^{2}(\mathbb{R}^{+};V\cap
H^{2}(\mathcal{O}))),\qquad\forall T>\tau.
\]
In addition, the solution $(v_{\delta},\eta_{\delta})$ of \eqref{eq5-10} is
continuous with respect to the initial value $(v_{0,\delta},\eta_{0,\delta})$
for all $t\in[\tau,T]$ in $\mathcal{V}$.
\end{enumerate}
\end{theorem}

Now, thanks to transformation \eqref{eq5-2}, we obtain the well-posedness of
problem \eqref{eq5-1}.

\begin{theorem}
\label{thm5-4} Assume \eqref{eq1-2}, \eqref{eq2-1} and $(h_{1})$-$(h_{2})$
hold, $\phi\in V\cap H^{2}(\mathcal{O})\cap L^{2p}(\mathcal{O})$ such that
$\Delta\phi\in L^{2p}(\mathcal{O})$. Let {$h\in H$} and $a$ be a locally
Lipschitz function. Then, for every $\tau\in\mathbb{R}$ and $\omega\in\Omega$,
it holds:

\begin{enumerate}
\item[$(i)$] For any initial value $u_{0,\delta}\in H$ and initial function
$\varphi_{\delta}\in L_{V}^{2}$, there exists a unique solution $u_{\delta}$
to problem \eqref{eq5-1} in the weak sense, fulfilling
\[
u_{\delta}\in L^{\infty}(\tau,T;H)\cap L^{2}(\tau,T;V)\cap L^{2p}%
(\tau,T;L^{2p}(\mathcal{O})),\qquad\forall T>\tau.
\]
Furthermore, the solution $u_{\delta}$ of \eqref{eq5-1} is continuous with
respect to the initial values $(u_{0,\delta}, \varphi_{\delta})$ for all
$t\in[\tau,T]$ in $H$;

\item[$(ii)$] For any initial value $u_{0,\delta}\in V$ and initial function
$\varphi_{\delta}\in L_{V\cap H^{2}(\mathcal{O})}^{2}$, the unique solution
$u_{\delta}$ to problem \eqref{eq5-1} satisfies,
\[
u_{\delta}\in L^{\infty}(\tau,T;V)\cap L^{2}(\tau,T;V\cap H^{2}(\mathcal{O}%
)),\qquad\forall T>\tau.
\]
In addition, the solution $u_{\delta}$ of \eqref{eq5-1} is continuous with
respect to the initial values $(u_{0,\delta},\varphi_{\delta})$ for all
$t\in[\tau,T]$ in $V$.
\end{enumerate}
\end{theorem}

Next, we can define a continuous cocycle in $X$ associated to the solutions of
problem \eqref{eq5-1}. Let $\tau=0$, $\Xi_{\delta}:\mathbb{R}^{+}\times
\Omega\times X\rightarrow X$ be a mapping defined, for every $t\in
\mathbb{R}^{+}$, $\omega\in\Omega$ and $(u_{0,\delta},\varphi_{\delta})\in X$,
by
\[
\Xi_{\delta}(t,\omega,(u_{0,\delta},\varphi_{\delta}))=(u_{\delta}%
(t;0,\omega,(u_{0,\delta},\varphi_{\delta})),u_{\delta,t}(\cdot;0,\omega
,(u_{0,\delta},\varphi_{\delta}))).
\]
Here and in the sequel, we denote $u_{\delta,t}(s)=u_{\delta}(t+s)$ for
$s\leq0$. Moreover, problem \eqref{eq5-10} also generates a random dynamical
system $\Phi_{\delta}$ on the phase space $H\times L_{\mu}^{2}(\mathbb{R}%
^{+};V)$, which is defined, for every $t\in\mathbb{R}^{+}$, $\omega\in\Omega$
and $(v_{0,\delta},\eta_{0,\delta})\in H\times L_{\mu}^{2}(\mathbb{R}^{+};V)$,
by
\[
\Phi_{\delta}(t,\omega,(v_{0,\delta},\eta_{0,\delta}))=(v_{\delta}%
(t;0,\omega,(v_{0,\delta},\eta_{0,\delta})),\eta_{\delta}^{t}(\cdot
;0,\omega,(v_{0,\delta},\eta_{0,\delta}))),
\]
where the right-hand side of the above equality denotes the solution of
\eqref{eq5-10} for $\tau=0$, the initial value $(v_{0,\delta},\eta_{0,\delta
})\in H\times L_{\mu}^{2}(\mathbb{R}^{+};V)$ and $\eta_{0,\delta}$ is given in
\eqref{eq5-10}. Thanks to Dafermos' transformation, we can obtain a random
dynamical system $\Psi_{\delta}:\mathbb{R}^{+}\times\Omega\times X\rightarrow
X$ generated by \eqref{eq5-8} which is given, for every $t\in\mathbb{R}^{+}$,
$\omega\in\Omega$ and $(v_{0,\delta},\psi_{\delta})\in X$, by
\[
\Psi_{\delta}(t,\omega,(v_{0,\delta},\psi_{\delta}))=(v_{\delta}%
(t;0,\omega,(v_{0,\delta},(\mathcal{J}\psi_{\delta}))),v_{\delta,t}%
(\cdot;0,\omega,(v_{0,\delta},(\mathcal{J}\psi_{\delta})))).
\]
Then, on account of the random transformation \eqref{eq5-2}, for
$(u_{0,\delta},\varphi_{\delta})\in X$, we deduce
\begin{align}
&  \quad\Xi_{\delta}(t,\omega,(u_{0,\delta},\varphi_{\delta}%
))\nonumber\label{eq5-11}\\
&  =(u_{\delta}(t;0,\omega,(u_{0,\delta},\varphi_{\delta})),u_{\delta,t}%
(\cdot;0,\omega,(u_{0,\delta},\varphi_{\delta})))\nonumber\\
&  =\big(v_{\delta}(t;0,\omega,(u_{0,\delta}-\phi y_{\delta}(\omega
),(\mathcal{J}_{\omega,0}^{\delta}{\varphi_{\delta}})))+\phi y_{\delta}%
(\theta_{t}\omega),\\
&  ~~~~v_{\delta,t}(\cdot;0,\omega,(u_{0,\delta}-\phi y_{\delta}%
(\omega),(\mathcal{J}_{\omega,0}^{\delta}{\varphi_{\delta}})))+\phi y_{\delta
}(\theta_{t+\cdot}\omega)\big)\nonumber\\
&  =\Psi_{\delta}(t,\omega,(u_{0,\delta}-\phi y_{\delta}(\omega),{\varphi
_{v,\delta}}))+(\phi y_{\delta}(\theta_{t}\omega),\phi y_{\delta}%
(\theta_{t+\cdot}\omega)).\nonumber
\end{align}
Therefore, similar to Section 4, it is easy to check that $\Xi_{\delta}$ and
$\Psi_{\delta}$ are conjugated.

We will prove now the existence of absorbing sets for the cocycle $\Xi
_{\delta}$.

\begin{lemma}
\label{lem5-5} Under the assumptions of Theorem \ref{thm5-4}, there exists
$B_{\delta}\in\mathcal{D}$ which is $\mathcal{D}$-pullback absorbing for the
RDS $\Xi_{\delta}$. In other words, for any given $\omega\in\Omega$ and
$D_{\delta}\in\mathcal{D}$, there exists $t_{0,\delta}:=t_{0,\delta}%
(\omega,D_{\delta})\geq0$, such that
\[
\Xi_{\delta}(t,\theta_{-t}\omega,D_{\delta}(\theta_{-t}\omega))\subset
B_{\delta}(\omega),\ \quad\text{for all}\ t\geq t_{0,\delta}(\omega,D_{\delta
}),
\]
where $B_{\delta}(\omega)$ is the ball in $X$ centered $0$ with radius
$\rho_{\delta}(\omega)$,
\begin{equation}
\rho_{\delta}^{2}(\omega)=1+2K_{2}\int_{-\infty}^{0}e^{\gamma s}%
\Theta_{1,\delta}(\theta_{s}\omega)ds+2|\phi|^{2}|y_{\delta}(\omega
)|^{2}+2\Vert\phi y_{\delta}(\theta_{\cdot}\omega)\Vert_{L_{V}^{2}}^{2},
\label{thetadelta}%
\end{equation}
{where $\Phi_{1,\delta}(\omega)$ is defined in \eqref{theta2} below.} In
addition, for every $\omega\in\Omega$,
\begin{equation}
\lim_{\delta\rightarrow0}\rho_{\delta}^{2}(\omega)=1+2K_{2}\int_{-\infty}%
^{0}e^{\gamma s}\Theta_{1}(\theta_{s}\omega)ds+2|\phi|^{2}|z_{\ast}%
(\omega)|^{2}+2\Vert\phi z_{\ast}(\theta_{\cdot}\omega)\Vert_{L_{V}^{2}}%
^{2}:=\rho^{2}(\omega). \label{rhoconvergence}%
\end{equation}

\end{lemma}

\textbf{Proof.} Let us first pick $(u_{0,\delta},\varphi_{\delta})\in
D_{\delta}$ and thanks to \eqref{eq5-11}, we have
\begin{align*}
&  \quad\Xi_{\delta}(t,\omega,(u_{0,\delta},\varphi_{\delta}))\\[0.4ex]
&  =(v_{\delta}(t;0,\omega,(u_{0,\delta}-\phi y_{\delta}(\omega),(\mathcal{J}%
_{\omega,0}{\varphi_{\delta}}))),v_{\delta,t}(\cdot;0,\omega,(u_{0,\delta
}-\phi y_{\delta}(\omega),(\mathcal{J}_{\omega,0}{\varphi_{\delta}%
})))\\[0.4ex]
&  ~~+(\phi y_{\delta}(\theta_{t}\omega),\phi y_{\delta}(\theta_{t+\cdot
}\omega))\\[0.4ex]
&  =\Psi_{\delta}(t,\omega,(u_{0,\delta}-\phi y_{\delta}(\omega),{\varphi
_{v,\delta}}))+(\phi y_{\delta}(\theta_{t}\omega),\phi y_{\delta}%
(\theta_{t+\cdot}\omega)).
\end{align*}
For the sake of simplicity, denote by $z_{\delta}(t):=(v_{\delta}%
(t),\eta_{\delta}^{t}(s))$ the solution to \eqref{eq5-10} with initial value
$(v_{0,\delta},\eta_{0,\delta})=(u_{0,\delta}-\phi y_{\delta}(\omega
),(\mathcal{J}_{\omega,0}{\varphi_{\delta}}))$. Now, for every $\omega
\in\Omega$, we multiply the first equation of \eqref{eq5-10} by $v_{\delta
}(t)$ in $H$ and the second equation of \eqref{eq5-10} by $\eta_{\delta}^{t}$
in $L_{\mu}^{2}(\mathbb{R}^{+};V)$, respectively. Then, by means of the same
estimates as in the proof of Theorem \ref{thm3-3} (cf. \eqref{eq3-11}) and the
Poincar\'{e} inequality, we obtain
\begin{equation}
\frac{d}{dt}\Vert z_{\delta}(t)\Vert_{\mathcal{H}}^{2}+\frac{m\lambda_{1}}%
{2}|v_{\delta}(t)|^{2}+\frac{m}{2}\Vert v_{\delta}(t)\Vert^{2}+2((\eta
_{\delta}^{t}(s),(\eta_{\delta}^{t}(s))^{\prime}))_{\mu}+\frac{f_{0}}{2D}\Vert
v_{\delta}(t)\Vert_{2p}^{2p}\leq\Theta_{1,\delta}(\theta_{t}\omega),
\label{eq5-12}%
\end{equation}
where we have used the notation
\begin{equation}
\Theta_{1,\delta}(\omega)=\frac{4M^{2}}{m}|y_{\delta}(\omega)|^{2}\Vert
\phi\Vert^{2}+\frac{4}{m\lambda_{1}}|y_{\delta}(\omega)|^{2}|\phi|^{2}%
+2\alpha|\mathcal{O}|+{C_{1,\delta}(\omega)(1+}\Vert\phi\Vert_{2p}^{2p}%
)+\frac{C_{2,\delta}^{2}(\omega)}{m}\Vert\phi\Vert^{2}+\frac{{4}}{m\lambda
_{1}}\left\vert h\right\vert ^{2}, \label{theta2}%
\end{equation}
and $C_{1,\delta}(\omega),\ C_{2,\delta}(\omega),\ D$ are the same as the ones
in (\ref{eq3-11})-\eqref{theta1} but replacing $z_{\ast}(\omega)$ by
$y_{\delta}(\omega)$. Taking into account condition $(h_{2})$ and {recalling
that $0<\gamma<\min\{\frac{m\lambda_{1}}{2},\varpi\}$}, the above inequality,
together with \eqref{eq5-12}, implies that
\begin{equation}
\frac{d}{dt}\Vert z_{\delta}(t)\Vert_{\mathcal{H}}^{2}+\gamma\Vert z_{\delta
}(t)\Vert_{\mathcal{H}}^{2}+\frac{m}{2}\Vert v_{\delta}(t)\Vert^{2}%
+\frac{f_{0}}{2D}\Vert v_{\delta}(t)\Vert_{2p}^{2p}\leq\Theta_{1,\delta
}(\theta_{t}\omega). \label{eq5-13}%
\end{equation}
Next, multiplying the above inequality by $e^{\gamma t}$ and integrating over
$(0,t)$, neglecting the last term on the left hand side of \eqref{eq5-13}, by
doing similar computations as in (\ref{eq4-2})-(\ref{eq4-3}), we find
\begin{equation}%
\begin{split}
\frac{m}{2}\Vert v_{\delta,t}\Vert_{L_{V}^{2}}^{2}  &  =\frac{m}{2}%
\int_{-\infty}^{0}e^{-\gamma(t-s)}\Vert\varphi_{v,\delta}(s)\Vert^{2}%
ds+\frac{m}{2}\int_{0}^{t}e^{-\gamma(t-s)}\Vert v_{\delta}(s)\Vert
^{2}ds\\[0.8ex]
&  \leq\frac{m}{2}e^{-\gamma t}\Vert\varphi_{v,\delta}\Vert_{L_{V}^{2}}%
^{2}+\Vert z_{0,\delta}\Vert_{\mathcal{H}}^{2}e^{-\gamma t}+\int_{0}%
^{t}e^{-\gamma(t-s)}\Theta_{1,\delta}(\theta_{s}\omega)ds.
\end{split}
\label{eq5-14}%
\end{equation}
On account of Corollary \ref{cor3-2b}, we have
\begin{equation}%
\begin{split}
\Vert z_{0,\delta}\Vert_{\mathcal{H}}^{2}  &  =|v_{0,\delta}|^{2}%
+\Vert\mathcal{J}(\varphi_{\delta}-\phi y_{\delta}(\theta_{\cdot}\omega
))\Vert_{L_{\mu}^{2}(\mathbb{R}^{+};V)}^{2}\\
&  \leq|u_{0,\delta}-\phi y_{\delta}(\omega)|^{2}+2K_{\mu}\left(  \Vert
\varphi_{\delta}-\phi y_{\delta}(\theta_{\cdot}\omega)\Vert_{L_{V}^{2}}%
^{2}\right)  .
\end{split}
\label{eq5-15}%
\end{equation}
Hence, collecting \eqref{eq5-14}-\eqref{eq5-15}, we arrive at
\[%
\begin{split}
&  ~\quad\Vert\Psi_{\delta}(t,\omega,(u_{0,\delta}-\phi y_{\delta}%
(\omega),\mathcal{J}(\varphi_{\delta}-\phi y_{\delta}(\theta_{\cdot}%
\omega))))\Vert_{X}^{2}=|v_{\delta}(t)|^{2}+\Vert v_{\delta,t}\Vert_{L_{V}%
^{2}}^{2}\leq\Vert z_{\delta}(t)\Vert_{\mathcal{H}}^{2}+\Vert v_{\delta
,t}\Vert_{L_{V}^{2}}^{2}\\[0.8ex]
&  \leq K_{1}e^{-\gamma t}\left(  |u_{0,\delta}-\phi y_{\delta}(\omega
)|^{2}+\Vert\varphi_{\delta}-\phi y_{\delta}(\theta_{\cdot}\omega)\Vert
_{L_{V}^{2}}^{2}\right)  +K_{2}\int_{0}^{t}e^{-\gamma(t-s)}\Theta_{1,\delta
}(\theta_{s}\omega)ds,
\end{split}
\]
where $K_{1}$, $K_{2}>0$ are the same constants as in \eqref{eq48}. Replacing
$\omega$ by $\theta_{-t}\omega$ in the above inequality, we obtain
\[%
\begin{split}
&  ~\quad\Vert\Psi_{\delta}(t,\theta_{-t}\omega,(u_{0,\delta}-\phi y_{\delta
}(\theta_{-t}\omega),\mathcal{J}(\varphi_{\delta}-\phi y_{\delta}%
(\theta_{-t+\cdot}\omega))))\Vert_{X}^{2}\\[0.4ex]
&  \leq K_{1}e^{-\gamma t}\left(  |u_{0,\delta}-\phi y_{\delta}(\theta
_{-t}\omega)|^{2}+\Vert\varphi_{\delta}-\phi y_{\delta}(\theta_{-t+\cdot
}\omega)\Vert_{L_{V}^{2}}^{2}\right)  +K_{2}\int_{-\infty}^{0}e^{\gamma
s}\Theta_{1,\delta}(\theta_{s}\omega)ds.
\end{split}
\]
Therefore, for any $(u_{0,\delta},\varphi_{\delta})\in D_{\delta}(\theta
_{-t}\omega)$, we have
\[%
\begin{split}
&  \quad~\Vert\Xi_{\delta}(t,\theta_{-t}\omega,(u_{0,\delta},\varphi_{\delta
}))\Vert_{X}^{2}\\[0.8ex]
&  =\Vert\Psi_{\delta}(t,\theta_{-t}\omega,(u_{0,\delta}-\phi y_{\delta
}(\theta_{-t}\omega),\mathcal{J}(\varphi_{\delta}-\phi y_{\delta}%
(\theta_{-t+\cdot}\omega))))+(\phi y_{\delta}(\omega),\phi y_{\delta}%
(\theta_{\cdot}\omega))\Vert_{X}^{2}\\[0.4ex]
&  \leq2K_{1}e^{-\gamma t}\Vert\widetilde{D}_{\delta}(\theta_{-t}\omega
)\Vert^{2}+2K_{2}\int_{-\infty}^{0}e^{\gamma s}\Theta_{1,\delta}(\theta
_{s}\omega)ds+2|\phi|^{2}|y_{\delta}(\omega)|^{2}+2\Vert\phi y_{\delta}%
(\theta_{\cdot}\omega)\Vert_{L_{V}^{2}}^{2}.
\end{split}
\]
Consequently, let
\[
\rho_{\delta}^{2}(\omega)=1+2K_{2}\int_{-\infty}^{0}e^{\gamma s}%
\Theta_{1,\delta}(\theta_{s}\omega)ds+2|\phi|^{2}|y_{\delta}(\omega
)|^{2}+2\Vert\phi y_{\delta}(\theta_{\cdot}\omega)\Vert_{L_{V}^{2}}^{2}.
\]
Taking into account the temperedness of $\widetilde{D}_{\delta}$ and the
properties of the Ornstein-Uhlenbeck process, it is straightforward to check
that the set,
\[
B_{\delta,\rho}(\omega)=\{(u_{0,\delta},\varphi_{\delta})\in X:\Vert
(u_{0,\delta},\varphi_{\delta})\Vert_{X}\leq\rho_{\delta}(\omega)\},
\]
is tempered, i.e., $B_{\delta,\rho}(\omega)$ belongs to $\mathcal{D}$, and
that it is pullback absorbing for this universe $\mathcal{D}$.

We now prove \eqref{rhoconvergence}. To this end, suitable estimates of the
functions $y_{\delta},C_{1,\delta},C_{2,\delta}$ are needed. We see first that
(\ref{eq5-6}) implies that there exist $r<0$ and $\delta_{0}>0$, such that for
all $0<\delta<\delta_{0}$,
\begin{equation}
|y_{\delta}(\theta_{t}\omega)|\leq|t|,\quad\forall t\leq r. \label{eq518}%
\end{equation}
Let us analyze the functions $C_{1,\delta}(\omega)$ and $C_{2,\delta}(\omega
)$. It is clear from (\ref{eq5-5}) that {when $\delta\rightarrow0$,}%
\begin{equation}
{C_{1,\delta}(\theta_{t}\omega)=\widetilde{C}_{2}({1+}}|y_{\delta}(\theta
_{t}\omega)|^{4p^{2}})\rightarrow{\widetilde{C}_{2}({1+}}|z_{\ast}(\theta
_{t}\omega)|^{4p^{2}})={C_{1}(\theta_{t}\omega)}\text{ uniformly in }%
[\tau,T],\ \tau<T. \label{ConvergC1delta}%
\end{equation}
Also, (\ref{eq518}) gives that%
\begin{equation}
\left\vert C_{1,\delta}(\theta_{t}\omega)\right\vert \leq\widetilde{C}%
_{2}\left(  1+\left\vert t\right\vert ^{4p^{2}}\right)  ,\ \forall t\leq
r,\ 0<\delta<\delta_{0}. \label{eq518b}%
\end{equation}
Furthermore, let us consider the function,%
\[
C_{2,\delta}(\theta_{t}\omega)=2\left(  M_{1}\int_{0}^{1}|y_{\delta}%
(\theta_{t-s}\omega)|ds+\frac{\mu(1)e^{\varpi}}{\varpi}\int_{1}^{\infty
}e^{-\varpi s}|y_{\delta}(\theta_{t-s}\omega)|ds\right)  .
\]
From (\ref{eq5-5})-(\ref{eq5-6}), it is easy to see that%
\begin{equation}
{C_{2,\delta}(\theta_{t}\omega)}\rightarrow{C_{2}(\theta_{t}\omega)}\text{
uniformly in }[\tau,T],\ \tau<T,~~\mbox{as}~{\delta\rightarrow0}.
\label{ConvergC2delta}%
\end{equation}
Finally, (\ref{eq518}) implies that for $t\leq r,\ 0<\delta<\delta_{0}$,
\begin{align}
\left\vert C_{2,\delta}(\theta_{t}\omega)\right\vert  &  \leq2\left(
M_{1}\int_{0}^{1}\left\vert t-s\right\vert ds+\frac{\mu(1)e^{\varpi}}{\varpi
}\int_{1}^{\infty}e^{-\varpi s}|t-s|ds\right) \nonumber\\
&  \leq2\left(  M_{1}(\left\vert t\right\vert +1)+\frac{\mu(1)e^{\varpi}%
}{\varpi}\left(  \frac{1}{\varpi}\left\vert t\right\vert +\int_{1}^{\infty
}e^{-\varpi s}sds\right)  \right)  . \label{eq518c}%
\end{align}
On the one hand, notice that
\[%
\begin{split}
&  ~~\quad K_{2}\int_{-\infty}^{0}e^{\gamma s}\Theta_{1,\delta}(\theta
_{s}\omega)ds+|\phi|^{2}|y_{\delta}(\omega)|^{2}+\Vert\phi y_{\delta}%
(\theta_{\cdot}\omega)\Vert_{L_{V}^{2}}^{2}\\
&  =K_{2}\int_{-\infty}^{r}e^{\gamma s}\Theta_{1,\delta}(\theta_{s}%
\omega)ds+K_{2}\int_{r}^{0}e^{\gamma s}\Theta_{1,\delta}(\theta_{s}%
\omega)ds+|\phi|^{2}|y_{\delta}(\omega)|^{2}+\Vert\phi y_{\delta}%
(\theta_{\cdot}\omega)\Vert_{L_{V}^{2}}^{2}.
\end{split}
\]
For all $0<\delta<\delta_{0}$, we obtain by \eqref{eq518} and (\ref{eq518c})
that%
\begin{equation}%
\begin{split}
&  \qquad e^{\gamma s}\Theta_{1,\delta}(\theta_{s}\omega)ds\\
&  =e^{\gamma s}\left(  \frac{4M^{2}}{m}|y_{\delta}(\theta_{s}\omega
)|^{2}\Vert\phi\Vert^{2}+\frac{4}{m\lambda_{1}}|y_{\delta}(\theta_{s}%
\omega)|^{2}|\phi|^{2}+2\alpha|\mathcal{O}|+C_{1,\delta}(\theta_{s}%
\omega)(1+\Vert\phi\Vert_{2p}^{2p})\right. \\
&  ~~\left.  +\frac{C_{2,\delta}^{2}(\theta_{s}\omega)}{m}\Vert\phi\Vert
^{2}+\frac{{4}}{m\lambda_{1}}\left\vert h\right\vert ^{2}\right) \\
&  \leq e^{\gamma s}\left(  \frac{4M^{2}}{m}|s|^{2}\Vert\phi\Vert^{2}+\frac
{4}{m\lambda_{1}}|s|^{2}|\phi|^{2}+2\alpha|\mathcal{O}|+\widetilde{C}%
_{2}(1+\left\vert s\right\vert ^{4p^{2}})(1+\Vert\phi\Vert_{2p}^{2p})\right.
\\
&  ~~\left.  +\frac{2}{m}\left(  M_{1}(\left\vert s\right\vert +1)+\frac
{\mu(1)e^{\varpi}}{\varpi}\left(  \frac{1}{\varpi}\left\vert s\right\vert
+\int_{1}^{\infty}e^{-\varpi l}ldl\right)  \right)  \Vert\phi\Vert^{2}%
+\frac{{4}}{m\lambda_{1}}\left\vert h\right\vert ^{2}\right)  :=g(s),\text{
for }s\leq r,
\end{split}
\label{gEst}%
\end{equation}
where $g\in L^{1}\left(  -\infty,r\right)  $. By means of the above estimates,
\eqref{eq5-5}, (\ref{ConvergC1delta}), (\ref{ConvergC2delta}) and the
dominated convergence theorem, we find that
\begin{equation}
\lim_{\delta\rightarrow0}\int_{-\infty}^{r}e^{\gamma s}\Theta_{1,\delta
}(\theta_{s}\omega)ds=\int_{-\infty}^{r}e^{\gamma s}\Theta_{1}(\theta
_{s}\omega)ds. \label{eq519}%
\end{equation}
On the other hand, by \eqref{eq5-5}, (\ref{ConvergC1delta}) and
(\ref{ConvergC2delta}), we infer that%
\[%
\begin{split}
&  \quad~~\lim_{\delta\rightarrow0}\int_{r}^{0}e^{\gamma s}\Theta_{1,\delta
}(\theta_{s}\omega)ds\\
&  =\lim_{\delta\rightarrow0}\int_{r}^{0}e^{\gamma s}\left(  \frac{4M^{2}}%
{m}|y_{\delta}(\theta_{s}\omega)|^{2}\Vert\phi\Vert^{2}+\frac{4}{m\lambda_{1}%
}|y_{\delta}(\theta_{s}\omega)|^{2}|\phi|^{2}+2\alpha|\mathcal{O}%
|+C_{1,\delta}(\theta_{s}\omega)(1+\Vert\phi\Vert_{2p}^{2p})\right. \\
&  ~~\left.  +\frac{C_{2,\delta}^{2}(\theta_{s}\omega)}{m}\Vert\phi\Vert
^{2}+\frac{{4}}{m\lambda_{1}}\left\vert h\right\vert ^{2}\right)  ds\\
&  =\int_{r}^{0}e^{\gamma s}\left(  \frac{4M^{2}}{m}|z_{\ast}(\theta_{s}%
\omega)|^{2}\Vert\phi\Vert^{2}+\frac{4}{m\lambda_{1}}|z_{\ast}(\theta
_{s}\omega)|^{2}|\phi|^{2}+2\alpha|\mathcal{O}|+C_{1}(\theta_{s}%
\omega)(1+\Vert\phi\Vert_{2p}^{2p})\right. \\
&  ~~\left.  +\frac{C_{2}^{2}(\theta_{s}\omega)}{m}\Vert\phi\Vert^{2}%
+\frac{{4}}{m\lambda_{1}}\left\vert h\right\vert ^{2}\right)  ds\\
&  =\int_{r}^{0}e^{\gamma s}\Theta_{1}(\theta_{s}\omega)ds.
\end{split}
\]
Combining the above inequality with \eqref{eq519}, we deduce
\begin{equation}
\lim_{\delta\rightarrow0}\int_{-\infty}^{0}e^{\gamma s}\Theta_{1,\delta
}(\theta_{s}\omega)ds=\int_{-\infty}^{0}e^{\gamma s}\Theta_{1}(\theta
_{s}\omega)ds. \label{eq520}%
\end{equation}
Moreover, (\ref{eq5-5})-(\ref{eq5-6}) and the dominated convergence theorem
imply that,
\begin{equation}
\lim_{\delta\rightarrow0}\left(  |\phi|^{2}|y_{\delta}(\omega)|^{2}+\Vert\phi
y_{\delta}(\theta_{\cdot}\omega)\Vert_{L_{V}^{2}}^{2}\right)  =|\phi
|^{2}|z_{\ast}(\omega)|^{2}+\Vert\phi z_{\ast}(\theta_{\cdot}\omega
)\Vert_{L_{V}^{2}}^{2}. \label{eq521}%
\end{equation}
Therefore, \eqref{eq520}-\eqref{eq521} finish the proof of this lemma. $\Box$

\begin{remark}
\label{rem5-6} It follows from the proof of Lemma \ref{lem5-5} that, under the
assumptions of Theorem \ref{thm5-3}, there exists $\tilde{B}_{\delta}%
\in\mathcal{D}$ which is $\mathcal{D}$-pullback absorbing for the RDS
$\Psi_{\delta}$. In other words, for any given $\omega\in\Omega$ and
$D_{\delta}\in\mathcal{D}$, there exists $\tilde{t}_{0,\delta}:=\tilde
{t}_{0,\delta}(\omega,D_{\delta})\geq0$, such that
\[
\Psi_{\delta}(t,\theta_{-t}\omega,D(\theta_{-t}\omega))\subset\tilde
{B}_{\delta}(\omega),\ \quad\text{for all}\ t\geq\tilde{t}_{0,\delta}%
(\omega,D_{\delta}).
\]

\end{remark}

Next, by means of the same procedure and estimates as in the proof of Lemma
\ref{lem4-5}, before stating the asymptotic compactness of the cocycle
$\Xi_{\delta}$, we first need the following auxiliary lemma. Since the details
are similar to those in Lemma \ref{lem4-5}, we omit the proof here.

\begin{lemma}
\label{lem5-7} Assume the hypotheses in Theorem \ref{thm5-3} hold. Let
{$\{v_{0,\delta}^{n},\varphi_{v,\delta}^{n}\}$} be a sequence such that
$(v_{0,\delta}^{n},\varphi_{v,\delta}^{n})\rightarrow(v_{0,\delta}%
,\varphi_{v,\delta})$ weakly in $X$ as $n\rightarrow\infty$. Then, for every
$\omega\in\Omega$, $\Psi_{\delta}(t,\omega,(v_{0,\delta}^{n},\varphi
_{v,\delta}^{n}))=(v_{\delta}^{n}(t),v_{\delta,t}^{n}($\textperiodcentered$))$
fulfills:
\begin{equation}
v_{\delta}^{n}\rightarrow v_{\delta}\quad\mbox{in}\quad
C([r,T];H)~~\mbox{for all}~~0<r<T; \label{eq5-17}%
\end{equation}%
\begin{equation}
v_{\delta}^{n}\rightarrow v_{\delta}\quad\mbox{weakly in}\quad L^{2}%
(0,T;V)~~\mbox{for all}~~T>0; \label{eq5-18}%
\end{equation}%
\begin{equation}
v_{\delta}^{n}\rightarrow v_{\delta}\quad\mbox{in}\quad L^{2}%
(0,T;H)~~\mbox{for all}~~T>0; \label{eq5-19}%
\end{equation}%
\begin{equation}
\limsup_{n\rightarrow\infty}\Vert v_{\delta,t}^{n}-v_{\delta,t}\Vert
_{L_{V}^{2}}^{2}\leq Ke^{-\gamma t}\limsup_{n\rightarrow\infty}\left(
|v_{0,\delta}^{n}-v_{0,\delta}|^{2}+\Vert\varphi_{v,\delta}^{n}-\varphi
_{v,\delta}\Vert_{L_{V}^{2}}^{2}\right)  ~\mbox{for all}~~t\geq0,
\label{eq5-20}%
\end{equation}
where {$K=(1+\frac{2K_{\mu}}{m}+\frac{1}{m})$ and }$(v_{\delta}(t),v_{\delta
,t}($\textperiodcentered$))=\Psi_{\delta}(t,\omega,(v_{0,\delta}%
,\varphi_{v,\delta}))$. Moreover, if $(v_{0,\delta}^{n},\varphi_{v,\delta}%
^{n})\rightarrow(v_{0,\delta},\varphi_{v,\delta})$ strongly in $X$ as
$n\rightarrow\infty$, then
\begin{equation}
v_{\delta}^{n}\rightarrow v_{\delta}\quad\mbox{in}\quad L^{2}%
(0,T;V)~\mbox{for all}~~T>0; \label{eq5-21}%
\end{equation}%
\begin{equation}
v_{\delta,t}^{n}\rightarrow v_{\delta,t}\quad\mbox{in}\quad L_{V}%
^{2}~~\mbox{for
all}~~t\geq0. \label{eq5-22}%
\end{equation}

\end{lemma}

\begin{lemma}
\label{lem5-8} Suppose that the conditions of Theorem \ref{thm5-3} hold, then
the cocycle $\Xi_{\delta}$ is asymptotically compact.
\end{lemma}

As a consequence of Lemmas \ref{lem5-5}, \ref{lem5-7} and \ref{lem5-8}, we can
ensure the existence of a {random} attractor to problem \eqref{eq5-1}.

\begin{theorem}
Assume that \eqref{eq1-2}, \eqref{eq2-1} and $(h_{1})$-$(h_{2})$ hold,
$\phi\in V\cap H^{2}(\mathcal{O})\cap L^{2p}(\mathcal{O})$ is such that
$\Delta\phi\in L^{2p}(\mathcal{O})$. Let $h\in H$ and $a$ be a locally
Lipschitz function. Then the cocycle $\Xi_{\delta}$ of problem \eqref{eq5-1}
has a unique random attractor $\mathcal{A}_{\delta}=\{\mathcal{A}_{\delta
}(\omega):\omega\in\Omega\}$ in $H$.
\end{theorem}

Let us define the family $\overline{B}$ by
\[
\overline{B}(\omega)=\cup_{0<\delta\leq\delta_{0}}B_{\delta}(\omega),
\]
where $\delta_{0}\leq\widetilde{\sigma}.$ Then, the following lemma holds true.

\begin{lemma}
\label{BBarraTemp}$\overline{B}$ is tempered for some $\delta_{0}>0.$
\end{lemma}

\textbf{Proof. }It is enough to check that for any $c>0$ and $\omega\in\Omega
$, we have%
\begin{equation}
\lim_{t\rightarrow\infty}e^{-ct}\sup_{0<\delta\leq\delta_{0}}\rho_{\delta}%
^{2}(\theta_{-t}\omega)=0,\ \label{LimRadius}%
\end{equation}
where $\rho_{\delta}\left(  \omega\right)  $ is the radius of the absorbing
ball $B_{\delta}(\omega)$, defined in (\ref{thetadelta}). First, from
(\ref{eq518}), it is clear that%
\[
\lim_{t\rightarrow\infty}\ e^{-ct}2|\phi|^{2}\sup_{0<\delta\leq\delta_{0}%
}|y_{\delta}(\theta_{-t}\omega)|^{2}=0.
\]
Since%
\begin{align*}
2\Vert\phi y_{\delta}(\theta_{-t+\cdot}\omega)\Vert_{L_{V}^{2}}^{2}  &
=2\left\Vert \phi\right\Vert ^{2}\int_{-\infty}^{0}e^{\gamma s}\left\vert
y_{\delta}(\theta_{-t+s}\omega)\right\vert ^{2}ds\\
&  \leq4\left\Vert \phi\right\Vert ^{2}\left(  t^{2}\int_{-\infty}%
^{0}e^{\gamma s}ds+\int_{-\infty}^{0}e^{\gamma s}\left\vert s\right\vert
^{2}ds\right)  =R_{1}(1+t^{2}),
\end{align*}
for any $t\geq t_{0},\ 0<\delta\leq\delta_{0}$ and some positive constants
$t_{0}$, $\delta_{0}$ and $R_{1}$, we infer
\[
\lim_{t\rightarrow\infty}\ e^{-ct}2\sup_{0<\delta\leq\delta_{0}}\Vert\phi
y_{\delta}(\theta_{-t+\cdot}\omega)\Vert_{L_{V}^{2}}^{2}=0.
\]
Next, we need to analyze the integral term in $\rho_{\delta}^{2}(\omega)$. By
(\ref{gEst}), we deduce there are positive constants $R_{2}$, $R_{3}$ and
$t_{0}$, such that for $t\geq t_{0}$,
\begin{align*}
e^{-ct}\int_{-\infty}^{0}e^{\gamma s}\Theta_{1,\delta}(\theta_{-t+s}\omega)ds
&  \leq R_{2}e^{-ct}\int_{-\infty}^{0}e^{\gamma s}\left(  1+\left\vert
-t+s\right\vert ^{4p^{2}}ds\right) \\
&  \leq R_{3}e^{-ct}\left(  \left\vert t\right\vert ^{2p}\int_{-\infty}%
^{0}e^{\gamma s}ds+\int_{-\infty}^{0}e^{\gamma s}\left\vert s\right\vert
^{4p^{2}}ds\right)  \rightarrow0,
\end{align*}
as $t\rightarrow+\infty$. Thus, (\ref{LimRadius}) holds, namely, $\overline
{B}$ is tempered for some $\delta_{0}>0.$

\begin{corollary}
The set
\[
\overline{\mathcal{A}}(\omega)=\cup_{0<\delta\leq\delta_{0}}\mathcal{A}%
_{\delta}(\omega),
\]
is tempered for some $\delta_{0}>0$.
\end{corollary}

\section{Upper-semicontinuity}

\label{s6}

In this section, we consider the limiting behavior of the random attractor
$\mathcal{A}_{\delta}$ of the stochastic nonlocal PDEs with long time memory
driven by colored noise \eqref{eq5-1} when $\delta\rightarrow0$.

First, observe that the nonlinear term $f$ satisfies%
\begin{equation}
|f^{\prime}(u)|\leq\widetilde{\beta}(1+|u|^{2p-2}), \label{fcondition2}%
\end{equation}
for some constants $\widetilde{\beta}>0.$ The results of this section are
valid for a general function $f\in C^{1}(\mathbb{R})$ satisfying
(\ref{fcondition}), (\ref{eq630}) and (\ref{fcondition2}).

\begin{lemma}
\label{lem6-1} Assume that $(h_{1})$-$(h_{2})$ hold true, and let {$\phi\in
V\cap H^{2}(\mathcal{O})\cap L^{2p}(\mathcal{O})$} be such that $\Delta\phi\in
L^{2p}(\mathcal{O})$. Then, for every $\omega\in\Omega$, we have
\[
\Vert\eta_{0}-\eta_{0,\delta}\Vert_{L_{\mu}^{2}(\mathbb{R}^{+};V)}^{2}%
\leq2K_{\mu}\left(  \Vert\varphi-\varphi_{\delta}\Vert_{L_{V}^{2}}^{2}%
+\Vert(z_{\ast}(\theta_{\cdot}\omega)-y_{\delta}(\theta_{\cdot}\omega
))\phi\Vert_{L_{V}^{2}}^{2}\right)  ,
\]
where $\varphi$ and $\eta_{0}$ appear in \eqref{eq1-5} and $\varphi_{\delta}$
and $\eta_{0,\delta}$ appear in \eqref{eq5-10}, respectively. $K_{\mu}$ is the
same constant as in Corollary \ref{cor3-2}.
\end{lemma}

\textbf{Proof.} The proof follows the same lines of Corollary \ref{cor3-2}. We
omit the details here. $\Box$

\begin{lemma}
\label{lem6-2} Assume that \eqref{eq1-2}, \eqref{eq2-1} and $(h_{1})$%
-$(h_{2})$ hold. Let $\phi\in V\cap H^{2}(\mathcal{O})\cap L^{2p}%
(\mathcal{O})$ be such that $\Delta\phi\in L^{2p}(\mathcal{O})$, let $a$ be a
locally Lipschitz function. Suppose $u_{\delta}$ and $u$ be solutions to
problems \eqref{eq5-1} and \eqref{eq1-1} with initial data $u_{0,\delta}$ and
$u_{0}$ in $H$, and the initial functions $\varphi_{\delta}$ and $\varphi$ in
$L_{V}^{2}$, respectively. Then, for every $\omega\in\Omega$, $T>0$ and
$\varepsilon\in(0,1)$, there exists $\delta_{0}=\delta_{0}(\omega
,T,\varepsilon)$ and $c=c(\omega,T,\varepsilon,\sup_{t\in\lbrack0,T]}|z_{\ast
}(\theta_{t}\omega)|,\phi)$ such that, for all $0<\delta<\delta_{0}$ and
$t\in\lbrack0,T]$,
%% \red{(Jiaohui, thanks to (6.14) the
%%$\sup_{t\in[0,T]}|y_\delta(\theta_t\omega)|$ can be bounded by $\sup_{t\in[0,T]}|\z|$ and $\varepsilon$, so the constant $c$ depends only on $\sup_{t\in[0,T]}|\z|$.)}%
\begin{equation}%
\begin{split}
&  ~\quad\Vert(u_{\delta}(t;\omega,(u_{0,\delta},\varphi_{\delta}%
)),u_{\delta,t}(\cdot;\omega,(u_{0,\delta},\varphi_{\delta})))-(u(t;\omega
,(u_{0},\varphi)),u_{t}(\cdot;\omega,(u_{0},\varphi)))\Vert_{X}^{2}\\[1ex]
&  \leq c\Vert(u_{0,\delta},\varphi_{\delta})-(u_{0},\varphi)\Vert_{X}%
^{2}\\[0.8ex]
&  +c\varepsilon\left(  1+|u_{0}|^{2}+|u_{0,\delta}|^{2}+\Vert\varphi
\Vert_{L_{V}^{2}}^{2}+\Vert\varphi_{\delta}\Vert_{L_{V}^{2}}^{2}+\int_{0}%
^{t}\Theta_{1}(\theta_{r}\omega)dr+\int_{0}^{t}\Theta_{1,\delta}(\theta
_{r}\omega)dr\right)  ,
\end{split}
\label{eq6-1}%
\end{equation}
where $\Theta_{1}(\omega)$ and $\Theta_{1,\delta}(\omega)$ are the same
constants as in \eqref{theta1} and \eqref{theta2}, respectively.
\end{lemma}

\textbf{Proof.} Let $\xi_{\delta}=v-v_{\delta}$, $\theta_{\delta}^{t}=\eta
^{t}-\eta_{\delta}^{t}$ and $q_{\delta}=(\xi_{\delta},\theta_{\delta}^{t})$
with $q_{0,\delta}=(v_{0}-v_{0,\delta},\eta_{0}-\eta_{0,\delta})$. By
\eqref{eq1-5} and \eqref{eq5-10}, we obtain
\begin{equation}%
\begin{split}
\frac{1}{2}\frac{d}{dt}|\xi_{\delta}|^{2}  &  +\frac{1}{2}\frac{d}{dt}%
\Vert\theta_{\delta}^{t}\Vert_{\mu}^{2}+\left(  -a(l(v)+l(\phi)z_{\ast}%
(\theta_{t}\omega))\Delta v+a(l(v_{\delta})+l(\phi)y_{\delta}(\theta_{t}%
\omega))\Delta v_{\delta},\xi_{\delta}\right) \\[0.8ex]
&  ~~+\left(  -a(l(v)+l(\phi)z_{\ast}(\theta_{t}\omega))z_{\ast}(\theta
_{t}\omega)\Delta\phi+a(l(v_{\delta})+l(\phi)y_{\delta}(\theta_{t}%
\omega))y_{\delta}(\theta_{t}\omega)\Delta\phi,\xi_{\delta}\right) \\[0.8ex]
&  ~~+(f(v+\phi z_{\ast}(\theta_{t}\omega))-f(v_{\delta}+\phi y_{\delta
}(\theta_{t}\omega)),\xi_{\delta})\\[0.8ex]
&  =(z_{\ast}(\theta_{t}\omega)-y_{\delta}(\theta_{t}\omega))(\phi,\xi
_{\delta})+(z_{k}^{\phi}(t,\omega)-z_{k,\delta}^{\phi}(t,\omega),\xi_{\delta
})-(((\theta_{\delta}^{t})^{\prime},\theta_{\delta}^{t}))_{\mu}.
\end{split}
\label{eq6-2}%
\end{equation}
For the third term of left hand side of the above equality, we have
\begin{equation}%
\begin{split}
&  \quad\left(  -a(l(v)+l(\phi)z_{\ast}(\theta_{t}\omega))\Delta
v+a(l(v_{\delta})+l(\phi)y_{\delta}(\theta_{t}\omega))\Delta v_{\delta}%
,\xi_{\delta}\right) \\[0.6ex]
&  =\left(  -a(l(v)+l(\phi)z_{\ast}(\theta_{t}\omega))\Delta v+a(l(v)+l(\phi
)z_{\ast}(\theta_{t}\omega))\Delta v_{\delta},\xi_{\delta}\right) \\[0.6ex]
&  ~~+\left(  -a(l(v)+l(\phi)z_{\ast}(\theta_{t}\omega))\Delta v_{\delta
}+a(l(v_{\delta})+l(\phi)y_{\delta}(\theta_{t}\omega))\Delta v_{\delta}%
,\xi_{\delta}\right) \\[0.6ex]
&  \geq a(l(v)+l(\phi)z_{\ast}(\theta_{t}\omega))\Vert\xi_{\delta}\Vert
^{2}-\left\vert \left(  -a(l(v)+l(\phi)z_{\ast}(\theta_{t}\omega
))+a(l(v_{\delta})+l(\phi)y_{\delta}(\theta_{t}\omega))\right)  \right\vert
\Vert v_{\delta}\Vert\Vert\xi_{\delta}\Vert.
\end{split}
\label{eq6-3}%
\end{equation}
For the fourth term of left hand side of equation \eqref{eq6-2}, by the
Lipschitz condition of the function $a$ and \eqref{eq1-2}, we deduce
\begin{equation}%
\begin{split}
&  ~\quad\left(  -a(l(v)+l(\phi)z_{\ast}(\theta_{t}\omega))z_{\ast}(\theta
_{t}\omega)\Delta\phi+a(l(v_{\delta})+l(\phi)y_{\delta}(\theta_{t}%
\omega))y_{\delta}(\theta_{t}\omega)\Delta\phi,\xi_{\delta}\right) \\[0.6ex]
&  =\left(  -a(l(v)+l(\phi)z_{\ast}(\theta_{t}\omega))z_{\ast}(\theta
_{t}\omega)\Delta\phi+a(l(v)+l(\phi)z_{\ast}(\theta_{t}\omega))y_{\delta
}(\theta_{t}\omega)\Delta\phi,\xi_{\delta}\right) \\[0.6ex]
&  ~~+\left(  -a(l(v)+l(\phi)z_{\ast}(\theta_{t}\omega))y_{\delta}(\theta
_{t}\omega)\Delta\phi+a(l(v_{\delta})+l(\phi)z_{\ast}(\theta_{t}%
\omega))y_{\delta}(\theta_{t}\omega)\Delta\phi,\xi_{\delta}\right) \\[0.6ex]
&  ~~+\left(  -a(l(v_{\delta})+l(\phi)z_{\ast}(\theta_{t}\omega))y_{\delta
}(\theta_{t}\omega)\Delta\phi+a(l(v_{\delta})+l(\phi)y_{\delta}(\theta
_{t}\omega))y_{\delta}(\theta_{t}\omega)\Delta\phi,\xi_{\delta}\right)
\\[0.6ex]
&  =a(l(v)+l(\phi)z_{\ast}(\theta_{t}\omega))(z_{\ast}(\theta_{t}%
\omega)-y_{\delta}(\theta_{t}\omega))(\nabla\phi,\nabla\xi_{\delta})\\[0.6ex]
&  ~~+\left(  -a(l(v)+l(\phi)z_{\ast}(\theta_{t}\omega))y_{\delta}(\theta
_{t}\omega)\Delta\phi+a(l(v_{\delta})+l(\phi)z_{\ast}(\theta_{t}%
\omega))y_{\delta}(\theta_{t}\omega)\Delta\phi,\xi_{\delta}\right) \\[0.6ex]
&  ~~+\left(  -a(l(v_{\delta})+l(\phi)z_{\ast}(\theta_{t}\omega))y_{\delta
}(\theta_{t}\omega)\Delta\phi+a(l(v_{\delta})+l(\phi)y_{\delta}(\theta
_{t}\omega))y_{\delta}(\theta_{t}\omega)\Delta\phi,\xi_{\delta}\right)
\\[0.6ex]
&  \leq M|z_{\ast}(\theta_{t}\omega)-y_{\delta}(\theta_{t}\omega)|\Vert
\phi\Vert\Vert\xi_{\delta}\Vert+L_{a}(R)|l||y_{\delta}(\theta_{t}\omega
)|{|\xi_{\delta}|}\Vert\phi\Vert\Vert{\xi_{\delta}}\Vert\\[0.6ex]
&  ~~+L_{a}(R)|l||\phi||z_{\ast}(\theta_{t}\omega)-y_{\delta}(\theta_{t}%
\omega)||y_{\delta}(\theta_{t}\omega)|\Vert\phi\Vert\Vert\xi_{\delta}\Vert,
\end{split}
\label{6-4}%
\end{equation}
for some $R>0$, which is chosen in a similar way as in the proof of Lemma
\ref{lem4-5} as $v_{\delta}$ are bounded in $C([0,T],H)$ (see Lemma
\ref{lem5-5}). For the last term of left hand side of equation \eqref{eq6-2},
we have
\begin{equation}%
\begin{split}
&  \quad~~\int_{\mathcal{O}}\left(  f(v+\phi z_{\ast}(\theta_{t}%
\omega))-f(v_{\delta}+\phi y_{\delta}(\theta_{t}\omega))\right)  \xi_{\delta
}dx\\
&  =\int_{\mathcal{O}}\left(  f(v+\phi z_{\ast}(\theta_{t}\omega))-f(v+\phi
y_{\delta}(\theta_{t}\omega))\right)  \xi_{\delta}dx\\
&  ~~+\int_{\mathcal{O}}\left(  f(v+\phi y_{\delta}(\theta_{t}\omega
))-f(v_{\delta}+\phi y_{\delta}(\theta_{t}\omega))\right)  \xi_{\delta}dx\\
&  =\int_{\mathcal{O}}f^{\prime}(\vartheta_{1}(x))\xi_{\delta}\phi(z_{\ast
}(\theta_{t}\omega)-y_{\delta}(\theta_{t}\omega))dx+\int_{\mathcal{O}%
}f^{\prime}(\vartheta_{2}(x))|\xi_{\delta}|^{2}dx,
\end{split}
\label{eq6-5}%
\end{equation}
where $\vartheta_{1}(x)=v(x)+\theta_{1}(x)\phi(x)z_{\ast}(\theta_{t}%
\omega)+(1-\theta_{1}(x))\phi(x)y_{\delta}(\theta_{t}\omega)$, $\vartheta
_{2}(x)=\phi(x)y_{\delta}(\theta_{t}\omega)+\theta_{2}(x)v_{\delta
}(x)+(1-\theta_{2}(x))v(x)$, $\theta_{i}(x)\in\lbrack0,1]$. Then, recalling
that $f^{\prime}(u)\geq-\frac{\sigma}{2}$ for some $\sigma>0$ (cf.
\eqref{eq630}), by means of \eqref{fcondition2} and the Young inequality, we
obtain that there exists constants $\widetilde{c}$, $c>0$, such that%
\begin{align}
&  \int_{\mathcal{O}}\left(  f(v+\phi z_{\ast}(\theta_{t}\omega))-f(v_{\delta
}+\phi y_{\delta}(\theta_{t}\omega))\right)  \xi_{\delta}dx\nonumber\\
&  \geq-\widetilde{c}\int_{\mathcal{O}}\left(  1+\left\vert v\right\vert
^{2p-2}+|\phi z_{\ast}(\theta_{t}\omega)|^{2p-2}+|\phi y_{\delta}(\theta
_{t}\omega)|^{2p-2}\right)  \left\vert \xi_{\delta}\right\vert \left\vert
\phi\right\vert |z_{\ast}(\theta_{t}\omega)-y_{\delta}(\theta_{t}%
\omega)|dx-\frac{\sigma}{2}\left\vert \xi_{\delta}\right\vert ^{2}%
\nonumber\\[0.8ex]
&  \geq-c|\xi_{\delta}||\phi||z_{\ast}(\theta_{t}\omega)-y_{\delta}(\theta
_{t}\omega)|-c\left(  |z_{\ast}(\theta_{t}\omega)|^{2p-2}+|y_{\delta}%
(\theta_{t}\omega)|^{2p-2}+1\right) \nonumber\\[0.8ex]
&  \times\left(  \Vert\phi\Vert_{2p}^{2p}+\Vert\xi_{\delta}\Vert_{2p}%
^{2p}+\left\Vert v\right\Vert _{2p}^{2p}\right)  |z_{\ast}(\theta_{t}%
\omega)-y_{\delta}(\theta_{t}\omega)|-\frac{\sigma}{2}\left\vert \xi_{\delta
}\right\vert ^{2}. \label{eq6-5b}%
\end{align}
For the first term of right hand side of \eqref{eq6-2}, by the Young
inequality, we infer
\begin{equation}
(z_{\ast}(\theta_{t}\omega)-y_{\delta}(\theta_{t}\omega))(\phi,\xi_{\delta
})\leq c|z_{\ast}(\theta_{t}\omega)-y_{\delta}(\theta_{t}\omega)|^{2}%
|\phi|^{2}+|\xi_{\delta}|^{2}. \label{eq6-6}%
\end{equation}
For the second term of right hand side of \eqref{eq6-2}, we deduce that
\begin{equation}
(z_{k}^{\phi}(t,\omega)-z_{k,\delta}^{\phi}(t,\omega),\xi_{\delta}%
)=\int_{-\infty}^{t}k(t-s)(z_{\ast}(\theta_{s}\omega)-y_{\delta}(\theta
_{s}\omega))ds(\Delta\phi,\xi_{\delta}). \label{eq6-7}%
\end{equation}
It follows from \eqref{zproperty} that for any $\varepsilon>0$, there exists
$T_{1}<0$ such that for all $t\leq T_{1}$,
\[
|z_{\ast}(\theta_{t}\omega)|\leq\varepsilon|t|.
\]
Similarly, \eqref{eq5-6} implies that for any $\varepsilon>0$, there exists
$T_{2}<0$ such that for all $t\leq T_{2}$, $0<\delta\leq\widetilde{\sigma},$%
\[
|y_{\delta}(\theta_{t}\omega)|\leq\varepsilon|t|.
\]
Notice that
\begin{equation}%
\begin{split}
&  ~~\quad\int_{-\infty}^{t}k(t-s)|z_{\ast}(\theta_{s}\omega)-y_{\delta
}(\theta_{s}\omega)|ds\\
&  =\int_{-\infty}^{T}k(t-s)|z_{\ast}(\theta_{s}\omega)-y_{\delta}(\theta
_{s}\omega)|ds+\int_{T}^{t}k(t-s)|z_{\ast}(\theta_{s}\omega)-y_{\delta}%
(\theta_{s}\omega)|ds.
\end{split}
\label{eq6-8}%
\end{equation}
On the one hand, let $T=\min\{T_{1},T_{2}\}$. We can assume that $t-T\geq1$.
By Remark \ref{rem2-1}, we arrive at
\begin{equation}%
\begin{split}
\int_{-\infty}^{T}k(t-s)|z_{\ast}(\theta_{s}\omega)-y_{\delta}(\theta
_{s}\omega)|ds  &  \leq\int_{-\infty}^{T}k(t-s)|z_{\ast}(\theta_{s}%
\omega)|ds+\int_{-\infty}^{T}k(t-s)|y_{\delta}(\theta_{s}\omega)|ds\\
&  \leq2\varepsilon\int_{-\infty}^{T}k(t-s)|s|ds=2\varepsilon\int%
_{t-T}^{\infty}k(s)|t-s|ds\\
&  \leq2\varepsilon\int_{1}^{\infty}\frac{\mu(1)e^{-\varpi(s-1)}}{\varpi
}|t-s|ds\leq c\varepsilon.
\end{split}
\label{eq6-9}%
\end{equation}
On the other hand, by the continuity of $z_{\ast}(\theta_{t}\omega)$ and
$y_{\delta}(\theta_{t}\omega)$ with respect to $t$, together with Remark
\ref{rem2-1}, we obtain%
\begin{equation}
\int_{T}^{t}k(t-s)|z_{\ast}(\theta_{s}\omega)-y_{\delta}(\theta_{s}%
\omega)|ds<\infty. \label{eq6-10}%
\end{equation}
Collecting \eqref{eq6-8}-\eqref{eq6-10}, it is obvious that \eqref{eq6-7} can
be bounded by
\begin{equation}
(z_{k}^{\phi}(t,\omega)-z_{k,\delta}^{\phi}(t,\omega),\xi_{\delta})\leq
c\varepsilon\Vert\phi\Vert\Vert\xi_{\delta}\Vert+\int_{T}^{t}k(t-s)|z_{\ast
}(\theta_{s}\omega)-y_{\delta}(\theta_{s}\omega)|ds\Vert\phi\Vert\Vert
\xi_{\delta}\Vert. \label{eq6-11}%
\end{equation}
Finally, for the last term of \eqref{eq6-2}, similar to \eqref{eq44}, we find
\begin{equation}
-(((\theta_{\delta}^{t})^{\prime},\theta_{\delta}^{t}))_{\mu}=\int_{0}%
^{\infty}\mu^{\prime}(s)|\nabla\theta_{\delta}^{t}(s)|^{2}ds\leq-\frac{\varpi
}{2}\Vert\theta_{\delta}^{t}\Vert_{\mu}^{2}. \label{eq6-12}%
\end{equation}
Substituting \eqref{eq6-3}-\eqref{eq6-6} and \eqref{eq6-11}-\eqref{eq6-12}
into \eqref{eq6-2}, by \eqref{eq1-2}, we have%
\[%
\begin{split}
&  ~~\quad\frac{1}{2}\frac{d}{dt}\Vert q_{\delta}\Vert_{\mathcal{H}}%
^{2}+m\Vert\xi_{\delta}\Vert^{2}+\frac{\varpi}{2}\Vert\theta_{\delta}^{t}%
\Vert_{\mu}^{2}\\[0.6ex]
&  \leq|-a(l(v)+l(\phi)z_{\ast}(\theta_{t}\omega))+a(l(v_{\delta}%
)+l(\phi)y_{\delta}(\theta_{t}\omega))|\Vert v_{\delta}\Vert\Vert\xi_{\delta
}\Vert+M|z_{\ast}(\theta_{t}\omega)-y_{\delta}(\theta_{t}\omega)|\Vert
\phi\Vert\Vert\xi_{\delta}\Vert\\[0.8ex]
&  ~~+\frac{\sigma}{2}|\xi_{\delta}|^{2}+L_{a}(R)|l||y_{\delta}(\theta
_{t}\omega)||\xi_{\delta}|\Vert\phi\Vert\Vert\xi_{\delta}\Vert+L_{a}%
(R)|l||\phi||z_{\ast}(\theta_{t}\omega)-y_{\delta}(\theta_{t}\omega
)||y_{\delta}(\theta_{t}\omega)|\Vert\phi\Vert\Vert\xi_{\delta}\Vert\\[0.8ex]
&  ~~+c|\xi_{\delta}||\phi||z_{\ast}(\theta_{t}\omega)-y_{\delta}(\theta
_{t}\omega)|+|\xi_{\delta}|^{2}+c\varepsilon\Vert\phi\Vert\Vert\xi_{\delta
}\Vert\\[0.6ex]
&  ~~+c\left(  |z_{\ast}(\theta_{t}\omega)|^{2p-2}+|y_{\delta}(\theta
_{t}\omega)|^{2p-2}+1\right)  \left(  \Vert\phi\Vert_{2p}^{2p}+\Vert
\xi_{\delta}\Vert_{2p}^{2p}+\Vert v\Vert_{2p}^{2p}\right)  |z_{\ast}%
(\theta_{t}\omega)-y_{\delta}(\theta_{t}\omega)|\\
&  ~~+c|z_{\ast}(\theta_{t}\omega)-y_{\delta}(\theta_{t}\omega)|^{2}|\phi
|^{2}+\int_{T}^{t}k(t-s)|z_{\ast}(\theta_{s}\omega)-y_{\delta}(\theta
_{s}\omega)|ds\Vert\phi\Vert\Vert\xi_{\delta}\Vert.
\end{split}
\]
By the Young inequality and the fact that $a$ is locally Lipschitz, we derive
\[%
\begin{split}
&  ~~\quad\frac{d}{dt}\Vert q_{\delta}\Vert_{\mathcal{H}}^{2}+2m\Vert
\xi_{\delta}\Vert^{2}+\varpi\Vert\theta_{\delta}^{t}\Vert_{\mu}^{2}\\[0.6ex]
&  \leq c|\xi_{\delta}|^{2}\Vert v_{\delta}\Vert^{2}+\frac{m}{4}\Vert
\xi_{\delta}\Vert^{2}+c|\phi|^{2}|z_{\ast}(\theta_{t}\omega)-y_{\delta}%
(\theta_{t}\omega)|^{2}\Vert v_{\delta}\Vert^{2}+\frac{m}{4}\Vert\xi_{\delta
}\Vert^{2}+c\Vert\phi\Vert^{2}|z_{\ast}(\theta_{t}\omega)-y_{\delta}%
(\theta_{t}\omega)|^{2}\\[0.6ex]
&  ~~+\frac{m}{4}\Vert\xi_{\delta}\Vert^{2}+c|\xi_{\delta}|^{2}|y_{\delta
}(\theta_{t}\omega)|^{2}\Vert\phi\Vert^{2}+\frac{m}{4}\Vert\xi_{\delta}%
\Vert^{2}+c|\phi|^{2}|y_{\delta}(\theta_{t}\omega)|^{2}\Vert\phi\Vert
^{2}|z_{\ast}(\theta_{t}\omega)-y_{\delta}(\theta_{t}\omega)|^{2}+\frac{m}%
{4}\Vert\xi_{\delta}\Vert^{2}\\[0.6ex]
&  ~~+\sigma|\xi_{\delta}|^{2}+c|\phi|^{2}|z_{\ast}(\theta_{t}\omega
)-y_{\delta}(\theta_{t}\omega)|^{2}+3|\xi_{\delta}|^{2}\\[0.6ex]
&  ~~+c\big(|z_{\ast}(\theta_{t}\omega)|^{2p-2}+|y_{\delta}(\theta_{t}%
\omega)|^{2p-2}+1\big)\left(  \Vert\phi\Vert_{2p}^{2p}+\Vert\xi_{\delta}%
\Vert_{2p}^{2p}+\left\Vert v\right\Vert _{2p}^{2p}\right)  |z_{\ast}%
(\theta_{t}\omega)-y_{\delta}(\theta_{t}\omega)|+\frac{m}{4}\Vert\xi_{\delta
}\Vert^{2}\\[0.6ex]
&  ~~+c\left(  \int_{T}^{t}k(t-s)|z_{\ast}(\theta_{s}\omega)-y_{\delta}%
(\theta_{s}\omega)|ds\right)  ^{2}\Vert\phi\Vert^{2}+\frac{m}{4}\Vert
\xi_{\delta}\Vert^{2}+c\varepsilon^{2}\left\Vert \phi\right\Vert ^{2}.
\end{split}
\]
Therefore,
\begin{equation}%
\begin{split}
&  \frac{d}{dt}\Vert q_{\delta}\Vert_{\mathcal{H}}^{2}+\frac{m}{4}\left\Vert
\xi_{\delta}\right\Vert ^{2}+\varpi\Vert\theta_{\delta}^{t}\Vert_{\mu}%
^{2}\\[0.6ex]
&  \leq\left(  c\Vert v_{\delta}\Vert^{2}+\sigma+3+c|y_{\delta}(\theta
_{t}\omega)|^{2}\Vert\phi\Vert^{2}\right)  |\xi_{\delta}|^{2}\\[0.6ex]
&  ~~+c\left(  |\phi|^{2}\Vert v_{\delta}\Vert^{2}+\Vert\phi\Vert^{2}%
+|\phi|^{2}|y_{\delta}(\theta_{t}\omega)|^{2}\Vert\phi\Vert^{2}+|\phi
|^{2}\right)  |z_{\ast}(\theta_{t}\omega)-y_{\delta}(\theta_{t}\omega
)|^{2}\\[0.6ex]
&  ~~+c\left(  |z_{\ast}(\theta_{t}\omega)|^{2p-2}+|y_{\delta}(\theta
_{t}\omega)|^{2p-2}+1\right)  \left(  \Vert\phi\Vert_{2p}^{2p}+\Vert
\xi_{\delta}\Vert_{2p}^{2p}+\left\Vert v\right\Vert _{2p}^{2p}\right)
|z_{\ast}(\theta_{t}\omega)-y_{\delta}(\theta_{t}\omega)|\\[0.6ex]
&  ~~+c\varepsilon^{2}\Vert\phi\Vert^{2}+c\left(  \int_{T}^{t}k(t-s)|z_{\ast
}(\theta_{s}\omega)-y_{\delta}(\theta_{s}\omega)|ds\right)  ^{2}\Vert\phi
\Vert^{2}.
\end{split}
\label{eq6-132}%
\end{equation}
With the help of \eqref{eq5-5}, for every $\varepsilon>0$, there exists
$\delta_{0}=\delta_{0}(\omega,T,\varepsilon)>0$ such that for all
$0<\delta<\delta_{0}$ and $t\in\lbrack0,T]$,
\begin{equation}
|z_{\ast}(\theta_{t}\omega)-y_{\delta}(\theta_{t}\omega)|\leq\varepsilon.
\label{eq613}%
\end{equation}
Notice that, thanks to this fact, there exists a constant $c:=c(\omega
,T,\varepsilon)$ such that,
\begin{equation}
\sup_{t\in\lbrack0,T]}|y_{\delta}(\theta_{t}\omega)|\leq c\sup_{t\in
\lbrack0,T]}|z_{\ast}(\theta_{t}\omega)|,\quad\mbox{for all}\ 0<\delta
<\delta_{0}. \label{eq613b}%
\end{equation}
Hence, it follows from Remark \ref{rem2-1} that
\begin{equation}%
\begin{split}
\frac{d}{dt}\Vert q_{\delta}\Vert_{\mathcal{H}}^{2}+\frac{m}{4}\left\Vert
\xi_{\delta}\right\Vert ^{2}  &  \leq c\left(  \Vert v_{\delta}\Vert
^{2}+1+|y_{\delta}(\theta_{t}\omega)|^{2}\Vert\phi\Vert^{2}\right)  \Vert
q_{\delta}\Vert_{\mathcal{H}}^{2}\\[0.8ex]
&  ~~+c\varepsilon^{2}\left(  |\phi|^{2}\Vert v_{\delta}\Vert^{2}+\Vert
\phi\Vert^{2}+|\phi|^{2}|y_{\delta}(\theta_{t}\omega)|^{2}\Vert\phi\Vert
^{2}+|\phi|^{2}\right) \\[0.6ex]
&  ~~+c\varepsilon\left(  |z_{\ast}(\theta_{t}\omega)|^{2p-2}+|y_{\delta
}(\theta_{t}\omega)|^{2p-2}+1\right)  \left(  \Vert\phi\Vert_{2p}^{2p}%
+\Vert\xi_{\delta}\Vert_{2p}^{2p}+\left\Vert v\right\Vert _{2p}^{2p}\right)  .
\end{split}
\label{eq6-13}%
\end{equation}
Multiplying by $e^{-c\int_{0}^{t}(\Vert v_{\delta}\Vert^{2}+1+|y_{\delta
}(\theta_{\tau}\omega)|^{2}\Vert\phi\Vert^{2})d\tau}$ and integrating in
\eqref{eq6-13}, we deduce that, for all $0<\delta<\delta_{0}$ and $t\in
\lbrack0,T]$,
\begin{equation}%
\begin{split}
&  \Vert q_{\delta}\Vert_{\mathcal{H}}^{2}+\frac{m}{4}\int_{0}^{t}e^{c\int%
_{s}^{t}(\Vert v_{\delta}\Vert^{2}+1+|y_{\delta}(\theta_{\tau}\omega
)|^{2}\Vert\phi\Vert^{2})d\tau}\left\Vert \xi_{\delta}\right\Vert ^{2}ds\\
&  \leq e^{c\int_{0}^{t}(\Vert v_{\delta}\Vert^{2}+1+|y_{\delta}(\theta
_{s}\omega)|^{2}\Vert\phi\Vert^{2})ds}\bigg(\Vert q_{0,\delta}\Vert
_{\mathcal{H}}^{2}+\int_{0}^{t}e^{-c\int_{0}^{s}(\Vert v_{\delta}\Vert
^{2}+1+|y_{\delta}(\theta_{\tau}\omega)|^{2}\Vert\phi\Vert^{2})d\tau}\\
&  ~~\times\bigg(c\varepsilon^{2}\left(  |\phi|^{2}\Vert v_{\delta}\Vert
^{2}+\Vert\phi\Vert^{2}+|\phi|^{2}|y_{\delta}(\theta_{s}\omega)|^{2}\Vert
\phi\Vert^{2}+|\phi|^{2}\right) \\[0.6ex]
&  ~~+c\varepsilon\left(  |z_{\ast}(\theta_{s}\omega)|^{2p-2}+|y_{\delta
}(\theta_{s}\omega)|^{2p-2}+1\right)  \left(  \Vert\phi\Vert_{2p}^{2p}%
+\Vert\xi_{\delta}\Vert_{2p}^{2p}+\left\Vert v\right\Vert _{2p}^{2p}\right)
\bigg)ds\bigg).
\end{split}
\label{eq6-15}%
\end{equation}
In view of \eqref{eq5-13}, $\phi\in V$ and (\ref{eq613b}), there is
$c=c(T,\omega,\varepsilon,\Vert\phi\Vert)$ such that,
\[
\int_{0}^{T}\left(  \Vert v_{\delta}\Vert^{2}+1+|y_{\delta}(\theta_{s}%
\omega)|^{2}\Vert\phi\Vert^{2}\right)  ds\leq c,~~\text{ if }0<\delta
<\delta_{0}.
\]
By \eqref{eq6-15}, \eqref{eq4-12} and \eqref{eq5-13}, there exist $\delta
_{1}\in(0,\delta_{0})$ and $c:=c(\omega,T,\varepsilon,\phi,\sup_{t\in
\lbrack0,T]}|z_{\ast}(\theta_{t}\omega)|)$ such that, for all $0<\delta
<\delta_{1}$ and $t\in\lbrack0,T]$,
\[%
\begin{split}
\Vert q_{\delta}\Vert_{\mathcal{H}}^{2}  &  +\frac{m}{4}\int_{0}^{t}\left\Vert
\xi_{\delta}\right\Vert ^{2}ds\\
&  \leq c\left(  |v_{0}-v_{0,\delta}|^{2}+\Vert\eta_{0}-\eta_{0,\delta}%
\Vert_{\mu}^{2}\right) \\
&  +c\varepsilon\left(  1+|v_{0}|^{2}+|v_{0,\delta}|^{2}+\Vert\eta_{0}%
\Vert_{\mu}^{2}+\Vert\eta_{0,\delta}\Vert_{\mu}^{2}+\int_{0}^{t}\Theta
_{1}(\theta_{r}\omega)dr+\int_{0}^{t}\Theta_{1,\delta}(\theta_{r}%
\omega)dr\right)  .
\end{split}
\]
Notice that
\begin{equation}
u_{\delta}(t;\omega,u_{0,\delta})-u(t;\omega,u_{0})=v_{\delta}(t;\omega
,v_{0,\delta})-v(t;\omega,v_{0})+\phi y_{\delta}(\theta_{t}\omega)-\phi
z_{\ast}(\theta_{t}\omega), \label{eq618}%
\end{equation}
where $u_{0,\delta}=v_{0,\delta}+\phi y_{\delta}(\omega)$ and $u_{0}%
=v_{0}+\phi z_{\ast}(\omega)$. It follows from the above equations,
corollaries \ref{cor3-2} and \ref{cor3-2b}, Lemma \ref{lem6-1} and
\eqref{eq613}, that there exist $\delta_{2}\in(0,\delta_{1})$ and
$c:=c(\omega,T,\varepsilon,\phi,\sup_{t\in\lbrack0,T]}|z_{\ast}(\theta
_{t}\omega)|)$, such that for all $0<\delta<\delta_{2}$ and $t\in\lbrack0,T]$,%
\begin{equation}%
\begin{split}
&  |u_{\delta}-u|^{2}+\int_{0}^{t}\left\Vert \xi_{\delta}\right\Vert
^{2}ds\\[0.8ex]
&  \leq c\left(  |u_{0}-u_{0,\delta}|^{2}+\Vert\varphi-\varphi_{\delta}%
\Vert_{L_{V}^{2}}^{2}\right) \\
&  +c\varepsilon\left(  1+|u_{0}|^{2}+|u_{0,\delta}|^{2}+\Vert\varphi
\Vert_{L_{V}^{2}}^{2}+\Vert\varphi_{\delta}\Vert_{L_{V}^{2}}^{2}+\int_{0}%
^{t}\Theta_{1}(\theta_{r}\omega)dr+\int_{0}^{t}\Theta_{1,\delta}(\theta
_{r}\omega)dr\right)  .
\end{split}
\label{eq6-19}%
\end{equation}
By (\ref{eq613}) and \eqref{eq618}-\eqref{eq6-19}, we obtain
\begin{equation}%
\begin{split}
&  ~~\quad\int_{0}^{t}\Vert u_{\delta}(s)-u(s)\Vert^{2}ds\\[0.8ex]
&  \leq2\int_{0}^{t}\Vert\xi_{\delta}(s)\Vert^{2}ds+2\int_{0}^{t}\Vert
\phi\Vert^{2}|y_{\delta}(\theta_{s}\omega)-z_{\ast}(\theta_{s}\omega
)|^{2}ds\\[0.8ex]
&  \leq c\left(  |u_{0}-u_{0,\delta}|^{2}+\Vert\varphi-\varphi_{\delta}%
\Vert_{L_{V}^{2}}^{2}\right) \\
&  +c\varepsilon\left(  1+|u_{0}|^{2}+|u_{0,\delta}|^{2}+\Vert\varphi
\Vert_{L_{V}^{2}}^{2}+\Vert\varphi_{\delta}\Vert_{L_{V}^{2}}^{2}+\int_{0}%
^{t}\Theta_{1}(\theta_{r}\omega)dr+\int_{0}^{t}\Theta_{1,\delta}(\theta
_{r}\omega)dr\right)  .
\end{split}
\end{equation}
Hence, for every $\omega\in\Omega$ and $t\in\lbrack0,T]$, we have
\begin{equation}%
\begin{split}
&  \Vert u_{\delta,t}-u_{t}\Vert_{L_{V}^{2}}^{2}=\int_{-\infty}^{0}e^{\gamma
s}\Vert u_{\delta}(t+s)-u(t+s)\Vert^{2}ds\\[0.8ex]
&  =\int_{-\infty}^{0}e^{-\gamma(t-s)}\Vert\varphi_{\delta}(s)-\varphi
(s)\Vert^{2}ds+\int_{0}^{t}e^{-\gamma(t-s)}\Vert u_{\delta}(s)-u(s)\Vert
^{2}ds\\[0.8ex]
&  \leq\Vert\varphi_{\delta}-\varphi\Vert_{L_{V}^{2}}^{2}+\int_{0}^{t}\Vert
u_{\delta}(s)-u(s)\Vert^{2}ds\\[0.8ex]
&  \leq c\left(  |u_{0}-u_{0,\delta}|^{2}+\Vert\varphi-\varphi_{\delta}%
\Vert_{L_{V}^{2}}^{2}\right) \\
&  +c\varepsilon\left(  1+|u_{0}|^{2}+|u_{0,\delta}|^{2}+\Vert\varphi
\Vert_{L_{V}^{2}}^{2}+\Vert\varphi_{\delta}\Vert_{L_{V}^{2}}^{2}+\int_{0}%
^{t}\Theta_{1}(\theta_{r}\omega)dr+\int_{0}^{t}\Theta_{1,\delta}(\theta
_{r}\omega)dr\right)  ,
\end{split}
\end{equation}
which, together with \eqref{eq6-19}, finishes the proof of this lemma. $\Box$

\begin{remark}
The constant $c$ depends continuously on $\varepsilon$ in Lemma \ref{lem6-2}.
\end{remark}

As a consequence of Lemma \ref{lem6-2}, we obtain the following covergence of
solutions to \eqref{eq5-1} when $\delta$ approaches to zero.

\begin{corollary}
\label{cor6-3} Assume that \eqref{eq1-2}, \eqref{eq2-1}, $(h_{1})$-$(h_{2})$
hold and $\delta_{n}\rightarrow0$ as $n\rightarrow\infty$. Let $\phi\in V\cap
H^{2}(\mathcal{O})\cap L^{2p}(\mathcal{O})$ be such that $\Delta\phi\in
L^{2p}(\mathcal{O})$ and $a$ be a locally Lipschitz function. Suppose that
$u_{\delta_{n}}$ and $u$ are the solutions of \eqref{eq5-1} and \eqref{eq1-5}
with initial data $u_{0,\delta_{n}}$ and $u_{0}$ in $H$, and the initial
functions $\varphi_{\delta_{n}}$ and $\varphi$ in $L_{V}^{2}$, respectively.
If {$u_{0,\delta_{n}}\rightarrow u_{0}$ in $H$ and $\varphi_{\delta_{n}%
}\rightarrow\varphi$ in $L_{V}^{2}$ as $n\rightarrow\infty$}, then for every
$\omega\in\Omega$ and {$t>0$},
\[
u_{\delta_{n}}(t;\omega,(u_{0,\delta_{n}},\varphi_{\delta_{n}}))\rightarrow
u(t;\omega,(u_{0},\varphi))\quad\mbox{in}\quad H,\quad\mbox{as}~~n\rightarrow
\infty,
\]
and
\[
u_{\delta_{n},t}(\cdot;\omega,(u_{0,\delta_{n}},\varphi_{\delta_{n}%
}))\rightarrow u_{t}(\cdot;\omega,(u_{0},\varphi))\quad\mbox{in}\quad
L_{V}^{2},\quad\mbox{as}~~n\rightarrow\infty.
\]
The above convergence is uniform with respect to $t\in\lbrack0,T].$
\end{corollary}

We also need the following weak convergence of solutions to prove the
upper-semicontinuity of random attractors in this section.

\begin{lemma}
\label{lem1} Under assumptions of Corollary \ref{cor6-3}, suppose that
$\{\delta_{n}\}_{n=1}^{\infty}$ is a sequence such that $\delta_{n}%
\rightarrow0$ as $n\rightarrow\infty$. Let $v_{\delta_{n}}$ and $v$ be the
solutions of (\ref{eq5-8}) and (\ref{eq1-4}) with initial data $v_{0,\delta
_{n}}$ and $v_{0}$ in $H$, and the initial functions $\varphi_{v,\delta_{n}%
}=\varphi_{\delta_{n}}-\phi y_{\delta_{n}}(\theta_{\text{\textperiodcentered}%
}\omega)$ and $\varphi_{v}=\varphi-\phi z_{\ast}(\theta
_{\text{\textperiodcentered}}\omega)$, respectively. If $v_{0,\delta_{n}%
}\rightarrow v_{0}$ weakly in $H$ and $\varphi_{v,\delta_{n}}\rightarrow
\varphi_{v}$ weakly in $L_{V}^{2}$ as $n\rightarrow\infty$. Then for every
$\omega\in\Omega$,
\begin{equation}%
\begin{split}
&  ~\quad\left(  v_{\delta_{n}}(r;0,\omega,(v_{0,\delta_{n}},\mathcal{J}%
_{\omega,0}^{\delta_{n}}{\varphi_{\delta_{n}}})),v_{\delta_{n},r}%
(\cdot;0,\omega,(v_{0,\delta_{n}},{\mathcal{J}^{\delta_{n}}_{\omega,0}%
\varphi_{\delta_{n}}}))\right) \\[0.8ex]
&  \rightarrow\left(  v(r;0,\omega,(v_{0},\mathcal{J}_{\omega,0}{\varphi
})),v_{r}(\cdot;0,\omega,(v_{0},\mathcal{J}_{\omega,0}{\varphi}))\right)
~~\mbox{weakly in}~~X,~~\forall r\geq0;
\end{split}
\end{equation}%
\begin{equation}
{v_{\delta_{n}}(\cdot;0,\omega,(v_{0,\delta_{n}},\mathcal{J}_{\omega
,0}^{\delta_{n}}\varphi}_{\delta_{n}}{))\rightarrow v(\cdot;0,\omega
,(v_{0},\mathcal{J}_{\omega,0}\varphi))~~\mbox{weakly in}\quad L^{2p}%
(0,T;L^{2p}(\mathcal{O})),~~\forall T>0;}%
\end{equation}%
\begin{equation}
v_{\delta_{n}}(\cdot;0,\omega,(v_{0,\delta_{n}},\mathcal{J}_{\omega,0}%
^{\delta_{n}}\varphi_{\delta_{n}}))\rightarrow v(\cdot;0,\omega,(v_{0}%
,\mathcal{J}_{\omega,0}\varphi))~~\mbox{weakly in}\quad L^{2}(0,T;V),~~\forall
T>0;
\end{equation}%
\begin{equation}
\eta_{0,\delta_{n}}:=\mathcal{J}_{\omega,0}^{\delta_{n}}\varphi_{\delta_{n}%
}\rightarrow\eta_{0}:=\mathcal{J}_{\omega,0}\varphi~~\mbox{weakly in}~~L_{\mu
}^{2}(\mathbb{R}^{+};V),
\end{equation}
as $n\rightarrow\infty$.
\end{lemma}

\textbf{Proof.} The proof of this lemma follows the same arguments as
\cite[Lemma 3.5]{G4}, so we omit the details here. $\Box$

Recall that for each $\delta>0$, $\mathcal{A}_{\delta}$ is the {unique}
$\mathcal{D}$-random attractor of $\Xi_{\delta}$ in $X$. To establish the
upper semicontinuity of these attractors as $\delta\rightarrow0$, we need the
following compactness result.

\begin{lemma}
\label{lem6-4} Suppose that conditions of Corollary \ref{cor6-3} hold. Let
$\omega\in\Omega$ be fixed. If $\delta_{n}\rightarrow0$ as $n\rightarrow
\infty$ and $(u^{n},\varphi^{n})\in\mathcal{A}_{\delta_{n}}(\omega)$, then the
sequence $\{(u^{n},\varphi^{n})\}_{n=1}^{\infty}$ has a convergent subsequence
in $X$.
\end{lemma}

\textbf{Proof.} Since $\delta_{n}\rightarrow0$, by \eqref{rhoconvergence} we
find that for every $\omega\in\Omega$, there exists $N_{1}=N_{1}(\omega)$ such
that for all $n\geq N_{1}$,
\begin{equation}
\rho_{\delta_{n}}^{2}(\omega)\leq2\rho^{2}(\omega). \label{eq3}%
\end{equation}
Due to $(u^{n},\varphi^{n})\in\mathcal{A}_{\delta_{n}}(\omega)$ and
$\mathcal{A}_{\delta_{n}}(\omega)\subset B_{\delta_{n}}(0,\rho_{\delta_{n}%
}^{2}(\omega))$, by \eqref{eq3} we obtain that, for all $n\geq N_{1}$,
\begin{equation}
\Vert(u^{n},\varphi^{n})\Vert_{X}^{2}\leq2\rho^{2}(\omega). \label{eq4}%
\end{equation}
It follows from \eqref{eq4} that the sequence $\{(u^{n},\varphi^{n}%
)\}_{n=1}^{\infty}$ is bounded in $X$, hence, there exists $(u_{0},\varphi
_{0})\in X$ such that, up to a subsequence,
\begin{equation}
(u^{n},\varphi^{n})\rightarrow(u_{0},\varphi_{0}),\quad
\mbox{weakly in}~X~\mbox{as}~~n\rightarrow\infty. \label{eq5}%
\end{equation}

In what follows, we will prove that the weak convergence in \eqref{eq5} is
actually a strong one in $X$. Since $(u^{n},\varphi^{n})\in\mathcal{A}%
_{\delta_{n}}(\omega)$, by the invariance of $\mathcal{A}_{\delta_{n}}$, for
every $k\geq1$, there exists $(u^{n,k},\varphi^{n,k})\in\mathcal{A}%
_{\delta_{n}}(\theta_{-k}\omega)$ such that,
\begin{equation}%
\begin{split}
(u^{n},\varphi^{n})  &  =\Xi_{\delta_{n}}(k,\theta_{-k}\omega,(u^{n,k}%
,\varphi^{n,k}))\\[1ex]
&  =(u_{\delta_{n}}(0;-k,\omega,(u^{n,k},\varphi^{n,k})),u_{\delta_{n}%
,0}(\cdot;-k,\omega,(u^{n,k},\varphi^{n,k}))).
\end{split}
\label{eq6}%
\end{equation}
On the one hand, since $(u^{n,k},\varphi^{n,k})\in\mathcal{A}_{\delta_{n}%
}(\theta_{-k}\omega)$ and $\mathcal{A}_{\delta_{n}}(\theta_{-k}\omega)\subset
B_{\delta_{n}}(0,\rho_{\delta_{n}}(\theta_{-k}\omega))$, by \eqref{eq3}, we
infer that for each $k\geq1$ and $n\geq N_{1}(\theta_{-k}\omega)$,
\begin{equation}
\Vert(u^{n,k},\varphi^{n,k})\Vert_{X}^{2}\leq2\rho^{2}(\theta_{-k}\omega).
\label{eq7}%
\end{equation}
On the other hand, by \eqref{eq5-2}, {denoting $\varphi_{v}^{n,k}%
:=\varphi_{v,\delta_{n}}^{k}$}, we have
\begin{equation}%
\begin{split}
&  ~\quad\left(  v_{\delta_{n}}(0;-k,\omega,(v^{n,k},\varphi_{v}%
^{n,k})),v_{\delta_{n},0}(\cdot;-k,\omega,(v^{n,k},\varphi_{v}^{n,k}))\right)
\\[0.6ex]
&  =\left(  u_{\delta_{n}}(0;-k,\omega,(u^{n,k},\varphi^{n,k})),u_{\delta
_{n},0}(\cdot;-k,\omega,(u^{n,k},\varphi^{n,k}))\right)  -(\phi y_{\delta_{n}%
}(\omega),\phi y_{\delta_{n}}(\theta_{\cdot}\omega)),
\end{split}
\label{eq8}%
\end{equation}
where%
\begin{equation}
(v^{n,k},\varphi_{v}^{n,k})=(u^{n,k},\varphi^{n,k})-(\phi y_{\delta}{}%
_{n}(\theta_{-k}\omega),\phi y_{\delta}{}_{n}(\theta_{-k+\cdot}\omega)).
\label{eq9}%
\end{equation}
By \eqref{eq6} and \eqref{eq8}, we obtain
\begin{equation}
(u^{n},\varphi^{n})=\left(  v_{\delta_{n}}(0;-k,\omega,(v^{n,k},\varphi
_{v}^{n,k})),v_{\delta_{n},0}(\cdot;-k,\omega,(v^{n,k},\varphi_{v}%
^{n,k}))\right)  +(\phi y_{\delta_{n}}(\omega),\phi y_{\delta_{n}}%
(\theta_{\cdot}\omega)). \label{eq10}%
\end{equation}
By \eqref{eq7} and \eqref{eq9}, we have, for $n\geq N_{1}(\theta_{-k}\omega
)$,
\begin{equation}
\Vert(v^{n,k},\varphi_{v}^{n,k})\Vert_{X}^{2}\leq4\left(  \rho^{2}(\theta
_{-k}\omega)+|\phi|^{2}y_{\delta_{n}}^{2}(\theta_{-k}\omega)+\Vert\phi
y_{\delta_{n}}(\theta_{\cdot}\omega)\Vert_{L_{V}^{2}}^{2}\right)  .
\label{eq11}%
\end{equation}
Now, by \eqref{eq5-5}, Remark \ref{ConvergYdelta} and \eqref{eq11}, we find
that there exists $N_{2}=N_{2}(\omega,k)\geq N_{1}$, such that for every
$k\geq1$ and $n\geq N_{2}$,
\begin{equation}
\Vert(v^{n,k},\varphi_{v}^{n,k})\Vert_{X}^{2}\leq4\rho^{2}(\theta_{-k}%
\omega)+8|\phi|^{2}z_{\ast}^{2}(\theta_{-k}\omega)+8\Vert\phi z_{\ast}%
(\theta_{\cdot}\omega)\Vert_{L_{V}^{2}}^{2}. \label{eq12}%
\end{equation}
Note that \eqref{eq5}, \eqref{eq10}, \eqref{eq5-5} and Remark
\ref{ConvergYdelta} imply that, as $n\rightarrow\infty$,
\begin{equation}
(v_{\delta_{n}}(0;-k,\omega,(v^{n,k},\varphi_{v}^{n,k})),v_{\delta_{n}%
,0}(\cdot;-k,\omega,(v^{n,k},\varphi_{v}^{n,k})))\rightarrow(v_{0}%
,\varphi_{v,0})\quad\mbox{weakly in}~X, \label{eq13}%
\end{equation}
with
\begin{equation}
(v_{0},\varphi_{v,0})=(u_{0},\varphi_{0})-(\phi z_{\ast}(\omega),\phi z_{\ast
}(\theta_{\cdot}\omega)). \label{eq14}%
\end{equation}
By \eqref{eq12}, we find that for each fixed $k\geq1$, the sequence
$\{(v^{n,k},\varphi_{v}^{n,k})\}$ is bounded in $X$, and hence, there is a
subsequence (not relabeled) such that for every $k\geq1$, there exists
$(\tilde{v}^{k},\tilde{\varphi}_{v}^{k})\in X$ such that
\begin{equation}
(v^{n,k},\varphi_{v}^{n,k})\rightarrow(\tilde{v}^{k},\tilde{\varphi}_{v}%
^{k})\quad\mbox{weakly in}~~X~~\mbox{as}~n\rightarrow\infty. \label{eq15}%
\end{equation}
By \eqref{eq15} and Lemma \ref{lem1}, we find
\begin{equation}
v_{\delta_{n}}(0;-k,\omega,(v^{n,k},\varphi_{v}^{n,k}))\rightarrow
v(0;-k,\omega,(\tilde{v}^{k},\tilde{\varphi}_{v}^{k}))\quad
\mbox{weakly in}~H~\mbox{as}~n\rightarrow\infty, \label{eq16}%
\end{equation}
and
\begin{equation}
v_{\delta_{n},0}(\cdot;-k,\omega,(v^{n,k},\varphi_{v}^{n,k}))\rightarrow
v_{0}(\cdot;-k,\omega,(\tilde{v}^{k},\tilde{\varphi}_{v}^{k}))\quad
\mbox{weakly in}~L_{V}^{2}~\mbox{as}~n\rightarrow\infty. \label{eq17}%
\end{equation}
Now, by \eqref{eq13} and \eqref{eq16}-\eqref{eq17}, we have
\begin{equation}
(v_{0},\varphi_{v,0})=(v(0;-k,\omega,(\tilde{v}^{k},\tilde{\varphi}_{v}%
^{k}),v_{0}(\cdot;-k,\omega,(\tilde{v}^{k},\tilde{\varphi}_{v}^{k}))).
\label{eq18}%
\end{equation}
We need to prove that, up to a subsequence, the convergence (\ref{eq5}) is
also true with respect to the strong topology. We will do it in several steps.

\bigskip

\textbf{Statement 1. }We have%
\begin{equation}
v_{\delta_{n}}(0;-k,\omega,(v^{n,k},\varphi_{v}^{n,k}))\rightarrow
v(0;-k,\omega,(\tilde{v}^{k},\tilde{\varphi}_{v}^{k}))\quad\text{strongly
in}~H~\text{as}~n\rightarrow\infty. \label{eq16b}%
\end{equation}
In a similar way as in Lemma \ref{lem4-5}, we obtain that%
\begin{equation}
v_{\delta_{n}}(\text{\textperiodcentered};-k,\omega,(v^{n,k},\varphi_{v}%
^{n,k}))\rightarrow v(\text{\textperiodcentered};-k,\omega,(\tilde{v}%
^{k},\tilde{\varphi}_{v}^{k}))\quad\text{in }L^{2}(-k,0;H). \label{ConvergL2H}%
\end{equation}
Thus,
\[
v_{\delta_{n}}(t;-k,\omega,(v^{n,k},\varphi_{v}^{n,k}))\left(  x\right)
\rightarrow v(t;-k,\omega,(\tilde{v}^{k},\tilde{\varphi}_{v}^{k}))\left(
x\right)  \quad\text{for a.a. }t\in\left(  -k,0\right)  \times\mathcal{O}.
\]
Let us denote by
\[
z_{\delta_{n}}\left(  t\right)  =\left(  v_{\delta_{n}}(t),\eta_{\delta_{n}%
}^{t}\left(  s\right)  \right)  ,
\]
the solution to problem (\ref{eq5-10}) with initial condition $(v^{n,k}%
,\mathcal{J}_{\omega,0}^{\delta_{n}}{\varphi^{n,k}})$ on $t=-k$. Integrating
in (\ref{eq5-13}) over $\left(  -k,t\right)  $\ for $t\in\lbrack-k,0)$ we have%
\begin{equation}
\Vert z_{\delta_{n}}(t)\Vert_{\mathcal{H}}^{2}+\frac{m}{2}\int_{-k}%
^{t}\left\Vert v_{\delta_{n}}\left(  s\right)  \right\Vert ^{2}ds\leq\Vert
z_{\delta n}(-k)\Vert_{\mathcal{H}}^{2}+\int_{-k}^{t}\Theta_{1,\delta_{n}%
}(\theta_{s}\omega)ds\leq R_{1}\left(  \omega,k\right)  , \label{Ineqzdelta}%
\end{equation}
where $\Theta_{1,\delta_{n}}(\omega)$ is defined in (\ref{theta2}). Observe
that the existence of the bound $R_{1}\left(  \omega,k\right)  $ follows from
(\ref{eq15}), (\ref{eq5-5}), (\ref{ConvergC1delta}) and (\ref{ConvergC2delta}).

In view of (\ref{Ineqzdelta}), $v^{n}(t)=v_{\delta_{n}}(t;-k,\omega
,(v^{n,k},\varphi_{v}^{n,k}))$ is bounded in $C([-k,0],H)$, so the same
arguments as in Lemma \ref{lem4-5} show that%
\[
v^{n}\rightarrow v\text{ in }C([-k,0],V^{\ast}+L^{q}(\mathcal{O})),
\]
where $v(t)=v(t;-k,\omega,(\tilde{v}^{k},\tilde{\varphi}_{v}^{k}))$. Then if
$t_{n}\rightarrow t_{0}$, $t_{n}\in\lbrack-k,0]$, $t_{0}\in(-k,0]$, we obtain%
\[
v^{n}(t_{n})\rightarrow v(t_{0})\text{ weakly in }H,
\]
and
\[
\left\vert v(t_{0})\right\vert \leq\liminf_{n\rightarrow\infty}\ \left\vert
v_{n}(t_{n}))\right\vert .
\]
Let us prove that $v^{n}(t_{n})\rightarrow v(t_{0})$ strongly in $H$. Using
Lemma \ref{lemEquivInt}, we deduce that $v^{_{n}}$ are weak solutions to
problem (\ref{eq5-8}). Multiplying the equation by $v^{n}$, we have
\[%
\begin{split}
&  \quad~~\frac{1}{2}\frac{d}{dt}|v^{n}(t)|^{2}+m\Vert v^{n}\Vert^{2}%
+(f(v^{n}+\phi y_{\delta}(\theta_{t}\omega)),v^{n})\\[0.8ex]
&  \leq M|y_{\delta}(\theta_{t}\omega)|\Vert\phi\Vert\Vert v^{n}\Vert+\left(
\int_{-\infty}^{t}k(t-s)\Delta v^{n}(s)ds,v^{n}(t)\right)  +(h,v^{n}%
(t))\\[0.8ex]
&  ~~+|y_{\delta}(\theta_{t}\omega)||\phi||v^{n}|+(z_{k,\delta}^{\phi}%
(\theta_{t}\omega),v^{n}).
\end{split}
\]
By similar arguments as in Lemma \ref{lem5-5} and the Young inequality, we
obtain
\[%
\begin{split}
&  ~~\quad\frac{d}{dt}|v^{n}|^{2}+m\Vert v^{n}\Vert^{2}+\frac{f_{0}}{2D}\Vert
v^{n}\Vert_{2p}^{2p}\leq2\alpha|\mathcal{O}|+C_{1,\delta}(\theta_{t}%
\omega)(1+\Vert\phi\Vert_{2p}^{2p})\\[0.8ex]
&  ~~+\frac{4M^{2}}{m}|y_{\delta}(\theta_{t}\omega)|^{2}\Vert\phi\Vert
^{2}+\frac{4}{m\lambda_{1}}|h|^{2}+\frac{4}{m\lambda_{1}}|y_{\delta}%
(\theta_{t}\omega)|^{2}|\phi|^{2}\\[0.8ex]
&  ~~+\frac{C_{2,\delta}(\theta_{t}\omega)}{m}\left\Vert \phi\right\Vert
^{2}+2\int_{-\infty}^{t}k(t-s)\Vert v^{n}(s)\Vert ds\Vert v^{n}(t)\Vert.
\end{split}
\]
By the same arguments in Lemma \ref{lem4-5}, we have%
\[
\int_{-\infty}^{t}k(t-s)\Vert v^{n}(s)\Vert ds\leq\frac{M_{1}^{\frac{1}{2}}%
\mu^{\frac{1}{2}}(t)\Vert\varphi_{v}^{n,k}\Vert_{L_{V}^{2}}}{\varpi^{\frac
{1}{2}}(\varpi-\gamma)^{\frac{1}{2}}}+M_{1}M^{\prime\prime}\sqrt{t}.
\]

Therefore, using (\ref{eq272})-(\ref{eq273}), we deduce%
\begin{equation}%
\begin{split}
\frac{d}{dt}|v^{n}(t)|^{2}  &  +\frac{m}{2}\Vert v^{n}(t)\Vert^{2}+\frac
{f_{0}}{2D}\Vert v^{n}(t)\Vert_{2p}^{2p}\leq2\alpha|\mathcal{O}|+C_{1,\delta
}(\theta_{t}\omega)(1+\Vert\phi\Vert_{2p}^{2p})\\[0.8ex]
&  +\frac{4M^{2}}{m}|y_{\delta}(\theta_{t}\omega)|^{2}\Vert\phi\Vert^{2}%
+\frac{4}{m\lambda_{1}}|h|^{2}+\frac{4}{m\lambda_{1}}|y_{\delta}(\theta
_{t}\omega)|^{2}|\phi|^{2}\\[0.8ex]
&  +\frac{C_{2,\delta}(\theta_{t}\omega)}{m}\left\Vert \phi\right\Vert
^{2}+\frac{4M_{1}\mu(t)\Vert\varphi_{v}^{n,k}\Vert_{L_{V}^{2}}^{2}}{\varpi
m(\varpi-\gamma)}+\frac{4(M_{1})^{2}(M^{\prime\prime})^{2}T}{m}.
\end{split}
\label{Eqvn}%
\end{equation}
The function $v$ satifies the same inequality but replacing $y_{\delta}$ by
$z_{\ast}$, $C_{i,\delta}$ by $C_{i}$, $i=1,2$ and $\varphi_{v}^{n,k}$ by
$\widetilde{\varphi}_{v}^{k}$. We define the functions
\[%
\begin{split}
J_{n}(t)  &  =|v^{n}(t)|^{2}-2\alpha|\mathcal{O}|t-\frac{4(M_{1}%
)^{2}(M^{\prime\prime})^{2}T}{m}t-\int_{0}^{t}C_{1,\delta}(\theta_{r}%
\omega)(1+\Vert\phi\Vert_{2p}^{2p})dr\\[0.8ex]
&  ~~-\int_{0}^{t}\left(  \frac{4M^{2}}{m}|y_{\delta}(\theta_{r}\omega
)|^{2}\Vert\phi\Vert^{2}+\frac{4}{m\lambda_{1}}|z_{\ast}(\theta_{r}%
\omega)|^{2}|\phi|^{2}\right)  dr\\[0.8ex]
&  ~~-\frac{4\left\vert h\right\vert ^{2}}{m\lambda_{1}}t-\frac{4M_{1}%
\Vert\varphi_{v}^{n,k}\Vert_{L_{V}^{2}}^{2}}{\varpi m(\varpi-\gamma)}\int%
_{0}^{t}\mu(r)dr-\int_{0}^{t}\frac{C_{2}(\theta_{r}\omega)}{m}\left\Vert
\phi\right\Vert ^{2}dr,
\end{split}
\]%
\[%
\begin{split}
J_{t}(t)  &  =|v(t)|^{2}-2\alpha|\mathcal{O}|t-\frac{4(M_{1})^{2}%
(M^{\prime\prime})^{2}T}{m}t-\int_{0}^{t}C_{1}(\theta_{r}\omega)(1+\Vert
\phi\Vert_{2p}^{2p})dr\\[0.8ex]
&  ~~-\int_{0}^{t}\left(  \frac{4M^{2}}{m}|z_{\ast}(\theta_{r}\omega
)|^{2}\Vert\phi\Vert^{2}+\frac{4}{m\lambda_{1}}|z_{\ast}(\theta_{r}%
\omega)|^{2}|\phi|^{2}\right)  dr\\[0.8ex]
&  ~~-\frac{4\left\vert h\right\vert ^{2}}{m\lambda_{1}}t-\frac{4M_{1}%
\Vert\widetilde{\varphi}_{v}^{k}\Vert_{L_{V}^{2}}^{2}}{\varpi m(\varpi
-\gamma)}\int_{0}^{t}\mu(r)dr-\int_{0}^{t}\frac{C_{2}(\theta_{r}\omega)}%
{m}\left\Vert \phi\right\Vert ^{2}dr.
\end{split}
\]
From the regularity of $v$ and all $v^{n}$, together with
(\ref{ConvergC1delta}), (\ref{ConvergC2delta}), (\ref{ConvergL2H}) and
(\ref{Eqvn}), it holds that these functions $J$ and $J_{n}$ are continuous and
non-increasing on $[-k,0]$, and
\[
J_{n}(s)\rightarrow J(s)~~\mbox{a.e.}~~s\in\lbrack-k,0]~\mbox{as}~n\rightarrow
\infty.
\]
Then the same argument as in Lemma \ref{lem4-5} implies that $v^{n}%
(t_{n})\rightarrow v(t_{0})$ strongly in $H$, and thus (\ref{eq16b}) follows.

\textbf{Statement 2.} The following inequality holds true:%
\begin{equation}
{\limsup_{n\rightarrow\infty}}\ \left\Vert v_{\delta_{n},0}-v_{0}\right\Vert
_{L_{V}^{2}}^{2}\leq Me^{-\gamma k}\limsup_{n\rightarrow\infty}\ \left(
\left\vert u^{n,k}-\tilde{u}^{k}\right\vert ^{2}+\left\Vert \varphi
^{n,k}-\tilde{\varphi}^{k}\right\Vert _{L_{V}^{2}}^{2}\right)  ,
\label{LimSup}%
\end{equation}
for any $k>0$, where $M$ is a positive constant and $\tilde{u}^{k}=\tilde
{v}^{k}+\phi z_{\ast}(\theta_{-k}\omega)$.

Define the functions $\overline{z}_{\delta_{n}}=\left(  \overline{v}%
_{\delta_{n}},\overline{\eta}_{\delta_{n}}\right)  =z_{\delta_{n}}-z$, where
$z(t)=(v(t),\eta^{t}(s))$ is the solution to problem (\ref{eq1-5}) with
initial condition $(v^{k},\mathcal{J}_{\omega,0}{\varphi^{k}})$ on $t=-k$.
Arguing as in Lemma \ref{lem4-5}, we have%
\[%
\begin{split}
\frac{d}{dt}\Vert\overline{z}_{\delta_{n}}\Vert_{\mathcal{H}}^{2}  &
+2(((\bar{\eta}_{\delta_{n}}^{t})^{\prime},\eta_{\delta_{n}}^{t}))_{\mu}%
\leq-2\int_{\mathcal{O}}(f(v_{\delta_{n}}+y_{\delta_{n}}(\theta_{t}%
\omega))-f(v+\phi z_{\ast}(\theta_{t}\omega)))(v_{\delta_{n}}-v)dx\\[0.8ex]
&  -2\int_{\mathcal{O}}\left(  a(l(v_{\delta_{n}})+l(\phi)y_{\delta_{n}%
}(\theta_{t}\omega))\nabla{v_{\delta_{n}}}-a(l(v)+l(\phi)z_{\ast}(\theta
_{t}\omega))\nabla v\right)  \cdot\nabla(v_{\delta_{n}}-v)dx\\[0.8ex]
&  +2\int_{\mathcal{O}}(a(l(v_{\delta_{n}})+l(\phi)y_{\delta_{n}}(\theta
_{t}\omega))y_{\delta_{n}}\left(  \theta_{t}\omega\right)  -a(l(v)+l(\phi
)z_{\ast}(\theta_{t}\omega))z_{\ast}(\theta_{t}\omega))\Delta\phi
(v_{\delta_{n}}-v)dx\\
&  +\int_{\mathcal{O}}\left(  y_{\delta_{n}}(\theta_{t}\omega)-z_{\ast}\left(
\theta_{t}\omega\right)  \right)  \phi\left(  v_{\delta_{n}}-v\right)
dx+{2}\int_{\mathcal{O}}\left(  z_{k,\delta_{n}}^{\phi}(t,\omega)-z_{k}^{\phi
}(t,\omega)\right)  \left(  v_{\delta_{n}}-v\right)  dx.
\end{split}
\]
Since $a$ is a locally Lipschitz function, $\phi\in V$, $z_{\ast}\left(
\theta_{t}\omega\right)  $ is uniformly bounded for any $\omega\in\Omega$ on
$[-k,0]$, making use of (\ref{eq5-5}), we deduce that
\begin{align*}
\{l(v_{\delta_{n}}(t))+l(\phi)y_{n}(\theta_{t}\omega)\}_{t\in\lbrack
-k,0],0<\delta\leq\overline{\sigma}}  &  \subset\lbrack-R,R],\\
\{l(v(t))+l(\phi)z_{\ast}(\theta_{t}\omega)\}_{t\in\lbrack-k,0],0<\delta
\leq\delta_{0}}  &  \subset\lbrack-R,R],
\end{align*}
for some $R,\delta_{0}>0,$ and
\[
\left\vert a(l(v_{\delta_{n}})+l(\phi)y_{\delta_{n}}(\theta_{t}\omega
))-a(l(v)+l(\phi)z_{\ast}(\theta_{t}\omega))\right\vert \leq L_{a}%
(R)\left\vert l\right\vert (\left\vert v_{\delta_{n}}-v\right\vert +\left\vert
y_{\delta_{n}}(\theta_{t}\omega)-z_{\ast}(\theta_{t}\omega)\right\vert
{|\phi|}).
\]
Hence, by \eqref{eq1-2} and the Young inequality, we infer that
\[%
\begin{split}
&  \quad-2\int_{\mathcal{O}}\left(  a(l(v_{\delta_{n}})+l(\phi)y_{\delta_{n}%
}(\theta_{t}\omega))\nabla v_{\delta_{n}}-a(l(v)+l(\phi)z_{\ast}(\theta
_{t}\omega))\nabla v\right)  \cdot\nabla(v_{\delta_{n}}-v)dx\\[0.8ex]
&  \leq-2m\Vert v_{\delta_{n}}-v\Vert^{2}+2L_{a}(R)|l|\left(  |v^{n}%
-v|+\left\vert y_{\delta_{n}}(\theta_{t}\omega)-z_{\ast}(\theta_{t}%
\omega)\right\vert {|\phi|}\right)  \Vert v\Vert\Vert v_{\delta_{n}}%
-v\Vert\\[0.8ex]
&  \leq(\alpha-2m)\Vert v_{\delta_{n}}-v\Vert^{2}+\frac{2L_{a}^{2}(R)|l|^{2}%
}{\alpha}\left(  |v_{\delta_{n}}-v|^{2}+\left\vert y_{\delta_{n}}(\theta
_{t}\omega)-z_{\ast}(\theta_{t}\omega)\right\vert ^{2}{|\phi|^{2}}\right)
\Vert v\Vert^{2},
\end{split}
\]
where $\alpha\leq(m\lambda_{1}-\gamma)/\lambda_{1}$. By the above estimates,
we deduce that%
\[%
\begin{split}
&  ~\quad\frac{d}{dt}\Vert\overline{z}_{\delta_{n}}\Vert_{\mathcal{H}}%
^{2}+{\gamma\Vert\bar{z}_{\delta_{n}}\Vert_{\mathcal{H}}^{2}+m\Vert
v_{\delta_{n}}-v\Vert^{2}}\\[0.8ex]
&  \leq\frac{d}{dt}\Vert\overline{z}_{\delta_{n}}\Vert_{\mathcal{H}}%
^{2}+(2m-\alpha)\Vert v_{\delta_{n}}-v\Vert^{2}+\varpi\int_{0}^{\infty}%
\mu(s)|\nabla\bar{\eta}_{\delta_{n}}^{t}(s)|^{2}ds\\[0.8ex]
&  \leq\frac{2L_{a}^{2}(R)|l|^{2}}{\alpha}\left(  |v_{\delta_{n}}%
-v|^{2}+\left\vert y_{\delta_{n}}(\theta_{t}\omega)-z_{\ast}(\theta_{t}%
\omega)\right\vert ^{2}{|\phi|^{2}}\right)  \Vert v\Vert^{2}\\[0.8ex]
&  -2\int_{\mathcal{O}}(f(v_{\delta_{n}}+\phi y_{\delta_{n}}(\theta_{t}%
\omega))-f(v+\phi z_{\ast}(\theta_{t}\omega)))(v_{\delta_{n}}-v)dx\\
&  +2\int_{\mathcal{O}}(a(l(v_{\delta_{n}})+l(\phi)y_{\delta_{n}}(\theta
_{t}\omega))y_{\delta_{n}}\left(  \theta_{t}\omega\right)  -a(l(v)+l(\phi
)z_{\ast}(\theta_{t}\omega))z_{\ast}(\theta_{t}\omega))\Delta\phi
(v_{\delta_{n}}-v)dx\\
&  +2\int_{\mathcal{O}}\left(  y_{\delta_{n}}(\theta_{t}\omega)-z_{\ast
}\left(  \theta_{t}\omega\right)  \right)  \phi\left(  v_{\delta_{n}%
}-v\right)  dx+{2}\int_{\mathcal{O}}\left(  z_{k,\delta_{n}}^{\phi}%
-z_{k}^{\phi}\right)  \left(  v_{\delta_{n}}-v\right)  dx,
\end{split}
\]
where we have used that $0<\gamma\leq\min\{(m-\alpha)\lambda_{1},\delta\}$ by
the choice of $\alpha$. Multiplying by $e^{\gamma t}$ on both sides of the
above inequality and integrating over $(-k,0)$, we obtain%
\[%
\begin{split}
&  ~\quad\Vert\overline{z}_{\delta_{n}}(t)\Vert_{\mathcal{H}}^{2}+m\int%
_{-k}^{0}e^{\gamma s}\Vert v_{\delta_{n}}(s)-v(s)\Vert^{2}ds\\
&  \leq e^{-\gamma k}\Vert\overline{z}_{\delta_{n}}(-k)\Vert_{\mathcal{H}}%
^{2}+\frac{2L_{a}^{2}(R)|l|^{2}}{\alpha}\int_{-k}^{0}e^{\gamma s}%
(|v_{\delta_{n}}-v|^{2}+\left\vert y_{\delta_{n}}(\theta_{s}\omega)-z_{\ast
}(\theta_{s}\omega)\right\vert )\Vert v\Vert^{2}ds\\
&  ~~-2\int_{-k}^{0}e^{\gamma s}\int_{\mathcal{O}}(f(v_{\delta_{n}}+\phi
y_{\delta_{n}}(\theta_{s}\omega))-f(v+\phi z_{\ast}(\theta_{s}\omega
)))(v_{\delta_{n}}-v)dxds\\
&  ~~+2\int_{-k}^{0}e^{\gamma s}\int_{\mathcal{O}}(a(l(v_{\delta_{n}}%
)+l(\phi)y_{\delta_{n}}(\theta_{s}\omega))y_{\delta_{n}}\left(  \theta
_{s}\omega\right)  -a(l(v)+l(\phi)z_{\ast}(\theta_{s}\omega))z_{\ast}%
(\theta_{s}\omega))\Delta\phi(v_{\delta_{n}}-v)dxds\\
&  ~~+{2}\int_{-k}^{0}e^{\gamma s}\int_{\mathcal{O}}\left(  y_{\delta_{n}%
}(\theta_{s}\omega)-z_{\ast}\left(  \theta_{s}\omega\right)  \right)
\phi\left(  v_{\delta_{n}}-v\right)  dxds\\
&  ~~+{2}\int_{-k}^{0}e^{\gamma s}\int_{\mathcal{O}}\left(  z_{k,\delta_{n}%
}^{\phi}(s,\omega)-z_{k}^{\phi}(s,\omega)\right)  \left(  v_{\delta_{n}%
}-v\right)  dxds.
\end{split}
\]
By a similar argument as in Lemma \ref{lem4-5}, we obtain that the second and
fourth terms of the right-hand side of the above inequality converge to zero.
Also, by (\ref{eq5-5}) and (\ref{ConvergL2H}), it is easy to see that the
fifth term goes to zero as well. On the other hand, by (\ref{eq5-5}%
)-(\ref{eq5-6}), we deduce easily that $z_{k,\delta_{n}}^{\phi}(s,\omega
)\rightarrow z_{k}^{\phi}(s,\omega)$ uniformly on $[-k,0]$. Hence,
(\ref{ConvergL2H}) implies that the last term of the above inequality also
converges to zero.

It remains to analyze the third term. On the one hand, as in the proof of
Lemma \ref{lem4-5}, we obtain that $f(v^{n}+\phi y_{\delta_{n}}(\theta
_{t}\omega))\rightarrow f(v+\phi z_{\ast}(\theta_{t}\omega))$ weakly in
$L^{q}\left(  -k,0;L^{q}(\mathcal{O})\right)  $. Thus,%
\[
\int_{-k}^{0}e^{\gamma s}\int_{\mathcal{O}}(f(v_{\delta_{n}}+\phi
y_{\delta_{n}}(\theta_{s}\omega))-f(v+\phi z_{\ast}(\theta_{s}\omega
)))vdxds\rightarrow0.
\]
On the other hand, by the same calculations as in Lemma \ref{lem4-5}, we infer
that there are positive constants $\kappa_{i}\left(  \omega,\phi,k\right)  $,
$i=1,2$, such that%
\[
f(v^{n}(s,x,\omega)+\phi y_{n}(\theta_{s}\omega))v^{n}(s,x,\omega)\geq
-\kappa_{1}+\kappa_{2}\left\vert v^{n}(s,,x,\omega)+\phi y_{n}(\theta
_{s}\omega)\right\vert ^{2p}.
\]
Then the Fatou-Lebesgue lemma implies that%
\[%
\begin{split}
&  \qquad\limsup_{n\rightarrow\infty}\left(  -2\int_{-k}^{0}e^{\gamma s}%
\int_{\mathcal{O}}f(v_{\delta_{n}}+\phi y_{\delta_{n}}(\theta_{s}%
\omega))v_{\delta_{n}}dxds\right) \\[0.8ex]
&  \leq-2\int_{0}^{t}e^{\gamma s}\int_{\mathcal{O}}\liminf_{n\rightarrow
\infty}f(v_{\delta_{n}}+\phi y_{\delta_{n}}(\theta_{s}\omega))v_{\delta_{n}%
}dxds\\[0.8ex]
&  =-2\int_{0}^{t}e^{\gamma s}\int_{\mathcal{O}}f(v+\phi z_{\ast}(\theta
_{s}\omega))vdxds.
\end{split}
\]
Since $f(v+\phi z_{\ast}(\theta_{\cdot}\omega))\in L^{q}(-k,0;L^{q}%
(\mathcal{O}))$ and $v_{\delta_{n}}\rightarrow v$ weakly in $L^{2p}%
(-k,0;L^{2p}(\mathcal{O}))$, \eqref{eq4-17} holds true. Hence,
\[
\limsup_{n\rightarrow\infty}\left(  -2\int_{-k}^{0}e^{\gamma s}\int%
_{\mathcal{O}}(f(v_{\delta_{n}}+\phi y_{\delta_{n}}(\theta_{s}\omega
))-f(v+\phi z_{\ast}(\theta_{s}\omega)))(v_{\delta_{n}}-v)dxds\right)  {\leq
0}.
\]
By Lemma \ref{lem3-1}, the proof of Corollary \ref{cor3-2}, (\ref{eq5-5}),
Remark \ref{ConvergYdelta} and the above convergences, we arrive at%
\begin{align*}
&  \limsup_{n\rightarrow\infty}\ \int_{-k}^{0}e^{\gamma s}\Vert v_{\delta_{n}%
}(s)-v(s)\Vert^{2}ds\\
&  \leq\frac{1}{m}e^{-\gamma k}\limsup_{n\rightarrow\infty}\Vert\overline
{z}_{\delta_{n}}(-k)\Vert_{\mathcal{H}}^{2}\\
&  \leq\frac{1}{m}e^{-\gamma k}\limsup_{n\rightarrow\infty}\left(  \left\vert
v^{n,k}-\tilde{v}^{k}\right\vert ^{2}+\int_{0}^{\infty}\mu\left(  s\right)
\left\Vert \mathcal{J}_{\omega,0}^{\delta_{n}}{\varphi^{n,k}}-\mathcal{J}%
_{\omega,0}\widetilde{\varphi}^{k}\right\Vert ^{2}ds\right) \\
&  \leq\frac{1}{m}e^{-\gamma k}\limsup_{n\rightarrow\infty}\left(  \left\vert
v^{n,k}-\tilde{v}^{k}\right\vert ^{2}+2\int_{0}^{\infty}\mu\left(  s\right)
\left\Vert \mathcal{J}\left(  \varphi^{n,k}-\widetilde{\varphi}^{k}\right)
\right\Vert ^{2}ds\right. \\
&  \left.  +2\int_{0}^{\infty}\mu\left(  s\right)  \left(  \int_{-s}%
^{0}\left\vert y_{\delta_{n}}(\theta_{r}\omega)-z_{\ast}(\theta_{r}%
\omega)\right\vert \left\Vert \phi\right\Vert ds\right)  ^{2}\right) \\
&  \leq\frac{1}{m}e^{-\gamma k}\limsup_{n\rightarrow\infty}\left(  2\left\vert
u^{n,k}-\tilde{u}^{k}\right\vert ^{2}+2K_{\mu}\int_{-\infty}^{0}e^{\gamma
s}\left\Vert \varphi^{n,k}(s)-\widetilde{\varphi}^{k}(s)\right\Vert
^{2}ds\right. \\
&  \left.  +2\left\vert y_{\delta_{n}}(\theta_{-k}\omega)-z_{\ast}(\theta
_{-k}\omega)\right\vert ^{2}\left\vert \phi\right\vert ^{2}+2K_{\mu}%
\int_{-\infty}^{0}e^{\gamma s}\left\Vert \phi\right\Vert ^{2}\left\Vert
y_{\delta_{n}}(\theta_{s}\omega)-z_{\ast}(\theta_{s}\omega)\right\Vert
^{2}ds\right) \\
&  =\frac{1}{m}e^{-\gamma k}\limsup_{n\rightarrow\infty}\left(  2\left\vert
u^{n,k}-\tilde{u}^{k}\right\vert ^{2}+2K_{\mu}\int_{-\infty}^{0}e^{\gamma
s}\left\Vert \varphi^{n,k}(s)-\widetilde{\varphi}^{k}(s)\right\Vert
^{2}ds\right)  .
\end{align*}
Then, the same arguments as in Lemma \ref{lem4-5} imply (\ref{LimSup}) immediately.

\textbf{Statement 3.} There is a subsequence such that
\begin{equation}
(u^{n},\varphi^{n})\rightarrow(u_{0},\varphi_{0}),\quad\text{strongly in
}X~\text{as}~n\rightarrow\infty. \label{eq5Strong}%
\end{equation}
By Lemma \ref{BBarraTemp}, we know that the family $\overline{B}(\omega)$ is
tempered in $X$. By (\ref{LimSup}), $(u^{n,k},\varphi^{n,k}),\ (\tilde{u}%
^{k},\tilde{\varphi}^{k})\in\overline{B}(\theta_{-k}\omega)$ and choosing some
$0<c<\gamma$, there is a constant $R>0$ such that%
\begin{align*}
\limsup_{n\rightarrow\infty}\ \left\Vert v_{\delta_{n},0}-v_{0}\right\Vert
_{L_{V}^{2}}^{2}  &  \leq Me^{-(\gamma-c)k}e^{-ck}\limsup_{n\rightarrow\infty
}\ \left(  \left\vert u^{n,k}-\tilde{u}^{k}\right\vert ^{2}+\left\Vert
\varphi^{n,k}-\tilde{\varphi}^{k}\right\Vert _{L_{V}^{2}}^{2}\right)  \text{
}\\
&  \leq Re^{-(\gamma-c)k},
\end{align*}
for all $k\geq1.$ Further, for every $d\in\mathbb{N}$, there is $k_{0}\left(
d\right)  $ such that,%
\[
\limsup_{n\rightarrow\infty}\ \left\Vert v_{\delta_{n},0}-v_{0}\right\Vert
_{L_{V}^{2}}^{2}\leq\frac{1}{d},\qquad\forall k\geq k_{0}.
\]
Taking $d\rightarrow\infty$ and using a diagonal argument, we deduce the
existence of a subsequence $\{v_{\delta_{n_{d}},0}\}$ such that $v_{\delta
_{n_{d}},0}\rightarrow v_{0}$ in $L_{V}^{2}.$ Together with (\ref{eq16b}) and
$\phi y_{\delta_{n}}(\theta_{\cdot}\omega)\rightarrow\phi z_{\ast}%
(\theta_{\cdot}\omega)$ in $X$, it shows that (\ref{eq5Strong}) is true. The
proof of this lemma is complete. $\Box$

\begin{theorem}
Assume that \eqref{eq1-2}, \eqref{eq2-1} and $(h_{1})$-$(h_{2})$ hold. Let
$\phi\in V\cap H^{2}(\mathcal{O})\cap L^{2p}(\mathcal{O})$ be such that
$\Delta\phi\in L^{2p}(\mathcal{O})$, $h\in H$ and $a$ be a locally Lipschitz
function. Then, for all $\omega\in\Omega$,
\[
\lim_{\delta\rightarrow0}dist_{X}(\mathcal{A}_{\delta}(\omega),\mathcal{A}%
(\omega))=0.
\]

\end{theorem}

\textbf{Proof. }By \eqref{rhoconvergence}, we have
\[
\lim_{\delta\rightarrow0}\Vert B_{\delta}(\omega)\Vert_{X}=\Vert
B(\omega)\Vert_{X},\qquad\mbox{for all}~~\omega\in\Omega,
\]
where for a set $S\subset X$, we denote $\left\Vert S\right\Vert _{X}%
=\sup_{y\in S}\left\Vert y\right\Vert _{X}$. This, together with Corollary
\ref{cor6-3} and Lemma \ref{lem6-4}, finishes the proof of this theorem by
\cite[Theorem 3.1]{W2}. $\Box$

\section{Appendix}

Let us consider equation (\ref{eq1-1}) in the deterministic case, that is,%
\begin{equation}%
\begin{cases}
\dfrac{\partial u}{\partial t}-a(l(u))\Delta u-\displaystyle\int_{-\infty}%
^{t}k(t-s)\Delta u(x,s)ds+f(u)=h,\\
u(x,t)=0,\\
u(x,0)=u_{0}(x),\\
u(x,t)=\phi(x,t),
\end{cases}
\begin{aligned} &\mbox{in}~~\mathcal{O}\times (\tau,\infty),\\ &\mbox{on}~\partial\mathcal{O}\times(\tau,\infty),\\ &\mbox{in}~~\mathcal{O}\\ &\mbox{in}~~\mathcal{O}\times (-\infty,0],\\ \end{aligned} \label{EqDet}%
\end{equation}
where $\mathcal{O}\subset\mathbb{R}^{N}$ is a bounded domain with regular
boundary. We assume that (\ref{eq1-2}), (\ref{eq2-1}) and $(h_{1})$-$(h_{2})$
hold. Also, let $h\in H$ and $a$ be a locally Lipschitz function.

As before, we consider Dadermos' tranform%
\[
\eta^{t}(s,x)=\int_{t-s}^{t}u(r,x)dr,~~~\text{ for }s\geq0,
\]
which gives rise to the system,%
\begin{equation}%
\begin{cases}
\dfrac{\partial u}{\partial t}-a(l(u))\Delta u-\displaystyle\int_{0}^{\infty
}\mu(s)\Delta\eta^{t}(s)ds+f(u)=h,\\
\dfrac{\partial}{\partial t}\eta^{t}(s)=u-\dfrac{\partial}{\partial s}\eta
^{t}(s),\\
u(x,t)=\eta^{t}(x,s)=0,\\
u(x,0)=u_{0}(x),\\
\eta^{0}(x,s)=\eta_{0}(x,s),
\end{cases}
\begin{aligned} &\mbox{in}~~\mathcal{O}\times (\tau,\infty),\\ &\mbox{in}~~\mathcal{O}\times (\tau,\infty)\times\mathbb{R}^+,\\ &\mbox{on}~~\partial \mathcal{O}\times \mathbb{R}, s>0,\\ &\mbox{in}~~ \mathcal{O},\\ &\mbox{in}~~\mathcal{O}\times\mathbb{R}^+, \end{aligned} \label{EqDetDaf}%
\end{equation}
where%
\[
\eta^{0}(x,s)=\int_{-s}^{0}u(x,r)dr=\int_{-s}^{0}\phi(x,r)dr:=\eta_{0}(x,s).
\]

For any $\left(  u_{0},\phi\right)  \in X=H\times L_{V}^{2}$, there exists a
unique weak solution $z\left(  \text{\textperiodcentered}\right)  =\left(
u\left(  \text{\textperiodcentered}\right)  ,\eta^{\text{\textperiodcentered}%
}\right)  $ to problem (\ref{EqDetDaf}) \cite[Theorem 3.4]{X2}. This is also a
particular case of the result in Theorem \ref{thm3-3} with $\phi=0$.

\begin{lemma}
Let $\{\left(  u_{0}^{n},\phi^{n}\right)  \}$ be a sequence such that $\left(
u_{0}^{n},\phi^{n}\right)  \rightarrow\left(  u_{0},\phi\right)  $ weakly in
$X$. Then
\begin{equation}
u^{n}\rightarrow u\text{ in }C([r,T],H),~~\text{ for all }0<r<T,
\label{ConvergH}%
\end{equation}
where $z^{n}\left(  \text{\textperiodcentered}\right)  =\left(  u^{n}\left(
\text{\textperiodcentered}\right)  ,\eta_{n}^{\text{\textperiodcentered}%
}\right)  $, $z\left(  \text{\textperiodcentered}\right)  =\left(  u\left(
\text{\textperiodcentered}\right)  ,\eta^{\text{\textperiodcentered}}\right)
$ are the solutions to problem (\ref{EqDetDaf}) corresponding to $\left(
u_{0}^{n},\phi^{n}\right)  $ and $\left(  u_{0},\phi\right)  $, respectively.
\end{lemma}

\textbf{Proof. }It is known (see the proof of Lemma 3.9 in \cite{X2}) that,%
\begin{align*}
u^{n}  &  \rightarrow u~~\text{ weak-}\ast\text{ in }~~L^{\infty
}(0,T;H);\\[0.01in]
u^{n}  &  \rightarrow u~~\text{ weakly in }~~L^{2}(0,T;V);\\[1ex]
u^{n}  &  \rightarrow u~~\text{ strongly in }~~L^{2}(0,T;H);\\
\frac{du^{n}}{dt}  &  \rightarrow\frac{du}{dt}~~\text{ weakly in }%
~~L^{2}(0,T;V^{\ast})+L^{q}\left(  0,T;L^{q}\left(  \Omega\right)  \right)
;\\
u^{n}(t)  &  \rightarrow u(t)\text{ in }H~~\text{ for a.a. }t\in\left(
0,T\right)  .
\end{align*}
By the same arguments in the proof of Lemma \ref{lem4-5}, we obtain
\[
u_{n}\rightarrow u~~\text{ in }C([-k,0],V^{\ast}+L^{q}(\mathcal{O})).
\]
Then, if $t_{n}\rightarrow t_{0}$, $t_{n}\in\lbrack0,T],\ t_{0}\in(0,T]$, we
infer%
\[
v_{n}(t_{n})\rightarrow v(t_{0})\text{ weakly in }H,
\]
and
\[
\left\vert v(t_{0})\right\vert \leq\liminf_{n\rightarrow\infty}\ \left\vert
v_{n}(t_{n}))\right\vert .
\]
We need to prove that $v^{n}(t_{n})\rightarrow v(t_{0})$ strongly in $H$. By
Corollary \ref{CorExistSol}, we know $u^{n}$ are weak solutions of the
equation%
\[
u_{t}^{n}-a(l(u))\Delta u^{n}-\int_{-\infty}^{t}k(t-s)\Delta u^{n}%
ds+f(u^{n})=h.
\]
Then,%
\[
\frac{1}{2}\frac{d}{dt}|u^{n}(t)|^{2}+m\Vert u^{n}\Vert^{2}+(f(u^{n}%
),u^{n})\leq\left(  \int_{-\infty}^{t}k(t-s)\Delta u^{n}(s)ds,u^{n}(t)\right)
+(h,u^{n}(t)).
\]
By (\ref{fcondition}) and the Young inequality, we obtain
\[
\frac{d}{dt}|u^{n}|^{2}+m\Vert u^{n}\Vert^{2}+f_{0}\Vert u^{n}\Vert_{2p}%
^{2p}\leq2\alpha|\mathcal{O}|+\frac{1}{m\lambda_{1}}|h|^{2}+2\int_{-\infty
}^{t}k(t-s)\Vert u^{n}(s)\Vert ds\Vert u^{n}(t)\Vert.
\]
By the arguments in Lemma \ref{lem4-5}, we have%
\[
\int_{-\infty}^{t}k(t-s)\Vert u^{n}(s)\Vert ds\leq\frac{M_{1}^{\frac{1}{2}}%
\mu^{\frac{1}{2}}(t)\Vert\phi^{n}\Vert_{L_{V}^{2}}}{\varpi^{\frac{1}{2}%
}(\varpi-\gamma)^{\frac{1}{2}}}+M_{1}M^{\prime\prime}\sqrt{t}.
\]
Thus,%
\begin{align*}
&  \frac{d}{dt}|u^{n}(t)|^{2}+m\Vert u^{n}(t)\Vert^{2}+f_{0}\Vert
u^{n}(t)\Vert_{2p}^{2p}\\
&  \leq2\alpha|\mathcal{O}|+\frac{1}{m\lambda_{1}}|h|^{2}+2\left(  \frac
{M_{1}^{\frac{1}{2}}\mu^{\frac{1}{2}}(t)\Vert\phi^{n}\Vert_{L}{}_{V}^{2}%
}{\varpi^{\frac{1}{2}}(\varpi-\gamma)^{\frac{1}{2}}}+M_{1}M^{\prime}%
\prime\sqrt{T}\right)  \Vert u^{n}(t)\Vert.
\end{align*}
Therefore,%
\begin{align*}
&  \frac{d}{dt}|u^{n}(t)|^{2}+\frac{m}{2}\Vert u^{n}(t)\Vert^{2}+f_{0}\Vert
u^{n}(t)\Vert_{2p}^{2p}\\
&  \leq2\alpha|\mathcal{O}|+\frac{1}{m\lambda_{1}}|h|^{2}+\frac{4M_{1}%
\mu(t)\Vert\phi^{n}\Vert_{L_{V}^{2}}^{2}}{\varpi m(\varpi-\gamma)}%
+\frac{4(M_{1})^{2}(M^{\prime\prime})^{2}T}{m}.
\end{align*}
The function $u$ satifies the same inequality but replacing $\phi^{n}$ by
$\phi$. We define the functions
\[
J_{n}(t)=|u^{n}(t)|^{2}-2\alpha|\mathcal{O}|t-\frac{4(M_{1})^{2}%
(M^{\prime\prime})^{2}T}{m}t-\frac{\left\vert h\right\vert ^{2}}{m\lambda_{1}%
}t-\frac{4M_{1}\Vert\phi^{n}\Vert_{L_{V}^{2}}^{2}}{\varpi m(\varpi-\gamma
)}\int_{0}^{t}\mu(r)dr,
\]%
\[
J(t)=|u(t)|^{2}-2\alpha|\mathcal{O}|t-\frac{4(M_{1})^{2}(M^{\prime\prime}%
)^{2}T}{m}t-\frac{\left\vert h\right\vert ^{2}}{m\lambda_{1}}t-\frac
{4M_{1}\Vert\phi\Vert_{L_{V}^{2}}^{2}}{\varpi m(\varpi-\gamma)}\int_{0}^{t}%
\mu(r)dr.
\]
These functions are continuous and non-increasing on $[0,T]$, and
\[
J_{n}(s)\rightarrow J(s)~~\text{ for a.a. }s\in\lbrack0,T]~\text{as}%
~n\rightarrow\infty.
\]
Then the same argument as in Lemma \ref{lem4-5} ensures that $v^{n}%
(t_{n})\rightarrow v(t_{0})$ strongly in $H$, and therefore (\ref{ConvergH}) follows.

\bigskip

\begin{remark}
The convergence (\ref{ConvergH}) was stated in Lemma 3.9 from \cite{X2}.
However, the proof of this result is incorrect there and we provided here a
correct one.
\end{remark}

\section*{ Conflict of interest}

The authors declare that they have no conflict of interest.

\section*{ Data availability}

The authors declare that there are no data associated to this paper.

\section*{Acknowledgements}

The research has been partially supported by the National Natural Science
Foundation of China (No. 12301234), the Spanish Ministerio de Ciencia e
Innovaci\'{o}n (MCI), Agencia Estatal de Investigaci\'{o}n (AEI) and Fondo
Europeo de Desarrollo Regional (FEDER) under the projects PID2021-122991NB-C21
and\ PID2019-108654GB-I00, and by the Generalitat Valenciana, project PROMETEO/2021/063.


\bigskip

\begin{thebibliography}{99}                                                                                               %


\bibitem {delfour}A. Bensoussan, G. Da Prato, M. C. Delfour, S. K. Mitter,
Representation and Control of Infinite-Dimensional Systems, Vol. 1. Systems \&
Control: Foundations \& Applications. Birkh\"{a}user Boston, Inc., Boston, MA, 1992.

\bibitem {B1}Z. Brze\'zniak, M. Capi\'nski, F. Flandoli, A convergence result
for stochastic partial differential equations, Stochastics, 24 (1988), 423-445.

\bibitem {CR}T. Caraballo, J. Real, Attractors for 2D-Navier-Stokes models
with delays, J. Differential Equations, 205 (2004), 271-297.

\bibitem {C2}T. Caraballo, J. Real, I. D. Chueshov, Pullback attractors for
stochastic equations in materials with memory, Discrete Contin. Dyn. Syst.
Ser. B, 9 (2008), 525-539.

\bibitem {C10}T. Caraballo, M. J. Garrido-Atienza, B. Schmalfu\ss , J. Valero,
Global attractor for a non-autonomous integro-differential equation in
materials with memory, Nonlinear Anal., 73 (2010), 183-201.

\bibitem {C4}T. Caraballo, M. J. Garrido-Atienza, B. Schmalfu\ss , J. Valero,
Attractors for a random evolution equation with infinite memory: theoretical
results, Discrete Contin. Dyn. Syst. Ser. B, 22 (2017), 1779-1800.

\bibitem {CW}P. Y. Chen, R. H. Wang, X. P. Zhang, Long-time dynamics of
fractional nonclassical diffusion equations with nonlinear colored noise and
delay on unbounded domains, Bull. Sci. Math., 173 (2021), 103071, 52pp.

\bibitem {CZ}Z. Chen, L. Y. Li, D. D. Yang, Asymptotic behavior of random
coupled Ginzburg-Landau equation driven by colored noise on unbounded domains,
Adv. Difference Equ., (2021), No. 291, 9pp.

\bibitem {C3}I. D. Chueshov, M. Scheutzow, Inertial manifolds and forms for
stochastically perturbed retarded semilinear parabolic equations, J. Dynam.
Differential Equations, 13 (2001), 355-380.

\bibitem {C1}M. Conti, V. Pata, M. Squassina, Singular limit of differential
systems with memory, Indiana Univ. Math. J., 1 (2006), 169-215.

\bibitem {Cui1}J. Cui, S. Liu, H. Zhou, Stochastic Wasserstein Hamiltonian
Flows, J. Dynamics Differential Equations, 2023, https://doi.org/10.1007/s10884-023-10264-4.

\bibitem {Cui2}J. Cui, S. Liu, H. Zhou, Wasserstein Hamiltonian Flow With
Common Noise On Graph, SIAM Journal on Applied Mathematics, 83 (2023), 484-509.

\bibitem {Gajewski}H. Gajewski, K. Gr\"{o}ger and K. Zacharias, Nichtlineare
Operatorgleichungen und Operatordifferentialgleichungen, Academi-Verlag,
Berlin, 1974.

\bibitem {G3}S. Gatii, M. Grasselli, V. Pata, Exponential attractors for a
phased-filed model with memory and quadratic nonlinearities, Indiana Univ.
Math. J., 53 (2004), 719-753.

\bibitem {G1}C. Giorgi, V. Pata, A. Marzocchi, Asymptotic behavior of a
semilinear problem in heat conduction with memory, Nonlinear Differ. Equ.
Appl., 5 (1998), 333-354.

\bibitem {G2}M. Grasselli, V. Pata, Uniform Attractors of Nonautonomous
Dynamical Systems with Memory, Evolution Equations, Semigroups and Functional
Analysis, Progr. Nonlinear Differential Equations Appl., 50, Birkh\"{a}user,
Basel, 2002.

\bibitem {G6}A. H. Gu, B. X. Wang, Asymptotic behavior of random
Fitzhugh-Nagumo sysmtes driven by colored noise, Discrete Contin. Dyn. Syst.
Ser. B, 23 (2018), 1689-1720.

\bibitem {G5}A. H. Gu, K. N. Lu, B. X. Wang, Asymptotic behavior of random
Navier-Stokes equations driven by Wong-Zakai approximations, Discrete Contin.
Dyn. Syst., 39 (2019), 185-218.

\bibitem {G4}A. H. Gu, B. L. Guo, B. X. Wang, Long term behavior of random
Navier-Stokes equations driven by colored noise, Discrete Contin. Dyn. Syst.
Ser. B, 25 (2020), 2495-2532.

\bibitem {L3}D. S. Li, X. H. Wang, J. Y. L. Zhao, Limiting dynamical behavior
of random fractional Fitzhugh-Nagumo systems driven by a Wong-Zakai
approximation process, Commun. Pure Appl. Anal., 19 (2020), 2751-2776.

\bibitem {L1}J. L. Lions, Quelques M\'{e}thodes de R\'{e}solution des
Probl\`{e}mes aux Limites non Lin\'{e}aires, Gauthier-Villar, Paris, 1969.

\bibitem {L2}K. N. Lu, B. X. Wang, Wong-Zakai approximations and long term
behavior of stochastic partial differential equations, J. Dynam. Differential
Equations, 31 (2019), 1341-1371.

\bibitem {N1}A. Nowak, A Wong-Zakai type theorem for stochastic systems of
Burgers equations, Panam. Math. J., 16 (2006), 1-25.

\bibitem {P1}V. Pata, A. Zucchi, Attractors for a damped hyperbolic equation
with linear memory, Adv. Math. Sci. Appl., 11 (2001), 505-529.

\bibitem {R1}J. Robinson, Infinite Dimensional Dynamical Systems, Cambridge
University Press, Cambridge, 2001.

\bibitem {T1}R. Temam, Navier-Stokes Equations, Theory and Numerical Analysis,
North-Holland, Amsterdam, 1979.

\bibitem {T3}G. Tessitore, J. Zabczyk, Wong-Zakai approximations of stochastic
evolution equations, J. Evol. Equ., 6 (2006), 621-655.

\bibitem {T2}K. Twardowska, An approximation theorem of Wong-Zakai type for
nonlinear stochastic partial differential equations, Stoch. Anal. Appl., 13
(1995), 601-626.

\bibitem {W1}B. X. Wang, Sufficient and necessary criteria for existence of
pullback attractors for non-compact random dynamical systems, J. Differential
Equations, 253 (2012), 1544-1583.

\bibitem {W2}B. X. Wang, Existence and upper semicontinuity of attractors for
stochastic equations with deterministic non-autonomous terms, Stoch. Dyn., 14
(2014), 1-31.

\bibitem {W3}Y. J. Wang, J. Y. Wang, Pullback attractors for multi-valued
non-compact random dynamical systems generated by reaction-diffusion equations
on an unbounded domain, J. Differential Equations, 259 (2015), 728-776.

\bibitem {W4}R. H. Wang, Y. R. Li, B. X. Wang, Random dynamics of fractional
nonclassical diffusion equations driven by colored noise, Discrete Contin.
Dyn. Syst., 39 (2019), 4091-4126.

\bibitem {WZ}E. Wong, M. Zakai, On the convergence of ordinary integrals to
stochastic integrals, Ann. Math. Statist., 36 (1965), 1560-1564.

\bibitem {X0}J. H. Xu, Z. Zhang, T. Caraballo, Non-autonomous nonlocal partial
differential equations with delay and memory, J. Differential Equations, 270
(2021), 505-546.

\bibitem {X1}J. H. Xu, T. Caraballo, Dynamics of stochastic nonlocal partial
differential equations, Eur. Phys. J. Plus, (2021), 136:849.

\bibitem {X4}J. H. Xu, T. Caraballo, Long time behavior of stochastic nonlocal
partial differential equations and Wong-Zakai approximations, SIAM J. Math.
Anal., 54 (2022), 2792-2844.

\bibitem {X2}J. H. Xu, T. Caraballo, J. Valero, Asymptotic behavior of a
semilinear problem in heat conduction with long time memory and non-local
diffusion, J. Differential Equations, 327 (2022), 418-447.

\bibitem {X3}J. H. Xu, T. Caraballo, J. Valero, Asymptotic behavior of
nonlocal partial differential equations with long time memory, Discrete
Contin. Dyn. Syst. Ser. S, 15 (2022), 3059-3078.

\bibitem {Y1}L. Yang, Y. J. Wang, Attractors for 2D quasi-geostrophic
equations with and without colored noise in $W^{2\alpha^{-},p}(\mathbb{R}%
^{2})$, Stoch. Dyn., 21 (2021), Paper No. 2150017, 32pp.
\end{thebibliography}
\end{document}